\documentclass[reqno]{amsart}%
\usepackage{palatino, mathpazo}
\usepackage{amsfonts}
\usepackage{amsmath}
\usepackage{tikz,graphicx,capt-of}
\usepackage{tikz,capt-of}
\usepackage{tikz,amssymb,float}
\usepackage{amssymb,latexsym,xcolor}
\usepackage{graphicx}
\usepackage{tikz,amssymb}
\usepackage{harpoon}
\usepackage[mathscr]{eucal}
\usepackage{amssymb}%
\usepackage[%backref=page,
linktocpage=true,colorlinks,citecolor=magenta,linkcolor=blue,urlcolor=magenta]{hyperref}
\hypersetup{
 colorlinks=true,
 linkcolor=blue,
 citecolor=blue,
 urlcolor=blue,
}

\providecommand{\U}[1]{\protect \rule{.1in}{.1in}}

\usepackage{amsmath,amssymb,mathrsfs}

\usepackage{amsmath,amssymb,amsthm,mathtools}
\usepackage{esint}
\usepackage{microtype}
\usepackage{enumitem}
\usepackage{xurl}
\newtheorem{theorem}{Theorem}[section]
\newtheorem{proposition}[theorem]{Proposition}
\newtheorem{lemma}[theorem]{Lemma}
\newtheorem{cor}[theorem]{Corollary}
\theoremstyle{remark}
\newtheorem{rem}[theorem]{Remark}
\numberwithin{equation}{section}
\DeclareMathOperator{\Sym}{Sym}

\newcommand{\R}{\mathbb R}

\allowdisplaybreaks[2]

\title[Sharp universal bound on the multiplicity]{Sharp universal bound on the multiplicity of\\ Dirichlet, Neumann, and Robin eigenvalues\\ on bounded plane domains}
\author{Zhijie Chen}
\address{Department of Mathematical Sciences, Yau Mathematical Sciences Center,
Tsinghua University, Beijing, 100084, China }
\email{zjchen2016@tsinghua.edu.cn}

\author{
Chen-Gang Long} 
\address{Yau Mathematical Sciences Center,
Tsinghua University, Beijing, 100084, China}
\email{lcg21@tsinghua.org.cn}

\author{Shihong Zhang}
\address{Yau Mathematical Sciences Center,
Tsinghua University, Beijing, 100084, China}
\email{shihong-zhang@tsinghua.edu.cn}
\date{}

\begin{document}
	\begin{abstract}
		We prove that for every $k\ge3$ (resp. $k=2$), the largest possible multiplicity of the \(k\)th
		Dirichlet eigenvalue on any bounded connected plane domain is \(k\) (resp. $3=k+1$). This upper bound requires no boundary regularity. Equality is attained
		on a disk with finitely many circular holes, where the multiple
		eigenvalue first occurs at index \(k\). We also prove the same sharp multiplicity bounds for Neumann and Robin eigenvalues on any bounded connected plane domain with Lipschitz boundary.
	\end{abstract}
	\maketitle
	
	\section{Introduction}

\subsection{Main results}
Let
$\Omega\subset\mathbb R^2$ be a bounded connected open set, and let \begin{align}(\alpha,\beta)\in\mathbb R^2\setminus\{(0,0)\}.\end{align}
If $\alpha\ne0$, assume in addition that $\Omega$ has
Lipschitz boundary.
We denote by $H_\Omega^{(\alpha,\beta)}$ the minus Laplacian $-\Delta$ with the generalized
Robin boundary condition
\[
    \alpha\partial_\nu u+\beta u=0
    \qquad\text{on }\partial\Omega,
\]
understood in the variational sense.
When $\alpha=0$, this is the Dirichlet Laplacian. When $\beta=0$, this is the Neumann Laplacian. When $\alpha\beta\neq 0$, this is the Robin Laplacian.

Under the additional condition that $\Omega$ has
Lipschitz boundary for $\alpha\neq 0$, 
it is well known that the eigenvalues of $H_\Omega^{(\alpha,\beta)}$, numbered from one and counted with multiplicity,
satisfy
\[
    \lambda_1^{(\alpha,\beta)}(\Omega)
    <
    \lambda_2^{(\alpha,\beta)}(\Omega)
    \le
    \lambda_3^{(\alpha,\beta)}(\Omega)
    \le\cdots,
    \qquad
    \lambda_k^{(\alpha,\beta)}(\Omega)\longrightarrow+\infty.
\]
Define the multiplicity of the eigenvalue $\lambda_k^{(\alpha,\beta)}(\Omega)$ by$$m_k(\Omega)=m_k^{(\alpha,\beta)}(\Omega):=\dim\ker\left(H_\Omega^{(\alpha,\beta)}-\lambda_k^{(\alpha,\beta)}(\Omega)\right).$$
Clearly $m_1^{(\alpha,\beta)}(\Omega)=1$ because $\Omega$ is connected. The research of estimating $m_k^{(\alpha,\beta)}(\Omega)$ for $k\geq 2$ has a long history; see the recent monograph of B\'erard and Helffer \cite{BH} for detailed introductions and rich references. Here we only mention some results related to our concern.

For the Dirichlet and Neumann Laplacian on a smooth bounded plane domain, 
Nadirashvili \cite[Theorem~2]{Nadirashvili} proved in 1988 that
$m_k\le 2k-1$. In particular, $m_2\leq 3$. Moreover, the multiplicity three occurs on a multiply connected domain for the second Dirichlet eigenvalue.
This example is due to M. Hoffmann-Ostenhof, T. Hoffmann-Ostenhof
and Nadirashvili \cite[Theorem~2.2]{HHOnodal}, with a corrected
proof by Helffer, Hoffmann-Ostenhof, Jauberteau and L\'ena
\cite[Theorem~2.1 and Section~5]{HHL}. They \cite[Proposition~4.2]{HHL} also proved the Dirichlet bound $m_2\le 3$ for bounded plane domains with piecewise $C^{1,+}$ boundary and corner angles in $(0,2\pi]$.
On the other hand, the sharpness of $m_2=3$ for the Neumann case was established in \cite[Theorem~1.4 and \S 6 (d)]{CdV}. Therefore, for the Dirichlet case $\alpha=0$ and the Neumann case $\beta=0$, 
\begin{align}\label{eq:maink=20}
    &\max_{\Omega\; \text{bounded connected, piecewise $C^{1,+}$ boundary}}m_2^{(0,1)}(\Omega)\notag\\
    =& \max_{\Omega\; \text{smooth bounded connected}}m_2^{(1,0)}(\Omega)=3.
\end{align}
In this paper, we will improve \eqref{eq:maink=20} to cover the Robin case and all bounded domains (with Lipschitz boundary for the Neumann and Robin cases), see our main result Theorem \ref{thm:main} below.

For general $k\geq 3$, Hoffmann-Ostenhof, Michor and Nadirashvili
\cite[Theorem~A]{HMN} stated the sharper bound $m_k\le 2k-3$
for the Dirichlet case. Berdnikov \cite[Section~4]{Berdnikov} identified
a gap in the proof for multiply connected domains. Recently, 
B\'erard and Helffer \cite[Theorem~1.1]{BH} revisited these results and improved them to cover Neumann and Robin eigenvalues.

\medskip
\noindent{\bf Theorem A.} \cite[Theorem~1.1]{BH} \emph{Given any $(\alpha,\beta)\in\mathbb R^2\setminus\{(0,0)\}$ and $k\geq 3$.
\begin{itemize}
\item[(1)] $m_k^{(\alpha,\beta)}(\Omega)\le 2k-2$ holds on every smooth bounded plane domain $\Omega$.
\item[(2)] $m_k^{(\alpha,\beta)}(\Omega)\le 2k-3$ when the smooth bounded plane domain $\Omega$ is simply connected.\end{itemize}}
 \medskip

B\'erard and Helffer \cite[Theorem~1.1]{BH} also considered general operators $-\Delta+V(x)$ with $V\in C^\infty$, and here we only state their theorem for $-\Delta$ for our purpose.

In view of Theorem A, an open long-standing problem is what the sharp upper bound of $m_k$ should be for all bounded domains (with Lipschitz boundary for the Neumann and Robin cases, but not necessarily smooth). To our knowledge, there are even no explicit conjectural formulas for the sharp upper bound of $m_k$ in the literature.

Recently, we proved in \cite[Theorem 1.3]{CLZ} that for the Dirichlet case, $m_k(\Omega)\leq k$ for any $k\geq 2$ when $\Omega$ is a convex bounded plane domain, which improves Lin's classical result \cite{Lin1987} for $k=2$ to all $k\geq 3$. This result motivates us to guess that the sharp upper bound of $m_k^{(\alpha,\beta)}(\Omega)$ for any bounded plane domain might be also $k$ for any $k\geq 3$.

The main result of this paper is to determine the largest possible value of $m_k^{(\alpha,\beta)}(\Omega)$ for any $(\alpha,\beta)\in\mathbb R^2\setminus\{(0,0)\}$ and each
index $k\ge2$, which indicates that our guess is true, and also improves previous results to the optimal version.

\begin{theorem}\label{thm:main}
Given any $(\alpha,\beta)\in\mathbb R^2\setminus\{(0,0)\}$, we have
\begin{equation}\label{eq:maink=2}
   \max_{\Omega}m_2^{(\alpha,\beta)}(\Omega)=3,
\end{equation}
\begin{equation}\label{eq:main}
    \max_{\Omega}m_k^{(\alpha,\beta)}(\Omega)=k,\quad\forall k\geq 3,
\end{equation}
where $\Omega$ ranges over all bounded connected  domains in
$\mathbb R^2$ with Lipschitz boundary. When $\alpha=0$, i.e., Dirichlet case, the Lipschitz boundary assumption can be removed.
Moreover, for every $k\geq2$, the maximum is attained by a bounded
connected domain $\Omega_k$ with smooth boundary.
\end{theorem}

Again the connectedness assumption is essential, see Remark \ref{Vital Remark}.  Given any $N\geq k$, 
for the union of $N$ disjoint congruent disks, the ground state
has multiplicity $N$, so $m_k^{(\alpha,\beta)}=N$.

For each $k\ge2$, we will prove that the maximum is attained by a smooth domain $\Omega_k$. More precisely, when $\alpha=0$, $\Omega_k$ is a disk with finitely many pairwise disjoint closed interior disks removed, while when $\alpha\neq0$, $\Omega_k$ can be taken to consist of finitely many rooms connected by thin tubes; see Figure \ref{fig:sec72-neumann}. More specifically, we construct $\Omega_k$ such that, for $k\geq 3$,
\begin{equation}\label{eq:cluster}
    \lambda_{k-1}^{(\alpha,\beta)}(\Omega_k)
    <\lambda_k^{(\alpha,\beta)}(\Omega_k)=\cdots=\lambda_{2k-1}^{(\alpha,\beta)}(\Omega_k)
    <\lambda_{2k}^{(\alpha,\beta)}(\Omega_k).
\end{equation}
Thus, the eigenvalue of multiplicity $k$
first occurs at index $k$.

\begin{rem}
For $\Omega$ being any bounded connected open subset in $\mathbb R^N$ with $N\geq 3$, we can not expect  upper bounds for the multiplicity of Dirichlet and Neumann eigenvalues. Indeed, it was already proved by Colin de Verdi\`ere \cite[Theorem 1.4]{CdV}
in 1987 that for any $N\geq 3$ and the Neumann case,
$$\sup_{\Omega \;\text{bounded connected with piecewise $C^1$ boundary}}m_{k}^{(1,0)}(\Omega)=+\infty,\quad\forall k\geq 2.$$
And recently, we proved in \cite[Theorem 1.1]{CLZ} that for any $N\geq 3$ and the Dirichlet case,
$$\sup_{\Omega \;\text{smooth bounded star-shaped}}m_{k}^{(0,1)}(\Omega)=+\infty,\quad\forall k\geq 2.$$
This is the reason why we only consider plane domains in this paper.
\end{rem}

\subsection{Ideas of the proofs} The approach of proving Theorem A in \cite{BH} can not be applied to prove Theorem \ref{thm:main}, and we need to develop different ideas.

First of all, we transform the eigenvalue problem on a planar domain into a smooth Schr\"odinger operator $H=-\Delta_g+V$ on a Riemannian two-sphere $(\mathbb{S}^2,g)$, where $V\in C^\infty(\mathbb{S}^2;\mathbb R)$, such that zero is an eigenvalue with the same $p$ and $m$ as the given planar eigenvalue (see \eqref{eq:counts} below for the definition of $(p,m)$); see Proposition \ref{prop:transfer}. 
The argument  uses stereographic projection to transfer $H_\Omega^{(\alpha,\beta)}$ from 
the planar domain to $(\mathbb S^2,g)$, and then regularizes it to obtain a family
of smooth Schr\"odinger operators whose resolvents 
converge to the transferred planar resolvent. Furthermore, we introduce local potential perturbations and 
 apply Brouwer's fixed point theorem to obtain a smooth
operator $H=-\Delta_g+V$ with exactly the same $p,m$.
 We list the eigenvalues of $H$ in nondecreasing order, counting multiplicities,
\[
    \lambda_1(H)<\lambda_2(H)\le\lambda_3(H)\le\cdots.
\]
We count the eigenvalues starting from one and with multiplicity. Therefore, the sharp upper bound in Theorem \ref{thm:main} reduces to Theorem \ref{Thm 2} below.

For an arbitrary smooth metric $g$ and smooth real potential $V$ on
the sphere, Nadirashvili \cite[Theorem~1(A)]{Nadirashvili} proved
that $m_k(H)\le 2k-1$, where $$m_k(H)=\dim\ker(H-\lambda_k(H)).$$
Kokarev \cite[Theorem~1.2]{Kokarev} showed that the same bound holds
for real potentials $V\in L^q(M)$ with $q>1$.
For smooth potentials, M. Hoffmann-Ostenhof, T. Hoffmann-Ostenhof
and Nadirashvili
\cite[Theorem~1 and the subsequent Remark~(ii)]{HHOsphere}
improved the bound to $m_k(H)\le 2k-3$ for $k\ge3$.
In particular, $m_2(H)\le3$ and $m_3(H)\le3$.
These two estimates provide the starting point for our proof.
For an eigenvalue $\lambda$ of a self-adjoint operator $H$ that is
bounded below and has compact resolvent, set
\begin{equation}\label{eq:counts}
    p:=\#\{j:\lambda_j(H)<\lambda\},\qquad
    m:=\dim\ker(H-\lambda).
\end{equation}
If $\lambda=\lambda_k(H)$, then
    \begin{align}\label{Intro f-1}
        \lambda=\lambda_k(H)=\lambda_{p+1}(H)=\lambda_{p+2}(H)=\cdots=\lambda_{p+m}(H), 
    \end{align}
and $ p+1\leq k\leq p+m$. Theorem \ref{thm:main} will follow from the following result.
\begin{theorem}\label{Thm 2}
Let $(M,g)$ be a smooth closed Riemannian surface diffeomorphic to
$\mathbb S^2$, and let $H=-\Delta_g+V$ with
$V\in C^\infty(M;\mathbb R)$. Let $m$ and $p$ be defined by
\eqref{eq:counts}. Then
\begin{align}\label{Intro f-2}
     m\le p+1\qquad\mathrm{for}\qquad p\geq 2.
\end{align}
Consequently, $m_k(H)\leq k$ for every $k\geq3$.
\end{theorem}

\begin{rem}
The upper bound $m_k(H)\leq k$ in Theorem \ref{Thm 2} is also sharp. Indeed, by
Theorem~\ref{thm:main}, for each $k\geq3$, there exists a smooth bounded
connected planar domain $\Omega_k$ such that $m_k(\Omega_k)=k$.
Applying Proposition~\ref{prop:transfer}, we obtain a smooth
Riemannian metric $\widetilde g$ on $\mathbb S^2$ and a smooth
Schr\"odinger operator $\widetilde H=-\Delta_{\widetilde g}+\widetilde V$
with $\widetilde V\in C^\infty(\mathbb S^2;\mathbb R)$ such that
$m_k(\widetilde H)=k$.
\end{rem}

We complete the proof of Theorem \ref{Thm 2} in Sections \ref{sec:streams}--\ref{sec:upper}.
On the sphere, the estimate to prove is $m\le p+1$ for $p\ge3$. 
Let $u$ and $v$ be
two linearly independent eigenfunctions with the same eigenvalue $\lambda$. Inspired by \cite{NTT2014,TT}, we construct a stream
function $s:=S(u,v)$. At each global maximum $y$ of $s$,
both eigenfunctions $u$ and $v$ vanish to the same order $d_y$. The key idea to show in Section \ref{sec:streams} that this one
point already forces at least $2d_y-1$ eigenvalues to lie strictly
below $\lambda$, counted with multiplicity. To this  goal, we shift
$\lambda$ to zero and apply a nonsingular double Moutard
deformation~\cite{NTT2014,TT}. As the potential varies, the multiplicity
$m$ at zero and the number $p$ of negative eigenvalues stay fixed.
Near $y$, the difference $s(y)-s$ vanishes to order $2d_y$, and
rescaling the deformed potential gives a radial limit determined by
$d_y$. One radial test function and its products with $\cos(j\theta)$
and $\sin(j\theta)$ for $1\le j<d_y$ then yield a $(2d_y-1)$-dimensional space on which the deformed quadratic form is negative definite.
Since the test functions associated with distinct maxima can be
supported in disjoint neighborhoods, the min-max principle gives
\begin{align}\label{Intro f-3}
    \sum_{y\in X}(2d_y-1)\le p,
    \quad\text{where}\; X:=\left\{y\in M:s(y)=\max_M s\right\}.
\end{align}

In the following, we prove \eqref{Intro f-2} by contradiction. Suppose that there exists a Schr\"odinger operator $H=-\Delta_g+V$ such that $m(H)\geq p(H)+2$ with $p(H)\geq3$, where
$$p(H)=\#\{j:\lambda_j(H)<0\},\quad m(H)=\dim\ker H.$$
In Section~\ref{sec:descent},
we choose a Schr\"odinger operator $H_{\min}$ such that $p(H_{\min})$ is minimal among all such operators, with $m(H_{\min})\ge p(H_{\min})+2$ and $p(H_{\min})\ge3$.
If some stream $S$ has a global maximum $x$ of order at least two,
we can obtain a smooth potential $\widehat V$ with
\[
    \widehat p<p(H_{\min}),\qquad \widehat m-\widehat p\ge m(H_{\min})-p(H_{\min})+1.
\]
If $\widehat p\le2$, this contradicts the known bounds; see \eqref{eq:low-level-bounds}. If $\widehat p\ge3$, this contradicts the
minimality of $p(H_{\min})$. Thus, for every independent pair in
$\ker H_{\min}$, all global maxima of its stream have order one.

In Section~\ref{sec:upper}, by choosing an orthonormal basis of
$\ker H_{\min}$, we identify the space of ordered orthonormal pairs
$(u,v)$ in $\ker H_{\min}$ with $\mathcal F_m$. Write $w=u+iv$.
As shown in Section~\ref{sec:descent}, we have $d_y=1$ at each global
maximum. Hence, by \eqref{Intro f-3}, each stream has at most $p$
global maxima, and the order-one expansion gives $w_z\ne0$ at each
of them. We consider the probability measure obtained by normalizing
$e^{ts}\,\mathrm{d}A_g$, which concentrates near the global maxima
as $t\to+\infty$. Using this measure, we construct a map
$F_t\colon\mathcal F_m\to\mathbb C^p$ that is nowhere zero for all
sufficiently large $t$. Moreover, its components satisfy a natural
rotation-equivariance property. However, an application of Stokes'
theorem on $\mathcal F_m$ shows that these components must have a
common zero whenever $p\le m-2$, yielding the desired contradiction.

Finally, we need to establish the sharpness of the upper bound
$m_k(\Omega)\le k$. To this end, for each fixed $k\ge3$, we construct
a domain $\Omega_k$ such that $m_k(\Omega_k)=k$. Inspired by the work
of Colin de Verdi\`ere \cite{CdV}, we begin with the construction of
a suitable graph and then replace each vertex by a room and each
edge by a thin tube. The key step is to approximate the spectral
data of the graph by those of a planar domain and then use spectral
stability to obtain the prescribed multiplicity.

More recently, the authors \cite{CLZ} developed this approach further
to construct domains in $\mathbb{R}^N$, $N\geq3$, with arbitrarily
prescribed multiplicities of Dirichlet eigenvalues. As pointed out
in \cite{CLZ}, the approach there cannot be applied directly to
the case $N=2$. This requires an essential modification of the
underlying graph construction. To this end, in
Section~\ref{sec:graph}, we construct a planar graph and a
corresponding well-signed matrix with exactly $k-1$ negative
eigenvalues and a $k$-dimensional kernel.

In Section~\ref{sec:domain},
for the Dirichlet problem, we follow the strategy developed in our
previous work \cite{CLZ}: we first pass from the divergence operator to a
Schr\"odinger operator and then approximate its spectrum by the
classical method (see\cite{DM,KP}) of removing small circular holes in a disk.
For the Neumann and Robin problems, we follow the approach of
Colin de Verdi\`ere \cite{CdV} and directly approximate the spectrum
of the matrix by constructing domains consisting of rooms connected
by thin tubes. Finally, applying the spectral perturbation 
Lemma~\ref{lem:stability}, we obtain the prescribed multiplicities
and complete the proof.

	\section{Transformation from planar domains to the Riemann sphere}\label{sec:stability}

In this section, we transfer the eigenvalue of the Laplacian
$H_\Omega^{(\alpha,\beta)}$ with Dirichlet, Neumann, or Robin
boundary conditions to the Schr\"odinger operator on $\mathbb{S}^2$,
while preserving the eigenvalue count \eqref{eq:counts} (the corresponding eigenvalue $\lambda$ transforms into the eigenvalue \(0\)).
This spectral transfer allows us to obtain the eigenvalue multiplicity bounds
of $H_\Omega^{(\alpha,\beta)}$ by studying the spherical Schr\"odinger operator in the following sections.

We first establish the stability of preserving the eigenvalue count \eqref{eq:counts} in a continuous family of compact self-adjoint operators.
Then, in Proposition~\ref{prop:transfer}, we obtain a Schr\"odinger operator
$\widetilde H=-\Delta_g+\widetilde V$ on $\mathbb{S}^2$
that has the same eigenvalue count \eqref{eq:counts} as $H_\Omega^{(\alpha,\beta)}$.

 Suppose $T$ is a compact self-adjoint
 operator on a real Hilbert space $\mathcal H$,
 and  $\mu>0$ is an eigenvalue of $T$.
Let 
\begin{equation}\label{eq:T eigenvalue}
    m=\dim  \ker(T-\mu I),\quad
    p=\#\{j:\lambda_j(T)>\mu\},\quad d_m=m(m+1)/2,
\end{equation}
where the eigenvalues are counted with multiplicity.
Let $P:\mathcal H\to \ker(T-\mu I)$ be the orthogonal projection onto $\ker(T-\mu I)$.

Suppose that $T(t)$ and $\widetilde T(t)$ are families of compact self-adjoint
operators on $\mathcal H$, with $t\in\mathbb R^{d_m}$ and $T(0)=T$.
Assume that $T(t)$ is $C^1$ in the operator norm and
$\widetilde T(t)$ is norm-continuous such that
\begin{equation}\label{eq:first derivative isomorphism}
    L:\mathbb R^{d_m}
\longrightarrow\operatorname{Sym}(\ker(T-\mu I)),
\quad
h\longmapsto P\,DT(0)[h]\Big|_{\ker(T-\mu I)}
\end{equation}
is an isomorphism, where $D$ denotes the Fr\'echet derivative (i.e., $DT(0)[h]=\frac{d}{ds}T(sh)\big|_{s=0}$),
$|_{\ker(T-\mu I)}$ denotes the restriction of an operator to $\ker(T-\mu I)$,
and
$\operatorname{Sym}(E)$ denotes the space of
self-adjoint linear operators on $E=\ker(T-\mu I)$, equipped with the trace inner product
$\langle M,N\rangle_{\operatorname{Sym}(E)}
:=\operatorname{tr}_E(MN)$.
Here, for any orthonormal basis
$u_1,\ldots,u_m$ of $E$ and any linear operator
$S:E\to E$,
\[
\operatorname{tr}_E(S)
:=\sum_{i=1}^m
\langle Su_i,u_i\rangle_{\mathcal H}.
\]
Remark that, under an orthogonal change of basis, the matrix $A$ of $S$ becomes
$Q^{\mathsf T} A Q$, where $Q$ is the orthogonal change-of-basis matrix.
By cyclicity of the trace,
$\operatorname{tr}(Q^{\mathsf T} A Q)
=\operatorname{tr}(A Q Q^{\mathsf T})=\operatorname{tr}(A)$. Thus,
taking $S=MN$ shows that $\langle M,N\rangle_{\operatorname{Sym}(E)}$
is independent of the choice of orthonormal basis of $E$.

Now, we first prove the stability of preserving the eigenvalue count $m,p$ in
 \eqref{eq:T eigenvalue} in this family of compact self-adjoint operators.

\begin{lemma}\label{lem:stability}
Under the above assumptions on $T,T(t)$ and $\widetilde T(t)$,
there exists $r_0>0$ such that for every $0<r<r_0$,
there exists $\delta(r)>0$ with the following property:
if 
\[
\sup_{t\in \overline B_r(0)}\|\widetilde T(t)-T(t)\|<\delta(r),
\]
then there exists
$t_*\in \overline B_r(0)$ satisfying
\[
\begin{aligned}
\dim\ker(\widetilde T(t_*)-\mu I)=m,\quad
\#\{j:\lambda_j(\widetilde T(t_*))>\mu\}=p.
\end{aligned}
\]
Here $\overline B_r(0)$ denotes the closed ball of radius $r$
centered at $0$ in $\mathbb R^{d_m}$.

\end{lemma}

\begin{proof} The proof is inspired by the argument that the strong Arnold
hypothesis implies the weak Arnold hypothesis
in \cite[\S~1]{CdVArnold}.
    Let $E=\ker(T-\mu I), Q=I-P$ and define 
    \begin{equation*}
        R_0=(T-\mu I)|_{E^\perp}.
    \end{equation*}
Since \(\mu>0\) is an isolated eigenvalue of the compact self-adjoint operator \(T\),
there exists $\varepsilon>0$ such that 
$\sigma(R_0)\cap[-2\varepsilon,2\varepsilon]=\varnothing$, where $\sigma(R_0)$ denotes the spectrum of $R_0$.

For any compact self-adjoint
operator $A$ defined on $\mathcal{H}$ satisfying 
$\|A-T\|<\varepsilon$, 
from the orthogonal decomposition
$\mathcal H=E\oplus E^\perp$,  we have 
\begin{equation}\label{eq:A matrix com}
\begin{aligned}
    A-\mu I&=
\begin{pmatrix}
P(A-\mu I)P|_E&P(A-\mu I)Q|_{E^\perp}\\
Q(A-\mu I)P|_E&Q(A-\mu I)Q|_{E^\perp}
\end{pmatrix}\\
&=
\begin{pmatrix}
P(A-\mu I)P|_E&(QAP|_E)^*\\
QAP|_E& Q(A-\mu I)Q|_{E^\perp}
\end{pmatrix}.
\end{aligned}
\end{equation}
We define the standard Frobenius-Schur complement \cite[\S 1.6]{Tretter}
\begin{equation}\label{eq:schur F}
     R=Q(A-\mu I)Q|_{E^\perp},\quad
   F(A)=P(A-\mu I)P|_E-(QAP|_E)^*R^{-1}QAP|_E.
\end{equation}

Since $\sigma(R_0)\cap[-2\varepsilon,2\varepsilon]=\varnothing$,
the standard spectral theorem and $\|R-R_0\|\leq\|A-T\|<\varepsilon$ yield
\begin{equation}\label{eq:R inverse estimate}
    \|R-R_0\|\,\|R_0^{-1}\|
    \leq \frac{\|R-R_0\|}{2\varepsilon}<\frac{1}{2}<1.
\end{equation}
Thus $R=R_0+(R-R_0)$ is invertible by
\cite[Chap.~IV, Theorem~1.16]{Kato}.
So, \eqref{eq:A matrix com} yields
\begin{equation}\label{eq:schur}
\begin{pmatrix}
I&-(QAP|_E)^*R^{-1}\\
0&I
\end{pmatrix}
(A-\mu I)
\begin{pmatrix}
I&0\\
-R^{-1}QAP|_E&I
\end{pmatrix}
=
\begin{pmatrix}
F(A)&0\\
0&R
\end{pmatrix}.
\end{equation}
Therefore,
this formula 
\eqref{eq:schur} and the invertibility of $R$ yield
\begin{equation}\label{eq:schur count old}
   \left\{ \begin{aligned}
\dim\ker(A-\mu I)&=\dim\ker F(A), \\
 \#\{j:\lambda_j(A)>\mu\}
&=\#\{j:\lambda_j(F(A))>0\}+\#\{j:\lambda_j(R)>0\}.
\end{aligned}\right.
\end{equation}
Moreover, since $\sigma(R_0)\cap[-2\varepsilon,2\varepsilon]=\varnothing$ and 
$\|R-R_0\|<\varepsilon$, we have 
\begin{equation}\label{eq:R eigenvalue number}
    \#\{j:\lambda_j(R)>0\}=\#\{j:\lambda_j(R_0)>0\}=\#\{j:\lambda_j(T)>\mu\}=p.
\end{equation}
Therefore, \eqref{eq:schur count old} 
and
\eqref{eq:R eigenvalue number} 
 yield
\begin{equation}\label{eq:schur count}
 \begin{aligned}
&\dim\ker(A-\mu I)=\dim\ker F(A), \\
 &\#\{j:\lambda_j(A)>\mu\}
=\#\{j:\lambda_j(F(A))>0\}+p.
\end{aligned}
\end{equation}

Now we choose $A=T(t)$ and define
\[
\Phi(t)=F(T(t)).
\]
Since $T(0)=T$ and $QTP|_E=0$, we have $\Phi(0)=F(T)=0$ and 
direct calculation gives $DF(T)[\cdot]=P \cdot P|_E$.
Thus, 
for any $h\in \R^{d_m}$, we have 
\begin{equation}
    D\Phi(0)[h]=DF(T)[DT(0)[h]]=P\,DT(0)[h]|_E.
\end{equation}
Therefore, 
assumption \eqref{eq:first derivative isomorphism} and the inverse function theorem imply the existence of a local inverse \(\Psi\) and
neighborhoods $0\in U\subset\mathbb R^{d_m}$
and $0\in V\subset\operatorname{Sym}(E)$ such that 
\begin{equation}\label{eq:inverse theorem S}
    \Psi: V\to U,\quad
\Psi(0)=0,\quad
\Phi(\Psi(X))=X.
\end{equation}

Since $T(t)$ is $C^1$, there exists $r_0>0$ such that for any $r\leq r_0$,
in 
the closed ball
$\overline B_r(0)\subset U\subset\mathbb{R}^{d_m}$ centered at $0$, we have 
\[
\sup_{t\in \overline B_r(0)}\|T(t)-T\|
<\frac{\varepsilon}{2}.
\]
Using the inverse function $\Psi$ in \eqref{eq:inverse theorem S}, there exists
 $s_0(r)>0$ such that 
   $\Psi(\overline B^{E}_{s_0}(0))\subset \overline B_{r}(0)$, where 
   $\overline B^{E}_{s_0}(0)\subset\operatorname{Sym}(E)$ is the closed ball of radius $s_0$ centered at $0$ in
    $\operatorname{Sym}(E)$.

Since  $F$ is smooth, for sufficiently small $\varepsilon>0$,
the map $A\mapsto F(A)$ is Lipschitz
on the set $\{A:\|A-T\|<\varepsilon\}$ with constant $K_F>0$, i.e.,
\begin{equation}\label{eq:Lip f}
\|F(A)-F(\widehat{A})\|\leq
 K_F\|A-\widehat{A}\|,\quad\forall \|A-T\|,\|\widehat{A}-T\|<\varepsilon.
\end{equation}
Choose $\delta(r)=\min\{\frac{s_0(r)}{K_F} ,\frac{\varepsilon}{2}\}$. Then if
\begin{equation}
    \sup_{t\in \overline B_{r}(0)}
\|\widetilde T(t)-T(t)\|<\delta(r),
\end{equation}
for any 
$X\in \overline B^{E}_{s_0}(0)$, we have $T(\Psi(X)),\widetilde T(\Psi(X))\in \{A:\|A-T\|<\varepsilon\}$.
Thus, formulas \eqref{eq:inverse theorem S} and \eqref{eq:Lip f} yield
\[
\sup_{\|X\|\le s_0}
\|X-F(\widetilde T(\Psi(X)))\|=
\sup_{\|X\|\le s_0}
\|F(T(\Psi(X)))-F(\widetilde T(\Psi(X)))\|\le K_F\delta\leq s_0.
\]
By Brouwer's fixed point theorem, there exists
$X_*\in \overline B^{E}_{s_0}(0)$ such that $X_*=X_*-  F(\widetilde T(\Psi(X_*)))$. Consequently,
$F(\widetilde T(\Psi(X_*)))=0.$
Choose $t_*=\Psi(X_*)\in \overline B_{r}(0)$, then \eqref{eq:schur count}
implies that 
\[
\begin{aligned}
\dim\ker(\widetilde T(t_*)-\mu I)=m,\quad \#\{j:\lambda_j(\widetilde T(t_*))>\mu\}=p.
\end{aligned}
\]
The proof is complete.
\end{proof}

Next, we transfer the eigenvalue count \eqref{eq:counts} of the Laplace operator 
\(H_\Omega^{(\alpha,\beta)}\) 
on a bounded connected domain $\Omega$ with Dirichlet, Neumann and Robin boundary conditions
 to that of the Schr\"odinger operator on \(\mathbb S^2\).

 Suppose that \(H\) is a self-adjoint operator bounded below with compact resolvent
 and that its eigenvalues satisfy
     $\lambda_1(H)\le\lambda_2(H)\le\cdots
     \nearrow +\infty$. For an eigenvalue $\lambda$ of $H$, define
     the eigenvalue count
 \begin{equation}\label{eq:counts of general H}
     p=\#\{j:\lambda_j(H)<\lambda\},\quad
     m=\dim\ker(H-\lambda I).
 \end{equation}
Here the eigenvalues $\lambda_j(H)$ are counted with multiplicity.
In this section, we consider two special cases of $H$.

In the case of a bounded connected domain $\Omega$, let $H=H_\Omega^{(\alpha,\beta)}$
be the minus Laplace operator with boundary condition  
\begin{equation}\label{eq:unified-bc}
 \alpha\partial_\nu u+\beta u=0\quad\text{on }\partial\Omega,\quad (\alpha,\beta)\in\R^2\setminus\{(0,0)\}.
\end{equation}
This boundary condition includes the Dirichlet boundary condition ($\alpha=0$), the Neumann boundary condition ($\beta=0$),
and the Robin boundary condition ($\alpha\beta\neq 0$).
Moreover,
for the Neumann and Robin cases $\alpha\neq 0$, we further assume that $\Omega$ is a bounded
Lipschitz domain. 
This assumption ensures that
\(H^1(\Omega)\) is compactly embedded into \(L^2(\Omega)\) and that the trace inequality
holds on the boundary so that $H_\Omega^{(\alpha,\beta)}$ satisfies the conditions of $H$.

This operator $H=H_\Omega^{(\alpha,\beta)}$ with boundary condition \eqref{eq:unified-bc}
corresponds to the symmetric bilinear form
\begin{equation}\label{eq:energy functional I}
    I_\Omega(u,v)=\int_\Omega\nabla u\cdot\nabla v\,dx+
 \begin{cases}
 0,&\alpha\beta=0,\\
 \displaystyle\frac{\beta}{\alpha}
  \int_{\partial\Omega}uv\,dS,
       &\alpha\beta\ne0,
 \end{cases} 
\end{equation}
where, in the boundary integral, \(u, v\) are understood as Sobolev traces on $\partial\Omega$,
\(dS\) is the arc length element and $I_\Omega$ is defined 
on \(u,v\in\mathcal X_\Omega\) with 
\[
 \mathcal X_\Omega=
 \begin{cases}
 H_0^1(\Omega),&\alpha=0,\\
 H^1(\Omega),&\alpha\ne0.
 \end{cases}
\]

In the case of an oriented smooth Riemannian two-sphere $\mathbb S^2$
with Riemannian metric $g$, let $H=-\Delta_g+V$, where $V\in L^q(\mathbb S^2;\mathbb R)$,
$1<q<\infty$, and \(\Delta_g=\operatorname{div}_g\nabla_g\).
By writing $V=V_b+V_s$ with $V_b$ bounded and $\|V_s\|_{L^q}$ arbitrarily
small, H\"older's inequality and the embedding
$H^1(\mathbb S^2)\hookrightarrow L^{2q/(q-1)}(\mathbb S^2)$ give
\[
 \left|\int_{\mathbb S^2}Vu^2\,\mathrm dA_g\right|
 \le \tfrac12\|\nabla_g u\|_2^2+C\|u\|_2^2 .
\]
Thus $(
\int_{\mathbb S^2}
(|\nabla_g u|^2+(V+C+1)u^2)\,\mathrm dA_g
)^{1/2}$ is equivalent to the $H^1$ norm.
Together with the compactness of $H^1(\mathbb S^2)\hookrightarrow L^2(\mathbb S^2)$,
the standard spectral theorem implies that 
$-\Delta_g+V$ also satisfies the conditions of $H$.

For the above two special cases of $H$, we prove that there exists a smooth 
spherical Schr\"odinger operator on $\mathbb S^2$ that preserves the eigenvalue 
count \eqref{eq:counts of general H}.

\begin{proposition}\label{prop:transfer}
    ~

\begin{enumerate}
    \item Suppose $\Omega$ is a bounded connected domain,
    let \(H=H_\Omega^{(\alpha,\beta)}\) with boundary condition
    \eqref{eq:unified-bc}, and assume that $\partial\Omega$ is Lipschitz if $\alpha\neq 0$.
    For any eigenvalue $\lambda$ of \(H_\Omega^{(\alpha,\beta)}\),
    there exist a smooth metric $g$ and a smooth potential $\widetilde V$ on
    $\mathbb S^2$ such that the smooth Schr\"odinger operator
    \(\widetilde H=-\Delta_g+\widetilde V\) satisfies
    \[
        \#\{j:\lambda_j(\widetilde H)<0\}=\#\{j:\lambda_j(H_\Omega^{(\alpha,\beta)})<\lambda\},~
        \dim\ker\widetilde H= \dim\ker(H_\Omega^{(\alpha,\beta)}-\lambda I).
    \]
    \item Let $g$ be an arbitrary smooth Riemannian metric
    on $\mathbb S^2$, and 
    let \(H=-\Delta_g+V\), where \(V\in L^q(\mathbb S^2;\R)\) for some \(1<q<\infty\) and 
$V\in C^\infty(\mathbb S^2\setminus\Sigma)$
for some finite set $\Sigma\subset\mathbb S^2$.
    For any eigenvalue $\lambda$ of $H$, there exists a smooth potential $\widetilde V$ on
    $\mathbb S^2$ such that the smooth Schr\"odinger operator
    \(\widetilde H=-\Delta_g+\widetilde V\) satisfies 
    \[
        \#\{j:\lambda_j(\widetilde H)<0\}=\#\{j:\lambda_j(H)<\lambda\},~
        \dim\ker\widetilde H= \dim\ker(H-\lambda I).
    \]
\end{enumerate}

\end{proposition}

To prove Proposition \ref{prop:transfer} (1), we need to  
transfer the Euclidean metric on $\Omega$ to a Riemannian metric $g$ on \(\mathbb S^2\).
Indeed, the stereographic projection
$\sigma:\mathbb S^2\setminus\{N\}
\longrightarrow\mathbb R^2$ is given by
\[
\sigma(y_1,y_2,y_3)
=
\left(
\frac{y_1}{1-y_3},
\frac{y_2}{1-y_3}
\right),
\]
where $N=(0,0,1)$ is the north pole. Then we have 
$(\Omega,dx_1^2+dx_2^2)\cong(\sigma^{-1}(\Omega)$, $\sigma^*(dx_1^2+dx_2^2))$, which is
an isometric diffeomorphism under the map $\sigma$. 
Using a standard smooth cutoff function, we can extend the metric
 $\sigma^*(dx_1^2+dx_2^2)$ on $\sigma^{-1}(\Omega)$
to a smooth Riemannian metric $g$ on $\mathbb S^2$
satisfying $g|_{\sigma^{-1}(\Omega)}=\sigma^*(dx_1^2+dx_2^2)$. Therefore, the map
\begin{equation}\label{eq:def pull back metric g}
    \sigma^{-1}:(\Omega,dx_1^2+dx_2^2)\to (\sigma^{-1}(\Omega),g)
\end{equation}
is an isometric diffeomorphism. Thus, we can define the isometric embedding
\begin{equation}\label{eq:def J}
    J:L^2(\Omega)\longrightarrow L^2(\mathbb S^2,dA_g),
\quad
(Jf)(y):=
\begin{cases}
f(\sigma(y)),&y\in\sigma^{-1}(\Omega),\\
0,&y\notin\sigma^{-1}(\Omega),
\end{cases}
\end{equation}
where $dA_g$ is the area measure on $(\mathbb S^2,g)$.

To prove Proposition \ref{prop:transfer},
we first prove the following approximation lemma, which uses a sequence of smooth
 Schr\"odinger operators to approximate the translation of $H_\Omega^{(\alpha,\beta)}$.

\begin{lemma}\label{lem:planar-resolvent-approximation}
Let $H=H_\Omega^{(\alpha,\beta)}$ be as defined in Proposition~\ref{prop:transfer}(1), 
and let the metric $g$ and the map $J$ be as defined above in 
\eqref{eq:def pull back metric g} and \eqref{eq:def J}.
Then there exists \(s_0>0\) such that, for any \(s>s_0\), there exist
\(V_j\in C^\infty(\mathbb S^2;\mathbb R)\), \(j\ge1\), satisfying
\[
-\Delta_g+V_j+sI \ge3I,\]
\[
\left\|
(-\Delta_g+V_j+s)^{-1}
-J(H+sI)^{-1}J^*
\right\|_{\mathcal L(L^2(\mathbb S^2))}\to 0
\quad (j\to\infty).
\]
Here $\mathcal L(L^2(\mathbb S^2))$ denotes the space of bounded linear operators from $L^2(\mathbb S^2)$ to itself.
\end{lemma}

\begin{proof}
Let $\phi$ be a positive first eigenfunction with respect to the eigenvalue 
$\lambda_1(H)$. Under different boundary conditions, $\phi$ has 
different standard 
properties:
\begin{enumerate}
    \item In the Dirichlet case ($\alpha=0$), $\phi\in C^{\infty}(\Omega)\cap H_0^1(\Omega)\cap L^{\infty}(\Omega) $
    and $\phi>0~\text{in }\Omega$;
    \item In the Neumann case ($\beta=0$), $\lambda_1(H)=0$ and $\phi=1$;
    \item In the Robin case ($\alpha\beta\neq 0$), $\phi\in C^{\infty}(\Omega)\cap C(\overline\Omega)$
and  $0<a\le\phi\le b<\infty~\text{in }\Omega$ (see  \cite[Theorem 4.5]{Arendt}).
\end{enumerate}
We next prove the three boundary results simultaneously.

Choose \(s_0=\max\{1,3-\lambda_1(H)\}\) and let $s>s_0$, then we have 
$\lambda_1(H)+s\geq 3$ .
We define the translated symmetric bilinear form
with respect to operator  $H+sI$
by
\begin{equation}\label{eq:def I_s}
    \begin{aligned}
       &I_s(u,v)
=I_\Omega(u,v)+s\langle u,v\rangle_{L^2(\Omega)}\\
=&\int_\Omega\nabla u\cdot\nabla v\,dx+ s\int_\Omega u v\,dx  +
\begin{cases}
0,&\alpha\beta=0,\\
\displaystyle\frac{\beta}{\alpha}
 \int_{\partial\Omega}uv\,dS,
      &\alpha\beta\ne0,
\end{cases} 
~ \forall u,v\in \mathcal X_\Omega.  
    \end{aligned}
\end{equation}
Using the Sobolev inequality and trace inequality, the symmetric bilinear form $I_s$
is coercive on $\mathcal X_\Omega$.

For any $u\in C^{\infty}(\mathbb S^2)$, define the corresponding function $\mathcal Au\in \mathcal X_\Omega$
\begin{equation}\label{eq:def Au}
 (\mathcal Au)(x)
:=\phi(x)u(\sigma^{-1}(x)),
\qquad x\in\Omega.   
\end{equation}
Since $\sigma$ is an isometric diffeomorphism and $\phi$ is the first eigenfunction of $H$,
for any $u,v\in C^{\infty}(\mathbb S^2)$, we have the standard identity for all three boundary conditions in \eqref{eq:def I_s}
\begin{equation}\label{eq:unified-ground-form}
    I_s(\mathcal Au,\mathcal Av)
    =\int_{\sigma^{-1}(\Omega)}(J\phi)^2
    \bigl(\nabla_g u\cdot\nabla_g v+(\lambda_1(H)+s)uv\bigr)\,dA_g.
\end{equation}
Moreover, based on the coercivity of
$I_s$ and the boundedness of the eigenfunction $\phi$,
 by choosing $u=v$ in \eqref{eq:unified-ground-form}, 
 we have 
\begin{equation*}
    \|\mathcal Au\|_{\mathcal X_\Omega}^2\leq  C I_s(\mathcal Au,\mathcal Au)\leq \tilde{C}\|u\|_{H^1(\mathbb S^2)}^2.
\end{equation*}
Thus $\mathcal A$ can be extended continuously to $H^1(\mathbb S^2)=\overline{C^{\infty}(\mathbb S^2)}^{H^1(\mathbb{S}^2)}$
such that $\mathcal A u\in \mathcal X_\Omega$.
Besides,
\eqref{eq:def Au} and \eqref{eq:unified-ground-form} hold for all $u,v\in H^1(\mathbb S^2)$.

Furthermore,  we have 
\begin{equation}\label{eq:dense of A}
    \overline{\operatorname{ran}\mathcal A}
^{\,\mathcal X_\Omega}=\mathcal X_\Omega.
\end{equation}
Indeed, in the Dirichlet case, for any $\eta\in C_c^{\infty}(\Omega)$, 
choose $u=(\frac{\eta}{\phi})\circ\sigma\in C_c^{\infty}(\sigma^{-1}(\Omega))$ 
and extend it by zero to $C_c^{\infty}(\mathbb{S}^2)$, then $\mathcal{A}u=\eta$.
Thus $C_c^{\infty}(\Omega)\subset{\operatorname{ran}\mathcal A}$ and \eqref{eq:dense of A}
holds. 
In the Neumann and Robin cases, for any $\eta\in C^{\infty}(\overline{\Omega})$,
since the positive eigenfunction $\phi$ is bounded below, we choose $u=(\frac{\eta}{\phi})\circ\sigma\in H^1(\sigma^{-1}(\Omega))$ 
and apply the Sobolev extension to $H^1(\mathbb{S}^2)$, then $\mathcal{A}u=\eta$.
Thus,
$C^{\infty}(\overline{\Omega})\subset{\operatorname{ran}\mathcal A}$ and \eqref{eq:dense of A}
holds.

Now, we begin to construct a sequence of operators converging to $J(H+sI)^{-1}J^*$.
We modify $I_s$ in \eqref{eq:unified-ground-form} so that its integral is defined on $\mathbb{S}^2$.
For any $\varepsilon>0$, define
\[
\phi_\varepsilon(y)=(J\phi)(y)^2+\varepsilon^2,
\quad
\psi_\varepsilon(y)
=\varepsilon\mathbf1_{\mathbb{S}^2\setminus\sigma^{-1}(\Omega)}(y),~y\in\mathbb{S}^2.
\]
For any $u,v\in H^{1}(\mathbb S^2)$, define
\begin{equation}\label{eq:def tildeI varepsilon}
    \widetilde{I}_\varepsilon (u,v)
=\int_{\mathbb{S}^2}\phi_\varepsilon
\bigl(\nabla_g u\cdot\nabla_g v+(\lambda_1(H)+s)uv\bigr)\,dA_g
+\int_{\mathbb{S}^2}\psi_\varepsilon uv\,dA_g.
\end{equation}
Since $\phi_\varepsilon(x)\geq \varepsilon^2$, the symmetric bilinear
form $\widetilde{I}_\varepsilon$ is also coercive.
By the Lax-Milgram theorem, for any $f\in L^2(\mathbb S^2)$,
there exists a unique $u_\varepsilon\in H^1(\mathbb S^2)$
such that 
\begin{equation}\label{eq:def u varepsilon}
    \widetilde{I}_\varepsilon(u_\varepsilon,v)
=\langle f,\sqrt{\phi_\varepsilon}v\rangle_{L^2(\mathbb S^2)},
\quad \forall v\in H^1(\mathbb S^2),
\end{equation}
We define the operator with respect to $ \widetilde{I}_\varepsilon$ by 
\begin{equation}\label{eq:def K varepsilon}
    \begin{aligned}
        K_\varepsilon:   L^2(\mathbb S^2)&\to L^2(\mathbb S^2),\\
K_\varepsilon f&=\sqrt{\phi_\varepsilon}u_\varepsilon.
    \end{aligned}
\end{equation}
We claim that 
\begin{equation}\label{eq:weighted-resolvent-limit}
\|K_\varepsilon-J(H+sI)^{-1}J^*\|_{\mathcal L(L^2(\mathbb S^2))}\longrightarrow0,~\varepsilon\to 0.
\end{equation}

To prove claim \eqref{eq:weighted-resolvent-limit}, we only need to 
prove that for any sequence $f_\varepsilon\rightharpoonup f$ in $L^2(\mathbb S^2)$
as $\varepsilon\to0$, up to a subsequence, we have 
\begin{equation}\label{eq:prove the claim}
    K_\varepsilon f_\varepsilon
\longrightarrow J(H+sI)^{-1}J^*f \quad\text{in }L^2(\mathbb S^2).
\end{equation}
Indeed, suppose by contradiction that claim \eqref{eq:weighted-resolvent-limit} is not true.
Then there exist $\delta>0$,
$\varepsilon_n\to0$, and $\|f_{\varepsilon_n}\|_{L^2(\mathbb S^2)}=1$ such that
\[
\|(K_{\varepsilon_n}-J(H+sI)^{-1}J^*)f_{\varepsilon_n}\|_{L^2(\mathbb S^2)}\ge\delta>0.
\]
Up to a subsequence, $f_{\varepsilon_n}\rightharpoonup f$ in $L^2(\mathbb S^2)$
as $n\to\infty$. Since \(J(H+sI)^{-1}J^*\) is compact from \(L^2(\mathbb S^2)\) to \(L^2(\mathbb S^2)\), letting $n\to\infty$ and using
\eqref{eq:prove the claim}, we obtain
$$\|J(H+sI)^{-1}J^*f-J(H+sI)^{-1}J^*f\|_{L^2(\mathbb S^2)}\ge\delta>0,$$
which is a contradiction.
So we next prove \eqref{eq:prove the claim}.

For any sequence $f_\varepsilon\rightharpoonup f$ in $L^2(\mathbb S^2)$
as $\varepsilon\to0$, let $u_\varepsilon$ be defined as in \eqref{eq:def u varepsilon}
with respect to $f_{\varepsilon}$. Then formulas \eqref{eq:unified-ground-form},
\eqref{eq:def tildeI varepsilon}, and \eqref{eq:def u varepsilon} imply that 
\begin{equation}\label{eq:energy identify}
    \begin{aligned}
&I_s(\mathcal Au_\varepsilon,\mathcal Av)
+\varepsilon^2\int_{\mathbb{S}^2}
\bigl(\nabla_g  u_\varepsilon\cdot\nabla_g v
+(\lambda_1(H)+s) u_\varepsilon v\bigr)\,dA_g\\
&\quad\quad\qquad\quad\;+\int_{\mathbb{S}^2\setminus\sigma^{-1}(\Omega)}\varepsilon u_\varepsilon v\,dA_g\\
=&\widetilde{I}_\varepsilon (u_{\varepsilon},v)
=\langle f_\varepsilon,
\sqrt{\phi_\varepsilon}v\rangle_{L^2(\mathbb S^2)},\quad \forall v\in H^1(\mathbb S^2).
\end{aligned}
\end{equation}
Choosing $v=u_\varepsilon$ in \eqref{eq:energy identify} and using the coercivity of $I_s$
and expression \eqref{eq:def Au},
we obtain the energy estimate 
\begin{equation}\label{eq:energy estimate}
    \begin{aligned}
    &\|\mathcal Au_\varepsilon\|_{\mathcal X_\Omega}^2
+\|\varepsilon u_\varepsilon\|_{H^1(\mathbb S^2)}^2
+\varepsilon^{-1}
\|\varepsilon u_\varepsilon\|_{L^2(\mathbb S^2\setminus\sigma^{-1}(\Omega))}^2\\
\le &C \langle f_\varepsilon,
\sqrt{\phi_\varepsilon}u_\varepsilon\rangle_{L^2(\mathbb S^2)}
\leq C \|f_{\varepsilon}\|_{L^2(\mathbb S^2)} \|\sqrt{\phi_\varepsilon}u_\varepsilon\|_{L^2(\mathbb S^2)}\\\leq& 
C  \|f_{\varepsilon}\|_{L^2(\mathbb S^2)}\left(\|(J\phi) u_\varepsilon\|_{L^2(\mathbb S^2)}+\|\varepsilon u_\varepsilon\|_{L^2(\mathbb S^2)}\right)\\
\leq&C \left(\|J(\mathcal Au_\varepsilon)\|_{L^2(\mathbb S^2)}+\|\varepsilon u_\varepsilon\|_{L^2(\mathbb S^2)}\right)= C \left(\|\mathcal Au_\varepsilon\|_{L^2(\Omega)}+\|\varepsilon u_\varepsilon\|_{L^2(\mathbb S^2)}\right).
    \end{aligned}
\end{equation}
The constant \(C\) may vary from line to line. Therefore, energy estimate \eqref{eq:energy estimate} implies that there exists a constant
$C>0$ independent of $\varepsilon$ such that 
\begin{equation}\label{eq:final energy estimate}
    \|\mathcal Au_\varepsilon\|_{\mathcal X_\Omega}^2
+\|\varepsilon u_\varepsilon\|_{H^1(\mathbb S^2)}^2
+\varepsilon^{-1}
\|\varepsilon u_\varepsilon\|_{L^2(\mathbb S^2\setminus\sigma^{-1}(\Omega))}^2
\le C.
\end{equation}
Since $\mathcal X_\Omega$ and $H^1(\mathbb S^2)$ are Hilbert spaces and the embeddings
$\mathcal X_\Omega\hookrightarrow L^2(\Omega)$ and
$H^1(\mathbb S^2)\hookrightarrow L^2(\mathbb S^2)$ are compact,
up to a subsequence,
\begin{equation}\label{eq:u limit}
    \mathcal Au_\varepsilon\rightharpoonup h
    \text{~in }\mathcal X_\Omega,\quad   \mathcal Au_\varepsilon\to h
\text{~in }L^2(\Omega),\quad
\varepsilon u_\varepsilon\to w
\text{~in }L^2(\mathbb S^2),\quad  \text{as }\varepsilon\to 0.
\end{equation}

Based on the energy estimate \eqref{eq:final energy estimate}, we have
$\|\varepsilon u_\varepsilon\|_{L^2(\mathbb S^2\setminus\sigma^{-1}(\Omega))}^2\leq C\varepsilon$,
so $w=0$ a.e. on $\mathbb S^2\setminus\sigma^{-1}(\Omega)$.
Moreover, for any $a>0$, the definition \eqref{eq:def Au} and 
estimate \eqref{eq:final energy estimate} imply that
\begin{equation*}
\|\varepsilon u_\varepsilon\|_{L^2(\sigma^{-1}(\{\phi\ge a\}))}
\leq \left\|\frac{\varepsilon}{a} (J\phi) u_{\varepsilon}\right\|_{L^2(\sigma^{-1}(\{\phi\ge a\}))}
\le\frac{\varepsilon}{a}
\|\mathcal Au_\varepsilon\|_{L^2(\Omega)}\leq \frac{C\varepsilon}{a}
\longrightarrow0.
\end{equation*}
Therefore, $w=0$ a.e. on $\sigma^{-1}(\{\phi\ge a\}),\forall a>0$.
Since $\phi>0$ in $\Omega$, we conclude that
$w=0$ a.e. on $\sigma^{-1}(\Omega)$.
Consequently,
\begin{equation}\label{eq:varepsilon u converge to 0}
    \varepsilon u_\varepsilon\longrightarrow0
\quad\text{in }L^2(\mathbb S^2),\quad \text{as }\varepsilon\to 0.
\end{equation}
Therefore, definitions \eqref{eq:def Au} and \eqref{eq:def K varepsilon} yield
\begin{equation}\label{eq:Kf energy estimate}
    \|K_\varepsilon f_\varepsilon
-J(\mathcal Au_\varepsilon)\|_{L^2(\mathbb S^2)}
=\|(\sqrt{\phi_\varepsilon}-J\phi)
u_\varepsilon\|_{L^2(\mathbb S^2)}\le\|\varepsilon u_\varepsilon\|_{L^2(\mathbb S^2)}
\to 0,\quad \varepsilon\to 0.
\end{equation}
Together with \eqref{eq:u limit} and the fact that $J$ is an isometric embedding, we have 
\begin{equation}\label{eq:Kf convergence}
    K_\varepsilon f_\varepsilon\longrightarrow J h
\quad\text{in }L^2(\mathbb S^2),\quad \text{as }\varepsilon\to 0.
\end{equation}

Passing to the limit $\varepsilon\to 0$ in \eqref{eq:energy identify} and using formulas  
\eqref{eq:final energy estimate}, \eqref{eq:varepsilon u converge to 0}, and 
$f_\varepsilon\rightharpoonup f$ in $L^2(\mathbb S^2)$, we obtain that 
for any $v\in H^1(\mathbb S^2)$,
\begin{equation}
    \begin{aligned}
       & I_s(\mathcal Au_\varepsilon,\mathcal Av)
+\varepsilon^2\int_{\mathbb{S}^2}
\bigl(\nabla_g  u_\varepsilon\cdot\nabla_g v
+(\lambda_1(H)+s) u_\varepsilon v\bigr)\,dA_g\\
&\qquad\qquad\quad\;+\int_{\mathbb{S}^2\setminus\sigma^{-1}(\Omega)}\varepsilon u_\varepsilon v\,dA_g\\
      = &(I_s(h,\mathcal Av)+o(1))+O(\varepsilon)+o(1)=\langle f_\varepsilon,
    \sqrt{\phi_\varepsilon}v\rangle_{L^2(\mathbb S^2)}\\
    \to&
    \langle f,(J\phi)v\rangle_{L^2(\mathbb S^2)}=\langle J^*f,\mathcal Av\rangle_{L^2(\Omega)},\quad \varepsilon\to 0.
    \end{aligned}
\end{equation}
Therefore, letting $\varepsilon\to 0$, we have $I_s(h,\mathcal Av)=\langle J^*f,\mathcal Av\rangle_{L^2(\Omega)}$ for all
$v\in H^1(\mathbb S^2)$. Since \eqref{eq:dense of A} shows the density of $\operatorname{ran}\mathcal A$ in $\mathcal X_\Omega$,
we have  
\[
I_s(h,\eta)=\langle J^*f,\eta\rangle_{L^2(\Omega)},
~\forall \eta\in\mathcal X_\Omega,
\]
which means $h=(H+sI)^{-1}J^*f$. Thus \eqref{eq:Kf convergence} can be rewritten as 
\begin{equation*}
    K_\varepsilon f_\varepsilon
\longrightarrow J(H+sI)^{-1}J^*f \quad\text{in }L^2(\mathbb S^2).
\end{equation*}
To sum up, we prove the claim \eqref{eq:weighted-resolvent-limit}.

Now, we have obtained that the operator $K_{\varepsilon}$ on
 $L^2(\mathbb S^2)$ converges to $J(H+sI)^{-1}J^*$. However,
  $K_{\varepsilon}$ may not be smooth, so we need to make a smooth 
  perturbation.
For any $k\in\mathbb{Z}_+$, let 
\[
\phi_{\varepsilon,k}
=e^{t_k\Delta_g}\phi_\varepsilon,~\psi_{\varepsilon,k}
=e^{t_k\Delta_g}\psi_\varepsilon,
\quad t_k\downarrow0,
\]
where $e^{t\Delta_g}$ is the standard heat semigroup operator on 
$\mathbb S^2$. By the standard heat kernel theory 
\cite[\S 7]{Gr09},
we have 
$\phi_{\varepsilon,k},\psi_{\varepsilon,k}
\in C^\infty(\mathbb S^2;\mathbb R)$ such that as $k\to\infty$,
\begin{equation}\label{eq:limit phi varepsilon k}
    \phi_{\varepsilon,k}\to\phi_\varepsilon,
\quad
\psi_{\varepsilon,k}\to\psi_\varepsilon
\quad\text{a.e.~}x\in\mathbb{S}^2.
\end{equation}
Moreover, the maximum principle yields
\begin{equation}\label{eq:positive Perturbation}
    \varepsilon^2\le\phi_{\varepsilon,k}
\le\|\phi\|_{L^\infty(\Omega)}^2+\varepsilon^2,
\quad
0\le\psi_{\varepsilon,k}\le\varepsilon,\quad\forall k\in\mathbb{Z}_+,\varepsilon>0.
\end{equation}

We regularize the coefficients of the symmetric bilinear form \eqref{eq:def tildeI varepsilon} by 
\begin{align}\label{eq:perturb I}
    \widetilde I_{\varepsilon,k}(u,v)
=&{}\int_{\mathbb S^2}\phi_{\varepsilon,k}
\bigl(\nabla_g u\cdot\nabla_g v+(\lambda_1(H)+s)uv\bigr)\,dA_g\notag\\
&+\int_{\mathbb S^2}\psi_{\varepsilon,k}uv\,dA_g,\quad u,v\in H^{1}(\mathbb S^2).
\end{align}
Similar to \eqref{eq:def u varepsilon} and
\eqref{eq:def K varepsilon}, for any $f\in L^2(\mathbb S^2)$,
there exists a unique $u_{\varepsilon,k}\in H^1(\mathbb S^2)$
such that 
\begin{equation}\label{eq:def u varepsilon k}
    \widetilde{I}_{\varepsilon,k}(u_{\varepsilon,k},v)
=\langle f,\sqrt{\phi_{\varepsilon,k}}v\rangle_{L^2(\mathbb S^2)},
\quad \forall v\in H^1(\mathbb S^2).
\end{equation}
Then the operator $K_{\varepsilon,k}$ corresponding to $\widetilde I_{\varepsilon,k}$ is given by
\begin{equation}
    \begin{aligned}
K_{\varepsilon,k} f:=\sqrt{\phi_{\varepsilon,k}}u_{\varepsilon,k}.
    \end{aligned}
\end{equation}

In \eqref{eq:perturb I}, for any $w\in H^{1}(\mathbb S^2)$, choosing 
$u=v=\frac{w}{\sqrt{\phi_{\varepsilon,k}}}$, we obtain that 
\begin{equation}\label{eq:expression I perturb}
    \begin{aligned}
&\widetilde I_{\varepsilon,k}
\left(
\frac{w}{\sqrt{\phi_{\varepsilon,k}}},
\frac{w}{\sqrt{\phi_{\varepsilon,k}}}
\right)\\
=&\int_{\mathbb S^2}
\left[
|\nabla_gw|^2+
\left(
\lambda_1(H)+s+\frac{\psi_{\varepsilon,k}}{\phi_{\varepsilon,k}}
+\frac{\Delta_g(\sqrt{\phi_{\varepsilon,k}})}
{\sqrt{\phi_{\varepsilon,k}}}
\right)w^2
\right]dA_g.
\end{aligned}
\end{equation}
We choose
\[
V_{\varepsilon,k}
=\lambda_1(H)
+\frac{\psi_{\varepsilon,k}}{\phi_{\varepsilon,k}}
+\frac{\Delta_g(\sqrt{\phi_{\varepsilon,k}})}
{\sqrt{\phi_{\varepsilon,k}}}
\in C^\infty(\mathbb S^2;\mathbb R).
\]
Therefore, \eqref{eq:def u varepsilon k}-\eqref{eq:expression I perturb} yield that 
\[
K_{\varepsilon,k}
=(-\Delta_g+V_{\varepsilon,k}+sI)^{-1}.
\]
We prove that these $K_{\varepsilon,k}
=(-\Delta_g+V_{\varepsilon,k}+sI)^{-1}$ are the operators we need.

In fact, for any $w\in H^{1}(\mathbb S^2)$,
formulas \eqref{eq:positive Perturbation}, \eqref{eq:perturb I} and \eqref{eq:expression I perturb} yield
\begin{equation*}
    \begin{aligned}
      &\langle(-\Delta_g+V_{\varepsilon,k}+sI)w,w\rangle_{H^{-1},H^1}=
    \widetilde I_{\varepsilon,k}\left(\frac{w}{\sqrt{\phi_{\varepsilon,k}}},\frac{w}{\sqrt{\phi_{\varepsilon,k}}}\right)\\
    \geq& \int_{\mathbb S^2}\phi_{\varepsilon,k}  (\lambda_1(H)+s) \left|\frac{w}{\sqrt{\phi_{\varepsilon,k}}}\right|^2\,\mathrm dA_g\geq 3\|w\|^2_{L^2(\mathbb{S}^2)}.
    \end{aligned}
\end{equation*}
This means $-\Delta_g+V_{\varepsilon,k}+sI\ge3I,\forall \varepsilon>0,k\ge 1$.

For any fixed $\varepsilon>0$, we prove that 
\begin{equation}\label{eq:smooth-weighted-limit}
\|K_{\varepsilon,k}-K_\varepsilon\|_{\mathcal L(L^2(\mathbb S^2))}
\longrightarrow0,\quad k\to\infty.
\end{equation}
Since $K_\varepsilon$ is compact,
similar to the proof of $\|K_\varepsilon-J(H+sI)^{-1}J^*\|\to 0$ in \eqref{eq:weighted-resolvent-limit},
to prove \eqref{eq:smooth-weighted-limit}, it suffices to prove that: 
for any sequence $f_{\varepsilon,k}\rightharpoonup f$ in $L^2(\mathbb S^2)$
as $k\to\infty$, there exists a subsequence such that
\begin{equation}\label{eq:K varepslion k converge}
    K_{\varepsilon,k} f_{\varepsilon,k}
\longrightarrow K_\varepsilon f\quad \text{in }L^2(\mathbb S^2),\text{ as }k\to\infty.
\end{equation}

Thus, we prove \eqref{eq:K varepslion k converge} in order to prove \eqref{eq:smooth-weighted-limit}.
Let $u_{\varepsilon,k}\in H^1(\mathbb S^2)$ be defined as in 
\eqref{eq:def u varepsilon k} with respect to $f_{\varepsilon,k}$.
Then 
\eqref{eq:positive Perturbation} and
\eqref{eq:perturb I} yield
\begin{equation*}
    \|u_{\varepsilon,k}\|^2_{H^1(\mathbb S^2)}\leq \varepsilon^{-2}  \widetilde I_{\varepsilon,k}(u_{\varepsilon,k},u_{\varepsilon,k})
    \leq C_{\varepsilon} \|f_{\varepsilon,k}\|_{L^2(\mathbb S^2)}\|u_{\varepsilon,k}\|_{H^1(\mathbb S^2)}.
\end{equation*}
Therefore, $\{u_{\varepsilon,k}\}_{k\geq 1}$ is bounded, and up to a subsequence,
\begin{equation}\label{eq:limit u varepsilon k}
   u_{\varepsilon,k}\rightharpoonup u
\text{~in }H^1(\mathbb S^2),\quad
u_{\varepsilon,k}\to u
\text{~in }L^2(\mathbb S^2). 
\end{equation}
Thus, letting $k\to\infty$ in \eqref{eq:def u varepsilon k} and using the dominated convergence theorem yields
\begin{equation}\label{eq:I limit result}
    \widetilde I_\varepsilon(u,v)
=\langle f,\sqrt{\phi_\varepsilon}v\rangle,
\quad \forall v\in H^1(\mathbb S^2).
\end{equation}
Together with \eqref{eq:def u varepsilon} and \eqref{eq:def K varepsilon}, this gives $K_\varepsilon f=\sqrt{\phi_\varepsilon}u$.
Thus, using the formula
\eqref{eq:limit phi varepsilon k}, \eqref{eq:positive Perturbation}, \eqref{eq:limit u varepsilon k} and 
dominated convergence theorem, we have
\[
\begin{aligned}
&\|K_{\varepsilon,k}f_{\varepsilon,k}-K_\varepsilon f\|_{L^2}=
\|\sqrt{\phi_{\varepsilon,k}}u_{\varepsilon,k}
-\sqrt{\phi_\varepsilon}u\|_{L^2}\\
\le&
\|\sqrt{\phi_{\varepsilon,k}}(u_{\varepsilon,k}-u)\|_{L^2}+
\|(\sqrt{\phi_{\varepsilon,k}}-\sqrt{\phi_\varepsilon})u\|_{L^2}\to 0, \text{as } k\to\infty.
\end{aligned}
\]
Therefore, we obtain \eqref{eq:K varepslion k converge}, and then \eqref{eq:smooth-weighted-limit} holds true.

Finally, we use the diagonal sequence argument:
first choose $\varepsilon_j\to0$, then choose $k_j$ sufficiently large such that
$\|K_{\varepsilon_j,k_j}
-K_{\varepsilon_j}\|_{\mathcal L(L^2(\mathbb S^2))}\leq \frac{1}{j},\quad \forall j\geq 1$.
Then, from \eqref{eq:weighted-resolvent-limit} and
\eqref{eq:smooth-weighted-limit}, we obtain
\[
\begin{aligned}
&\|(-\Delta_g+V_{\varepsilon_j,k_j}+sI)^{-1}
-J(H+sI)^{-1}J^*\|_{\mathcal L(L^2(\mathbb S^2))}\\
\le&
\|K_{\varepsilon_j,k_j}-K_{\varepsilon_j}\|_{\mathcal L(L^2(\mathbb S^2))}
+\|K_{\varepsilon_j}-J(H+sI)^{-1}J^*\|_{\mathcal L(L^2(\mathbb S^2))}
\longrightarrow0,\quad j\to\infty.
\end{aligned}
\]
Choosing $V_j=V_{\varepsilon_j,k_j}$ completes the proof.
\end{proof}

Now we use Lemma \ref{lem:stability} and Lemma 
\ref{lem:planar-resolvent-approximation} to prove Proposition \ref{prop:transfer}.

\begin{proof}[\bf Proof of Proposition~\ref{prop:transfer}]
We prove the two cases of Proposition \ref{prop:transfer} by a similar approach,
 with differences only in the operator construction.
Let the operator $H$ be defined as in Proposition \ref{prop:transfer} (1)-(2),
and let $\lambda$ be an eigenvalue of $H$. Define 
\begin{equation*}
    E=\ker(H-\lambda I),\quad \dim E=m,\quad d_m=m(m+1)/2.
\end{equation*}
Let $P:L^2(M)\to E$ 
be the orthogonal projection and $Q=I-P$, where the underlying space
in this proof is 
 $M=\Omega$ in part (1) and 
$M=\mathbb{S}^2$ in part (2).

To localize the perturbations away from the boundary
in part (1) and the singular set in part (2),
we choose a small open geodesic ball $B$ such that 
$B\Subset\Omega$ in part~(1) and
$B\Subset\mathbb S^2\setminus\Sigma$ in part~(2),
where $\Sigma$ denotes the finite singular set of $V$.

For any 
function $w\in C_c^\infty(B;\mathbb R)$ with $\|w\|_{L^{\infty}}$ sufficiently small,
we define the perturbation and Schur complement of $H-\lambda I$ as in equation \eqref{eq:schur F}:
\begin{equation}\label{eq:schur in prop 22}
    \begin{aligned}
 R(w)=Q(H+w-\lambda I)Q|_{E^\perp},\quad   
\mathcal F(w)
=PwP|_{E}
-PwQR(w)^{-1}QwP|_{E},
\end{aligned}
\end{equation}
where $L^2(M)=E\oplus E^{\perp}$.
Since \(\lambda\) is an isolated eigenvalue, if \(\|w\|_{L^{\infty}}\) is sufficiently small,
 as in the proof of \eqref{eq:R inverse estimate}, $R(w)$ is invertible. Thus, as in \eqref{eq:schur} and \eqref{eq:schur count old}, 
we obtain that 
\[
\dim \ker(H+w-\lambda I)=\dim \ker\mathcal F(w).
\]
In particular, if $\mathcal F(w)=0$ is the zero operator for a function $w$,
then $\dim \ker(H+w-\lambda I)=m$. 
We prove this proposition in three steps.

{\bf Step 1:} We find a suitable perturbation $w_*$ that preserves the eigenvalue count
\begin{equation}\label{eq:eigenvalue count of H+w}
    \dim \ker(H+w_*-\lambda I)=m, \quad
\#\{k:\lambda_k(H+w_*)<\lambda\}=p.
\end{equation}

Let $R(0)^{-1}$ be extended by zero on $E$.
Note that $\mathcal{F}:C_c^\infty(B;\mathbb R)\to\operatorname{Sym}(E)$,
direct differentiation yields that for any $\eta\in C_c^\infty(B;\mathbb R)$,
\[
\mathcal F(t\eta)
=tP\eta P|_E-t^2P\eta R(0)^{-1}\eta P|_E+O(t^3),\quad t\to 0^+.
\]
Therefore, we obtain the first and second derivatives at $0$:
\begin{equation}\label{eq:def first derivative F}
    L\eta :=D\mathcal F(0)[\eta]=P\eta P|_{E},\quad C_c^\infty(B;\mathbb R)\to\operatorname{Sym}(E),
\end{equation}
\begin{equation}\label{eq:def Q}
    \mathcal Q(\eta)
:=\frac12\pi_{(\operatorname{ran}L)^\perp}D^2\mathcal F(0)[\eta,\eta]=-\pi_{(\operatorname{ran}L)^\perp}P\eta R(0)^{-1} \eta P|_{E},
\quad \eta\in\ker L,
\end{equation}
where $\pi_{(\operatorname{ran}L)^\perp}$ is the orthogonal projection from 
$\operatorname{Sym}(E)$ to $(\operatorname{ran}L)^\perp$.
Here, we focus on the quadratic map $\mathcal Q(\eta)$ restricted to $\ker L\to (\operatorname{ran}L)^\perp$,
which implies the properties we need.

Now we state the following claim that 
\begin{equation}\label{eq:claim ind}
   \operatorname{ind}_-
\langle A,\mathcal Q(\cdot)\rangle_{\operatorname{Sym}(E)} \ge \dim(\operatorname{ran}L)^\perp,
\quad \forall 0\ne A\in(\operatorname{ran}L)^\perp,
\end{equation}
where $\operatorname{ind}_-$ denotes the negative index of the 
quadratic form $\langle A,\mathcal Q(\eta)\rangle_{\operatorname{Sym}(E)}=$  
$ -\operatorname{tr}_E
\bigl(AP\eta R(0)^{-1}\eta P|_E\bigr)$ for $\eta\in\ker L$.

If claim \eqref{eq:claim ind} holds true,
 when $\dim(\operatorname{ran}L)^\perp=0$, we choose $w_*=0$. Since $L$ is surjective,
there exists a
$d_m$-dimensional $W\subset C_c^\infty(B;\mathbb R)$ such that 
\begin{equation}\label{eq:F w property}
   \mathcal F(w_*)=0,\quad
D\mathcal F(w_*)|_W:
W{\longrightarrow}
\operatorname{Sym}(E) \text{~is an isomorphism}. 
\end{equation}
When $\dim(\operatorname{ran}L)^\perp>0$,
by \eqref{eq:claim ind} and the
finite-dimensional reduction in the proof of
\cite[Theorem~20.6]{AgrachevSachkov},
there exists a finite-dimensional
$Z\subset\ker L$ such that
\begin{equation}\label{eq:ind at Z}
    \operatorname{ind}_-
\bigl(\langle A,\mathcal Q(\cdot)\rangle|_Z\bigr)
\ge \dim(\operatorname{ran}L)^\perp,\quad \forall 0\ne A\in (\operatorname{ran}L)^\perp.
\end{equation}
Choose a finite-dimensional
$Y\subset C_c^\infty(B;\mathbb R)$
such that $L|_Y:Y\to\operatorname{ran}L$
is an isomorphism, and 
let \(F:=\mathcal F|_{Y\oplus Z}\) be the finite-dimensional restriction. 
Then 
\eqref{eq:def first derivative F} and \eqref{eq:def Q} yield
\begin{equation}\label{eq:restriction F Q}
\begin{aligned}
    \ker DF(0)=Z,\quad
\operatorname{coker}DF(0)\cong (\operatorname{ran}L)^\perp,\\
\operatorname{Hess}_0F(\eta)=\pi_{(\operatorname{ran}L)^\perp}D^2F(0)[\eta,\eta]=2\mathcal Q(\eta), ~\eta\in Z.
\end{aligned}
\end{equation}
For the quadratic mapping $\mathcal Q|_Z$ and \eqref{eq:ind at Z},
\cite[Lemma~20.8]{AgrachevSachkov} yields that there exists a regular zero $\eta_0\in Z$ such that 
$\mathcal Q(\eta_0)=0$ and 
$D(\mathcal Q|_Z)(\eta_0):
Z\to (\operatorname{ran}L)^\perp$
is surjective.
Then, based on \eqref{eq:restriction F Q}, \cite[Lemma~20.7]{AgrachevSachkov} yields that there exist $w_*\in Y\oplus Z\subset C_c^\infty(B;\mathbb R)$ with $\|w_*\|_{L^\infty}$ sufficiently small and a $d_m$-dimensional subspace
$W\subset Y\oplus Z\subset C_c^\infty(B;\mathbb R)$ that also satisfies \eqref{eq:F w property}.

Choose $w_*$ and $W$ as in \eqref{eq:F w property}, define 
\[
    E_*=\ker(H+w_*-\lambda I).
\]
\begin{equation*}
    U(w)u=u-R(w)^{-1}Q(wu)\quad \text{with}\quad w\in  C_c^\infty(B;\mathbb R),u\in E.
\end{equation*}
When $\|w\|_{L^\infty}$ is sufficiently small, based on \eqref{eq:schur in prop 22}
and the identity $P(H-\lambda I)=0$, we have 
\begin{equation}\label{eq:F expression}
\begin{aligned}
     (H+w-\lambda I)U(w)=&(P+Q)(H+w-\lambda I)U(w)\\
     =&\mathcal F(w)+0   =\mathcal F(w)
\end{aligned}
\end{equation}
as an operator from $E$ to $E\subset L^2(M)$.

Since 
$U(w)^*|_E=I_E$, formula \eqref{eq:F expression}
implies
$U(w)^*(H+w-\lambda I)U(w)=\mathcal F(w)$.
Differentiating this identity at $w_*$,
and using $\mathcal F(w_*)=0$ in \eqref{eq:F expression}, we have
\begin{equation}\label{eq:D F expression}
   D\mathcal F(w_*)[v]=U(w_*)^*vU(w_*) 
\end{equation}
and $U(w_*):E\to E_*$. Moreover, for any $u\in E_*$,  the identity
$Q(H+w_*-\lambda I)(Pu+Qu)=0$ yields
\begin{equation*}
    Qu=-R(w_*)^{-1}Qw_*Pu.
\end{equation*}
So $u=(P+Q)u=U(w_*)Pu$ and
 $U(w_*)$ is surjective. Together with 
$PU(w_*)=I_E$, 
we obtain that $U(w_*):E\to E_*$ is an isomorphism
and $\dim E_*=\dim E=m$.
Furthermore, when $w_*$ is sufficiently small, the spectral gap at 
$\lambda$ ensures
\[
\#\{k:\lambda_k(H+w_*)<\lambda\}=\#\{k:\lambda_k(H)<\lambda\}=p.
\]
Thus, we obtain \eqref{eq:eigenvalue count of H+w} for this choice of $w_*$.

{\bf  Step 2:} We use $w_*$ and  $W$ in \eqref{eq:F w property} to construct a 
convergent sequence of operators, and then use Lemma \ref{lem:stability}
 to find a smooth Schr\"odinger operator that preserves the eigenvalue count.

 In part (1), taking the Riemannian metric $g$, the isometric embedding
 $J$, and the smooth potentials $V_j$ as in Lemma \ref{lem:planar-resolvent-approximation},
 we have 
 \begin{equation}\label{eq:convergence of Kj}
K_j=(-\Delta_g+V_j+sI)^{-1}
\longrightarrow K=J(H+sI)^{-1}J^* \text{~in~} \mathcal L(L^2(\mathbb{S}^2)).\\
 \end{equation}
In part (2), keeping the metric $g$ fixed and choosing 
$V_j\in C^\infty(\mathbb S^2;\mathbb R)$ with
$V_j\to V$ in $L^q$, we can also obtain the convergence
\eqref{eq:convergence of Kj} with $J=I$.
Specifically, for any 
\(V\in L^q(M)\) and \(u,v\in H^1(M)\), we have 
\begin{equation}\label{eq:real-form}
\left|\int_M V uv\,dA_g\right|
\le \|V\|_{L^q}
   \|u\|_{L^{2q/(q-1)}}
   \|v\|_{L^{2q/(q-1)}}\le C_q\|V\|_{L^q}
   \|u\|_{H^1}
   \|v\|_{H^1}.
\end{equation}
Thus, \eqref{eq:real-form} and \cite[Theorem VI.3.6]{Kato}
yield the convergence
\eqref{eq:convergence of Kj} with $J=I$.
Moreover, by increasing $s_0$ if necessary,
we  further have
$$0\le K,K_j\le\frac12 I,\quad\forall s>s_0.$$

 Let $P_*$ be the orthogonal projection from $L^2(M)$ onto $E_*$. 
 Since \(P_*U(w_*)=U(w_*)\), we may write \eqref{eq:D F expression} as
 \begin{align}\label{eq:D F expression 2}
     D\mathcal F(w_*)[v]&=U(w_*)^*vU(w_*)=U(w_*)^*P_*vU(w_*)\notag\\
     &= U(w_*)^* (P_*v|_{E_*}) U(w_*).
 \end{align}
 Note that $U(w_*):E\to E_*$ is an isomorphism, so the map
 \[
 \operatorname{Sym}(E_*)\longrightarrow\operatorname{Sym}(E),
 \quad A\longmapsto U(w_*)^*AU(w_*)
 \]
 is also an isomorphism. Thus, 
 \eqref{eq:F w property} and \eqref{eq:D F expression 2} yield that the map
 \[
 W\longrightarrow\operatorname{Sym}(E_*),
 \quad v\longmapsto P_*v|_{E_*}
 \]
 is also an isomorphism. Since $\dim W=d_m$, there exists 
 a basis
 $w_1,\ldots,w_{d_m}$ of $W$ such that the map
 \[
 \mathbb R^{d_m}\longrightarrow
 \operatorname{Sym}(E_*),\quad
 h\longmapsto
 P_*\Bigl(\sum_{i=1}^{d_m}h_iw_i\Bigr)|_{E_*}
 \]
 is  an isomorphism.
 We choose a  family of perturbations
 \begin{equation}\label{eq:def wt}
     w(t)=w_*+
 \sum_{i=1}^{d_m}t_iw_i, 
 \end{equation}
 and define the perturbed operators of $K,K_j$
 \begin{equation}\label{eq:def final Tt}
    T(t)=J(H+w(t)+sI)^{-1}J^*,\quad
 \widetilde T_j(t)
 =(-\Delta_g+V_j+\widehat w(t)+sI)^{-1},
 \end{equation}
 where the function $\widehat w(t)=Jw(t)J^*$ is understood as a multiplication operator
 in part~(1), and $\widehat w(t)=w(t)$ in part~(2).
 Moreover, $\widehat w(t)\in C^{\infty}(\mathbb{S}^2)$.

 Let $\mu=(\lambda+s)^{-1}$, then formula \eqref{eq:eigenvalue count of H+w}
 yields that $T(0)=J(H+w_*+sI)^{-1}J^*$ satisfies
 \begin{equation}\label{eq:condition of Lem1}
     \dim \ker(T(0)-\mu I)=m,\quad  \#\{k:\lambda_k(T(0))>\mu\}=p.
 \end{equation}

For any $h\in\R^{d_m}$, the inverse operator differentiation formula and \eqref{eq:def wt}, \eqref{eq:def final Tt} yield
 \begin{equation}\label{eq:derivative of T0}
     DT(0)[h]
    =-J(H+w_*+sI)^{-1}\left(\sum_{i=1}^{d_m}h_i w_i\right)(H+w_*+sI)^{-1}J^*.
 \end{equation}
 Let $\tilde{E}=JE_*=\ker(T(0)-\mu I)$ and let \(\widetilde P=JP_*J^*\) be the orthogonal projection
 onto $\tilde{E}$.
 Since for any $u\in E_*$, $(H+w_*+sI)u=(\lambda+s)u=\mu^{-1}u$,
 based on expression \eqref{eq:derivative of T0}, we have 
 \begin{equation}\label{eq:condition of Lem2}
 \begin{aligned}
 \widetilde P DT(0)[h]|_{\widetilde E}
 &=-\widetilde PJ(H+w_*+sI)^{-1}
 \left(\sum_{i=1}^{d_m}h_i w_i\right)(H+w_*+sI)^{-1}J^*|_{\widetilde E}\\
 &=-\mu^2\widetilde PJ\left(\sum_{i=1}^{d_m}h_i w_i\right)J^*|_{\widetilde E}.
 \end{aligned}
 \end{equation}
 Therefore, $h\rightarrow  \tilde{P}DT(0)[h]|_{\tilde{E}}$
 is an isomorphism.

Moreover, since \(\|w_*\|_{L^\infty}\) can be chosen sufficiently small,
we can choose a  small closed ball
$\overline B_r(0)\subset
\mathbb R^{d_m}$ such that
$\|w(t)\|_{L^{\infty}}<1,\forall t\in \overline B_r(0)$. Then 
 for sufficiently large $j$, the resolvent identity $(A+B)^{-1}
 =(I+A^{-1}B)^{-1}A^{-1}$ and $0\le K,K_j\le\frac12 I$ 
yield
\begin{equation}\label{eq:condition of Lem3}
 \begin{aligned}
   &  \sup_{|t|\le r}
\|\widetilde T_j(t)-T(t)\|
= \sup_{|t|\le r}    \|(I+K_j \widehat w(t))^{-1}K_j-(I+K \widehat w(t))^{-1}K\|  \\
= & \sup_{|t|\le r}  \| (I+K_j\widehat w(t))^{-1}(K_j-K)(I+\widehat w(t)K )^{-1}\|\\
\leq & \sup_{|t|\le r}\frac{1}{1-\|K_j\widehat w(t)\|}  \|K_j-K\| \frac{1}{1-\|\widehat w(t) K\|}  \leq 4\|K_j-K\|\longrightarrow 0.
\end{aligned}
\end{equation}

 Therefore, formulas \eqref{eq:condition of Lem1}, \eqref{eq:condition of Lem2}, 
 and \eqref{eq:condition of Lem3} satisfy the conditions of
  Lemma \ref{lem:stability}, so there exist a sufficiently large $j$ and $t_*\in \overline B_r(0)$ 
 such that 
 \begin{equation}\label{eq:eigenvalue count of Tt*}
     \dim\ker(\widetilde T_j(t_*)-\mu I)=m,\quad
 \#\{k:\lambda_k(\widetilde T_j(t_*))>\mu\}=p.
 \end{equation}
 Finally, we choose 
 \[
 \widetilde V=V_j+\widehat w(t_*)-\lambda,
 \quad
 \widetilde H=-\Delta_g+\widetilde V.
 \]
 From \eqref{eq:def final Tt} and \eqref{eq:eigenvalue count of Tt*},
 we prove this proposition in both parts (1) and (2).

{\bf  Step 3:} Finally, we only need to prove \eqref{eq:claim ind} by constructing negative
subspaces of arbitrarily large dimension.

Fix $0\ne A\in(\operatorname{ran}L)^\perp$.
Write $T=R(0)^{-1}$, extended by zero on $E$, and choose a real
orthonormal basis $u_1,\ldots,u_m$ of $E$. On $C_c^\infty(B;\mathbb R)$, we set
\begin{equation}\label{eq: quadratic form perturbation}
 q_A(\eta)=-\sum_{i,j=1}^m A_{ij}
             \langle \eta u_i,T(\eta u_j)\rangle_{L^2(M)} .
\end{equation}
By \eqref{eq:def Q}, its restriction to $\ker L$ is
$\langle A,\mathcal Q(\cdot)\rangle_{\operatorname{Sym}(E)}$.

Take coordinate disks $D\Subset D_1\Subset B$ in isothermal coordinates.
On $D_1$ the area measure is $d\nu=a(x)\,dx$ and
$H-\lambda I=a^{-1}(-\Delta+v)$, where $a>0$ and $v$ are smooth.
In the planar case $a=1$ and $v=-\lambda$.
Set $$b(x,y)=\sum_{i,j} A_{ij}u_i(x)u_j(y).$$
Since $A\perp\operatorname{ran}L$, \eqref{eq:def first derivative F}
gives $$\int_D\eta(x)b(x,x)\,d\nu(x)=0,\quad\forall \eta\in C_c^\infty(D).$$
Thus $b(x,x)=0$ on $D$. The restrictions of the $u_i$ to $D$ are
linearly independent by unique continuation on the connected set
$\Omega$ or $\mathbb S^2\setminus\Sigma$; see \cite{Aronszajn}. In the
latter case $\Sigma$ has measure zero. Hence the vectors
$(u_1(x),\ldots,u_m(x))$, $x\in D$, span $\mathbb R^m$, and $A\ne0$
implies that $b$ is nonzero somewhere on $D\times D$.

Let $n$ be the least order for which a derivative in $y$ of $b(x,y)$
is nonzero at $y=x$ for some $x\in D$. Such an order exists. For at
least one $x$, $b(x,\cdot)$ is a nonzero solution of $$(-\Delta+v)b(x,\cdot)=0,$$
so strong unique continuation rules out a zero of infinite order at $y=x$.
Taylor's formula gives
\begin{equation}\label{eq:step3-taylor}
 b(x,x+z)=P_x(z)+O(|z|^{n+1}),\qquad
 P_x(z)=\left.\frac{1}{n!}(z\cdot\nabla_y)^n b(x,y)\right|_{y=x},
\end{equation}
locally uniformly in $x$. Here $P_{x_0}\ne0$ for some $x_0\in D$.
By symmetry, the same left-hand side equals
$$P_{x+z}(-z)+O(|z|^{n+1})=(-1)^nP_x(z)+O(|z|^{n+1}).$$
Thus $n=2\ell$ with $\ell\ge1$. Comparing terms of degree $n-2$ in
the equation for $b(x,\cdot)$ gives $\Delta_zP_x=0$.
The mean-value property now shows that $P_{x_0}$ takes both signs on
the unit circle. Choose a unit vector $\xi$ and a neighborhood
$D_0\Subset D$ of $x_0$ such that $(-1)^\ell P_x(\xi)>0$ on $D_0$.

We next determine the part of the kernel $G(x,y)$ of $T$, relative to
$d\nu$, that contributes to the index. On $D\times D$ it satisfies
\begin{equation}\label{eq:step3-log-kernel}
 G(x,y)=-\frac{1}{2\pi}\log|x-y|+O(1)\quad x\ne y,
\end{equation}
with a uniform remainder. Indeed, take $\zeta\in C_c^\infty(D_1)$ equal
to one near $\overline D$ and set
$$F_y(x)=-(2\pi)^{-1}\zeta(x)\log|x-y|,\quad y\in D,$$ extending by
zero outside $D_1$.
The planar fundamental solution gives
$$(H-\lambda I)F_y=\delta_y+r_y$$ in distributions relative to $d\nu$.
The $L^2(M,d\nu)$ norms of $F_y$ and $r_y$ are uniformly bounded
because the logarithm is locally square integrable and the derivatives
of $\zeta$ are supported away from $\overline D$. For sources supported
in $D$, the kernel of $T$ is therefore $(I-P)F_y-Tr_y$. To verify
this identity, integrate against $f(y)\,d\nu(y)$ for $f\in C_c^\infty(D)$.
The integral of $F_y$ is a smooth compactly supported function, to
which we can apply $T(H-\lambda I)=I-P$. On $D$, $PF_y$ is uniformly
bounded by finite dimensionality. Since $T$ is bounded on $L^2$ and
$(H-\lambda I)Tr_y=(I-P)r_y$, interior $H^2$ estimates and the
two-dimensional Sobolev embedding give a uniform bound for $Tr_y$ on $D$.
This proves \eqref{eq:step3-log-kernel}. 

The kernel is smooth off the
diagonal. Indeed, there each of $(H-\lambda I)_xG$ and $(H-\lambda I)_yG$
equals $-\sum_i u_i(x)u_i(y)$ by self-adjointness and
$(H-\lambda I)T=I-P$. The sum of these two differential operators is
elliptic on $D\times D$, so interior elliptic regularity applies.

It follows from \eqref{eq:step3-taylor} and \eqref{eq:step3-log-kernel}
that $K(x,y)=-b(x,y)G(x,y)$ extends continuously across the diagonal
by $K(x,x)=0$. Write
$$\mathcal B(\mu,\sigma)=\iint K(x,y)\,d\mu(x)\,d\sigma(y)$$ for finite
real measures compactly supported in $D$. Then
$\mathcal B(\eta\,d\nu,\eta\,d\nu)=q_A(\eta)$ for $\eta\in C_c^\infty(D;\mathbb R)$.
For $x\in D_0$ and small $h>0$, consider the finite difference
\begin{equation}\label{eq:step3-difference}
 \mu_{x,h}=h^{-\ell}\sum_{j=0}^{\ell}c_j\delta_{x+jh\xi},
 \qquad c_j=(-1)^{\ell-j}\binom{\ell}{j}.
\end{equation}
For $i\ne j$, \eqref{eq:step3-taylor} gives
$$b(x+ih\xi,x+jh\xi)=h^{2\ell}(j-i)^{2\ell}P_x(\xi)+O(h^{2\ell+1}).$$
The logarithm of the distance between these two points is
$\log h+\log|j-i|$. Since $K$ vanishes on the diagonal,
\eqref{eq:step3-log-kernel} yields
\begin{equation}\label{eq:step3-diagonal-pairing}
 \mathcal B(\mu_{x,h},\mu_{x,h})
 =\frac{(-1)^\ell(2\ell)!}{2\pi} P_x(\xi)\log h+O(1)
 \longrightarrow-\infty .
\end{equation}
Indeed, $\sum_{j=0}^\ell c_j j^k$ is zero for $0\le k<\ell$ and equals
$\ell!$ for $k=\ell$, by taking the $\ell$th difference of $t^k$.
Expanding $(j-i)^{2\ell}$ thus gives
$\sum_{i,j}c_i c_j(j-i)^{2\ell}=(-1)^\ell(2\ell)!$.

Now take any $N$ distinct points $x_1,\ldots,x_N$ in $D_0$.
For $r\ne s$, $\mathcal B(\mu_{x_r,h},\mu_{x_s,h})$ stays bounded
as $h\downarrow0$. It is a normalized $\ell$th difference in each
variable of $K$ on a product of disks with disjoint closures, where the
mixed derivatives of order $(\ell,\ell)$ are bounded. The bound
follows by applying the fundamental theorem of calculus to each
difference. By \eqref{eq:step3-diagonal-pairing} the resulting
$N\times N$ matrix is negative definite for all sufficiently small $h$.
Fix such an $h$ and replace each point mass in \eqref{eq:step3-difference}
by a nonnegative smooth density of unit mass relative to $d\nu$,
supported in a ball shrinking to that point. Continuity of $K$ gives
entrywise convergence of the matrices, so the resulting $N$ real test
functions in $C_c^\infty(D)$ span a negative definite subspace for $q_A$.
Its intersection with $\ker L$ has dimension at least $N-\operatorname{rank}L$.
Since $N$ was arbitrary, and replacing $A$ by $-A$ reverses the sign,
\begin{equation}\label{eq:infinity index}
 \operatorname{ind}_-(q_A|_{\ker L})
 =\operatorname{ind}_+(q_A|_{\ker L})=\infty .
\end{equation}
In particular, \eqref{eq:claim ind} follows. The proof is complete.
\end{proof}

\begin{rem}\label{Vital Remark}
The connectedness of $\Omega$ is essential in
Proposition~\ref{prop:transfer}(1).
Indeed, let
$\Omega=D_1\sqcup D_2$ be the union of two congruent disks with
disjoint closures, and have the same boundary condition
\eqref{eq:unified-bc} on both disks. Then the first eigenvalue $\lambda_1$ of $\Omega$ does not satisfy the conclusion of Proposition~\ref{prop:transfer}(1). Specifically, 
\[
    \#\{j:\lambda_j(H_\Omega^{(\alpha,\beta)})<\lambda_1\}=0,\quad
   \dim\ker(H_\Omega^{(\alpha,\beta)}-\lambda_1 I)=2.
\]
However, for any Schr\"odinger operator
\(\widetilde H=-\Delta_g+\widetilde V\) on $\mathbb S^2$, 
the first eigenvalue of $\widetilde H$ is simple, namely
\begin{equation*}
    \dim\ker(\widetilde H-\lambda_1(\widetilde H) I)=1<2=\dim\ker(H_\Omega^{(\alpha,\beta)}-\lambda_1 I).
\end{equation*}
Moreover, to prove Proposition~\ref{prop:transfer}(1), we used the connectedness of $\Omega$ in the proof of Lemma \ref{lem:planar-resolvent-approximation}. Indeed, the connectedness of $\Omega$ ensures that the first eigenfunction $\phi$ is strictly positive throughout
$\Omega$. This allows $\frac{1}{\phi}$ to be defined on all of 
$\Omega$ and yields the density of $\operatorname{ran}\mathcal A$
in $\mathcal X_\Omega$ as in \eqref{eq:dense of A}.

\end{rem}

\begin{rem}
In Proposition~\ref{prop:transfer}(1), the Dirichlet case requires no boundary regularity, and the Neumann and
Robin cases require only Lipschitz regularity on the boundary.
Indeed, the proof of Proposition~\ref{prop:transfer}(1)
requires the following properties:
\begin{enumerate}
    \item The Sobolev inequality and extension property
    for $\mathcal X_\Omega$;
    \item $\mathcal X_\Omega\hookrightarrow L^2(\Omega)$ is a compact embedding;
    \item In the Robin case, the Lipschitz boundary ensures that \(\inf_\Omega\phi>0\)
    (see \cite[Theorem 4.5]{Arendt}), and 
    the trace inequality ensures the coercivity of $I_s$ defined in
    \eqref{eq:def I_s}.
\end{enumerate}
All these properties hold under the stated assumptions,
without requiring additional smoothness of the boundary.

\end{rem}

\section{Stream functions and a change of potential}\label{sec:streams}
Let $M$ be an oriented smooth Riemannian two-sphere with metric $g$,
and let $H=-\Delta_g+V$, where $V$ is smooth and real-valued and
$\Delta_g f=\operatorname{div}_g(\nabla_g f)$.
For a one-form $\alpha$, we write
$\operatorname{div}_g\alpha:=\operatorname{div}_g(\alpha^\sharp)$,
where $\alpha^\sharp$ is the vector field dual to $\alpha$ with respect to $g$. Since $\nabla_g f=(\mathrm{d}f)^\sharp$ is the
vector field dual to $\mathrm{d}f$ with respect to $g$, we can write $\Delta_g=\operatorname{div}_g \mathrm{d}$. Write
\[
E=\ker H,\qquad m=\dim E,\qquad
p=\#\{j:\lambda_j(H)<0\}.
\]
Here the eigenvalues $\lambda_j(H)$ are counted with multiplicity.
Our main aim is to establish an inequality relating $m$ and $p$.

Our approach relies on the double Moutard transformation
in the spectral theory of planar Schr\"odinger operators.
Taimanov and Tsarev \cite{TT} used this transformation to construct
smooth potentials with a multidimensional $L^2$-kernel, and Novikov,
Taimanov, and Tsarev \cite{NTT2014} used it to construct a potential
with a multiple positive eigenvalue.

In our setting, the key idea is to use this transformation to relate
the local behavior of eigenfunctions to the negative index of $H$.
To this end, for $u,v\in E$, we observe that
$\operatorname{div}_g(u\,\mathrm{d}v-v\,\mathrm{d}u)=0$.
Since $M$ is simply connected, fixing a point $x_{\ast}\in M$, we can define
a unique stream function $S(u,v):M\to \mathbb R$ satisfying
\begin{equation}\label{Sec 3 def S}
    \mathrm{d}S(u,v)=\ast(u\,\mathrm{d}v-v\,\mathrm{d}u),
    \qquad S(u,v)(x_{\ast})=0,
\end{equation}
where $\ast$ denotes the Hodge star associated with $g$ and the
orientation of $M$. On $\{u\ne0\}$, $S(u,v)$ is the stream function associated with the
current density $u^2\nabla_g(v/u)$, where $u^2$ is the conductivity
and $v/u$ is the electric potential. For further discussion of its
physical interpretation and properties, see
\cite{BerryDennis2011,GN2019,TT}.
The main goal of this section is to establish an upper bound in terms
of $p$ for the sum of the weighted vanishing orders at the global
maxima of $S(u,v)$, namely,
\begin{align*}
    \sum_{y\in X}(2d_y-1)\le p,
\end{align*}
where $X$ is the set of global maxima of $S(u,v)$ and $d_y$ is the common
vanishing order of $u$ and $v$ at $y$.

To verify the stream function from the classical analysis, take positively oriented isothermal coordinates
$(z_1,z_2)$ near an arbitrary point of $M$, so that
$g=e^{2\omega_g}(\mathrm{d}z_1^2+\mathrm{d}z_2^2)$. We use the notations
$\nabla$, $\operatorname{div}$, and
$\Delta=\partial_{z_1}^2+\partial_{z_2}^2$ to denote the Euclidean
operators in these coordinates. Then
\[
\Delta_g=e^{-2\omega_g}\Delta,\qquad
\ast \mathrm{d}z_1=\mathrm{d}z_2,\qquad \ast \mathrm{d}z_2=-\mathrm{d}z_1.
\]
Moreover,
\[
\ast(u\,\mathrm{d}v-v\,\mathrm{d}u)
=(vu_{z_2}-uv_{z_2})\,\mathrm{d}z_1
 +(uv_{z_1}-vu_{z_1})\,\mathrm{d}z_2.
\]
Since $Hu=Hv=0$, we have
$\Delta u=e^{2\omega_g}Vu$ and $\Delta v=e^{2\omega_g}Vv$.
A direct computation gives
\begin{align*}
&\partial_{z_2}(vu_{z_2}-uv_{z_2})
-\partial_{z_1}(uv_{z_1}-vu_{z_1})\\
&\qquad=v\Delta u-u\Delta v
=e^{2\omega_g}(vVu-uVv)=0.
\end{align*}
Thus the one-form $\ast(u\,\mathrm{d}v-v\,\mathrm{d}u)$ is closed on $M$.
Since $M$ is simply connected, its integral along a piecewise smooth
path from $x_{\ast}$ to $z$ is independent of the path and defines a smooth
function satisfying \eqref{Sec 3 def S}. If another function $\widetilde S(u,v)$ satisfies
\eqref{Sec 3 def S}, then
$\mathrm{d}(\widetilde S(u,v)-S(u,v))=0$.
Since $M$ is connected and both functions vanish at $x_{\ast}$,
we have $\widetilde S(u,v)=S(u,v)$ on $M$.
This proves uniqueness.

In particular, in a small neighborhood of $x_{\ast}$, we choose positively oriented isothermal coordinates $(y_1,y_2)$, then
\begin{equation}\label{eq:stream}
S(u,v)(z)
=\int_{\gamma_{x_{\ast},z}}
\bigl[
(vu_{y_2}-uv_{y_2})\,\mathrm{d}y_1
+(uv_{y_1}-vu_{y_1})\,\mathrm{d}y_2
\bigr],
\end{equation}
where $\gamma_{x_{\ast},z}$ is any piecewise smooth path from $x_{\ast}$ to $z$
contained in that neighborhood.
Writing $S=S(u,v)$, we obtain, in every positively oriented
isothermal coordinate chart,
\begin{align}\label{Sec 3 S derivative}
S_{z_1}=vu_{z_2}-uv_{z_2},\qquad
S_{z_2}=uv_{z_1}-vu_{z_1}.
\end{align}
Furthermore,
the stream function $S=S(u,v)$ satisfies
the following three identities:
\begin{equation}\label{eq:stream-skew}
    S(u,v)=-S(v,u), \quad  S(u,v)=S(u,v-cu),  \end{equation}
\begin{equation}\label{Sec 3 Delta_g S}
    \Delta_g S=2\langle \mathrm{d}u,\ast \mathrm{d}v\rangle_g,
\end{equation}
\begin{equation}\label{Sec 3 identity S}
\begin{aligned}
    \operatorname{div}_g(u^{-2}\mathrm{d}S)&=0 \qquad \text{on } \{u\ne0\},\\
    \operatorname{div}_g(v^{-2}\mathrm{d}S)&=0 \qquad \text{on } \{v\ne0\}.
\end{aligned}
\end{equation}

Indeed, 
identity \eqref{eq:stream-skew} follows from the definition of \eqref{Sec 3 def S}. To prove \eqref{Sec 3 Delta_g S},
a direct computation in isothermal coordinates using
\eqref{Sec 3 S derivative} gives
\[
    \Delta S=-2(u_{z_1}v_{z_2}-v_{z_1}u_{z_2}),
\]
and hence
\[
    \Delta_g S=-2e^{-2\omega_g}(u_{z_1}v_{z_2}-v_{z_1}u_{z_2})
    =2\langle \mathrm{d}u,\ast \mathrm{d}v\rangle_g.
\]
On $\{u\ne0\}$, using $\langle \mathrm{d}u,\ast \mathrm{d}u\rangle_g=0$ and \eqref{Sec 3 Delta_g S} , we obtain
\begin{align*}
    \operatorname{div}_g(u^{-2}\mathrm{d}S)
    &=u^{-3}\bigl(u\Delta_g S-2\langle \mathrm{d}u,\mathrm{d}S\rangle_g\bigr)\\
    &=u^{-3}\bigl(u\Delta_g S
      -2\langle \mathrm{d}u,\ast(u\,\mathrm{d}v-v\,\mathrm{d}u)\rangle_g\bigr)=0.
\end{align*}
Interchanging $u$ and $v$ and using \eqref{eq:stream-skew}
gives the second identity in \eqref{Sec 3 identity S}.

For linearly independent $u,v\in E$, we first show that
at each local maximum of $s:=S(u,v)$, both $u$ and $v$ vanish to the same
order $d\ge1$, while $s-s(0)$ vanishes to order $2d$, see also \cite{BerryDennis2011}.

	\begin{lemma}\label{lem:maxima}
Let $u,v\in E$ be linearly independent and let $s:=S(u,v)$.
At each local maximum, choose positively oriented isothermal coordinates
$z=re^{i\theta}$ centered at that point, denoted by $0$. Then there exist
an integer $d\ge1$ and a real matrix
$P=\bigl(\begin{smallmatrix}a_1&a_2\\a_3&a_4\end{smallmatrix}\bigr)$
with $\det P>0$ such that
\[
    \begin{pmatrix} u(z) \\ v(z) \end{pmatrix}
    =  P\cdot
      \begin{pmatrix} \cos(d\theta) \\ \sin(d\theta) \end{pmatrix}r^d
      + O(r^{d+1})
\]
and
\begin{equation}\label{eq:max-exp}
    s(z)=s(0)+c_0|z|^{2d}+O(|z|^{2d+1}),
    \qquad c_0=-\frac{1}{2}\det P<0.
\end{equation}
\end{lemma}
\begin{proof}
Firstly, we claim $u(0)=v(0)=0$. 
If $u(0)\ne0$, then by \eqref{Sec 3 identity S},
$s$ satisfies $\operatorname{div}_g(u^{-2} \mathrm{d} s)=0$
in a neighborhood of $0$ where $u$ does not vanish.
The strong maximum principle implies that $s$ is constant
in a smaller connected neighborhood.
By \eqref{Sec 3 S derivative},
\[
\left(\frac{v}{u}\right)_{z_1}
=
\left(\frac{v}{u}\right)_{z_2}
=0.
\]
Thus $v=cu$ near $0$, and unique continuation
\cite{Aronszajn} gives $v=cu$ on $M$, a contradiction.
If $v(0)\ne0$, the same argument using the other identity
in \eqref{Sec 3 identity S} yields a contradiction.
Therefore, $u(0)=v(0)=0$.

In isothermal coordinates
$(z_1,z_2)$ near $0$, so that
$g=e^{2\omega_g}(\mathrm{d}z_1^2+\mathrm{d}z_2^2)$, then $\Delta u=e^{2\omega_g}Vu$ and
$\Delta v=e^{2\omega_g}Vv$.
Strong unique continuation and Taylor expansion give finite
vanishing orders and nonzero homogeneous harmonic leading terms
\cite{Aronszajn,Cheng}.
Writing $z=re^{i\theta}$, we have
\[
u=r^dA(\theta)+O(r^{d+1}),\qquad
v=r^\ell B(\theta)+O(r^{\ell+1}),
\]
where $d,\ell\ge1$, and $A,B$ are nonzero solutions of
\[
A''+d^2A=0,\qquad B''+\ell^2B=0.
\]
Therefore, there exist $a_i\in \R$ for $i=1,2,3,4$ such that
\begin{align}\label{Sec 3 A formula}
    A(\theta)=a_1\cos d\theta+a_2\sin d\theta,\qquad
B(\theta)=a_3\cos \ell\theta+a_4\sin \ell\theta.
\end{align}
Since $z_1=r\cos\theta$ and $z_2=r\sin\theta$, the chain rule gives
\[
\partial_r
=\cos\theta\,\partial_{z_1}+\sin\theta\,\partial_{z_2},
\qquad
\partial_\theta
=-r\sin\theta\,\partial_{z_1}
+r\cos\theta\,\partial_{z_2}.
\]
Combining these relations with \eqref{Sec 3 S derivative}, we obtain
\begin{align}\label{Sec 3 s_r}
\partial_rs
&=\cos\theta\,s_{z_1}+\sin\theta\,s_{z_2}\nonumber\\
&=v\bigl(\cos\theta\,u_{z_2}-\sin\theta\,u_{z_1}\bigr)
-u\bigl(\cos\theta\,v_{z_2}-\sin\theta\,v_{z_1}\bigr)\nonumber\\
&=\frac{v\,\partial_\theta u-u\,\partial_\theta v}{r}\nonumber\\
&=r^{d+\ell-1}(BA'-AB')+O(r^{d+\ell}).
\end{align}
Integrating along radial direction for \eqref{Sec 3 s_r} gives
\begin{equation}\label{Sec 3 asym S}
s(r,\theta)-s(0)
=\frac{r^{d+\ell}}{d+\ell}W(\theta)
+O(r^{d+\ell+1}),
\qquad W=BA'-AB',
\end{equation}
with the remainder uniform in $\theta$.

Suppose that $d\ne\ell$, by \eqref{Sec 3 A formula}, then $\int_0^{2\pi}W(\theta)\,\mathrm{d}\theta=0$.
If $W\not\equiv0$, it must change sign.
Choose $\theta^+$ such that $W(\theta^+)>0$.
Then \eqref{Sec 3 asym S} gives
$s(r,\theta^+)>s(0)$ for sufficiently small $r>0$,
contrary to the maximality of $s(0)$.
Hence $W\equiv0$.
However,
\[
W'=BA''-AB''=(\ell^2-d^2)AB,
\]
so $AB\equiv0$, which is impossible.
Therefore, $d=\ell$. We may now write
\[
A(\theta)=a_1\cos d\theta+a_2\sin d\theta,\qquad
B(\theta)=a_3\cos d\theta+a_4\sin d\theta.
\]
A direct computation gives the constant
$W=d(a_2a_3-a_1a_4)$.
By \eqref{Sec 3 asym S} and the maximality of $s(0)$,
we have $W\le0$.

Finally, we claim that $W<0$.
If $W=0$, then $a_2a_3=a_1a_4$ and we can get $B=cA$ for some constant $c$.
Set $\widetilde v=v-cu$.
Since $u$ and $v$ are linearly independent,
$\widetilde v$ is a nonzero element of $E$.
Its leading term of degree $d$ vanishes, so strong unique
continuation implies that its vanishing order
$\widetilde\ell$ is finite and satisfies $\widetilde\ell>d$.
On the other hand, from \eqref{eq:stream-skew}, we can see
\[
S(u,\widetilde v)=S(u,v-cu)=S(u,v)=s.
\]
Thus $0$ remains a maximum point of the stream function
associated with $u$ and $\widetilde v$.
Applying the preceding argument to this pair gives
$d=\widetilde\ell$, a contradiction.
Hence $W<0$. Substituting into \eqref{Sec 3 asym S}, we obtain
\eqref{eq:max-exp} with
\[
c_0=\frac{W}{2d}
=-\frac{a_1a_4-a_2a_3}{2}<0.
\]
The proof is complete.
\end{proof}

	We call \(d\) the order of the maximum. Thus \(d\) is the common
	vanishing order of \(u,v\), and \(s-s(0)\) vanishes to order \(2d\).
	
	\subsection{A change of potential at a fixed eigenvalue}
	We use the classical double Moutard transformation in the form given by the superposition formula of Taimanov and Tsarev \cite[Section~2.3 and Theorem~5]{TT2008}. We normalize all stream functions to vanish at the same point $x_{\ast}\in M$. Given $u,v\in E$, choose $c\in\mathbb R$ such that  $\tau:=c+S(u,v)$ has no zeros in $M$. For $f\in E$, define
\begin{equation}\label{eq:moutard}
    \widetilde f
    =f+\frac{uS(v,f)-vS(u,f)}{\tau},
    \qquad
    \widetilde V=V-2\Delta_g\log|\tau|.
\end{equation}
The following Lemma~\ref{lem:transform} shows that this transformation induces an isomorphism between the kernels of the corresponding Schr\"odinger operators.

\begin{lemma}\label{lem:transform}
If $\tau$ has no zeros on $M$, then $f\mapsto\widetilde f$ is a linear
isomorphism from $E:=\ker(-\Delta_g+V)$ onto
$\widetilde{E}:=\ker(-\Delta_g+\widetilde V)$. In particular, there holds $\dim\widetilde E=\dim E$.
\end{lemma}

\begin{proof}
The bilinearity of $S$ gives the linearity of $f\mapsto\widetilde f$.
Set $h=uS(v,f)-vS(u,f)$, so that $\widetilde f=f+h/\tau$.
By \eqref{Sec 3 identity S},
\begin{equation}\label{Sec 3 Identity-1}
    2\langle \mathrm{d}u,\mathrm{d}\tau\rangle_g=u\Delta_g\tau,
    \qquad
    2\langle \mathrm{d}v,\mathrm{d}\tau\rangle_g=v\Delta_g\tau.
\end{equation}
Since $\mathrm{d}S(v,f)=\ast(v\,\mathrm{d}f-f\,\mathrm{d}v)$, from \eqref{Sec 3 Delta_g S}, there holds
$\Delta_gS(u,v)=2\langle \mathrm{d}u,\ast \mathrm{d}v\rangle_g$; we obtain
\[
    2\langle \mathrm{d}u,\mathrm{d}S(v,f)\rangle_g
    =v\Delta_gS(u,f)-f\Delta_gS(u,v).
\]
Direct computation and using $S(v,u)=-S(u,v)$, we obtain
\begin{align}\label{Sec 3 Delta h}
    \Delta_g h
    &=Vh+2\langle \mathrm{d}u,\mathrm{d}S(v,f)\rangle_g
        -2\langle \mathrm{d}v,\mathrm{d}S(u,f)\rangle_g
        +u\Delta_gS(v,f)-v\Delta_gS(u,f)\nonumber\\
    &=Vh+\bigl(v\Delta_gS(u,f)-f\Delta_g\tau\bigr)
        -\bigl(u\Delta_gS(v,f)+f\Delta_g\tau\bigr)\nonumber
        \\
        &\qquad\qquad+u\Delta_gS(v,f)-v\Delta_gS(u,f)\nonumber\\
    &=Vh-2f\Delta_g\tau.
\end{align}
Moreover, by the definition \eqref{Sec 3 def S}, we compute  $u\,\mathrm{d}S(v,f)-v\,\mathrm{d}S(u,f)=-f\,\mathrm{d}\tau$, and hence
\begin{equation}\label{eq:moutard-dh}
    \mathrm{d}h=S(v,f)\,\mathrm{d}u-S(u,f)\,\mathrm{d}v-f\,\mathrm{d}\tau.
\end{equation}
Taking the inner product with $\mathrm{d}\tau$ and using
\eqref{Sec 3 Identity-1} gives
\begin{align}\label{Sec 3  Lem 2 euq-a}
    2\langle \mathrm{d}h,\mathrm{d}\tau\rangle_g
    =h\Delta_g\tau-2f|\mathrm{d}\tau|_g^2.
\end{align}
Therefore, using \eqref{Sec 3 Delta h} and \eqref{Sec 3  Lem 2 euq-a}, then we have
\begin{align*}
    \Delta_g\widetilde f
    &=Vf+\frac{\Delta_g h}{\tau}
      -\frac{2\langle \mathrm{d}h,\mathrm{d}\tau\rangle_g+h\Delta_g\tau}{\tau^2}
      +\frac{2h|\mathrm{d}\tau|_g^2}{\tau^3}\\
    &=\left(f+\frac h\tau\right)
      \left(V-2\frac{\Delta_g\tau}{\tau}
              +2\frac{|\mathrm{d}\tau|_g^2}{\tau^2}\right)
      =\widetilde V\widetilde f.
\end{align*}
Thus $\widetilde f\in\ker(-\Delta_g+\widetilde V)$.

To prove injectivity, suppose that $\widetilde f=0$. Then
$h=-\tau f$, and \eqref{eq:moutard-dh} gives
\[
    \tau f=vS(u,f)-uS(v,f),
    \qquad
    \tau\,\mathrm{d}f=S(u,f)\,\mathrm{d}v-S(v,f)\,\mathrm{d}u.
\]
Multiplying the second identity by $u$ and subtracting the first identity multiplied by $\mathrm{d}u$, the terms containing $uS(v,f)\,\mathrm{d}u$ cancel, and we obtain
\[
    \tau(u\,\mathrm{d}f-f\,\mathrm{d}u)=S(u,f)(u\,\mathrm{d}v-v\,\mathrm{d}u).
\]
Similarly, we obtain
\begin{align*}
    \tau(v\,\mathrm{d}f-f\,\mathrm{d}v)=S(v,f)(u\,\mathrm{d}v-v\,\mathrm{d}u).
\end{align*}
Applying $\ast$ to both identities, we obtain
\begin{equation}\label{eq:inject}
    \tau\,\mathrm{d}S(u,f)=S(u,f)\,\mathrm{d}\tau,
    \qquad
    \tau\,\mathrm{d}S(v,f)=S(v,f)\,\mathrm{d}\tau.
\end{equation}
Thus, $S(u,f)/\tau$ and $S(v,f)/\tau$ are constant on $M$. Since both vanish at $x_{\ast}$ and $\tau(x_{\ast})=c\ne0$, it follows that $h=0$, and hence $f=0$. Thus, the map $~\widetilde{}: f\to \widetilde{f}$ is injective, which leads to $\mathrm{dim}(E)\leq \mathrm{dim}(\widetilde{E})$.

 To prove the converse inequality $\mathrm{dim}(E)\geq \mathrm{dim}(\widetilde{E})$. We construct an injection
$\widehat{}:\widetilde E\to E$. Since $S(u,u)=S(v,v)=0$ and $S(v,u)=-S(u,v)$,
\eqref{eq:moutard} gives
\[
    \widetilde{u}=u-\frac{uS(u,v)}{\tau}=\frac{cu}{\tau},
    \qquad
    \widetilde{v}=v-\frac{vS(u,v)}{\tau}=\frac{cv}{\tau}.
\]
As $c=\tau(x_{\ast})\ne0$,  the functions
$a:=\widetilde u/c=u/\tau$ and $b:=\widetilde v/c=v/\tau$
belong to $\widetilde E$, so we can use this pair in the
transformation on $\widetilde E$.
Let $\widetilde S$ denote the stream on $\widetilde E$, normalized
at the same base point $x_{\ast}$. Then
\[
    \mathrm{d}\widetilde S(a,b)
    =\ast(a\,\mathrm{d}b-b\,\mathrm{d}a)
    =\frac{\ast(u\,\mathrm{d}v-v\,\mathrm{d}u)}{\tau^2}
    =\frac{\mathrm{d}\tau}{\tau^2}
    =\mathrm{d}\!\left(-\frac1\tau\right).
\]
Since $\widetilde S(a,b)(x_{\ast})=0$, we have
$\widetilde S(a,b)=-1/\tau+1/c$. Hence the same transformation, now applied to $\varphi\in\widetilde E$,
defines a map into $E$ given by
\[
    \widehat{\varphi}
    =\varphi+\frac{a\widetilde S(b,\varphi)
                       -b\widetilde S(a,\varphi)}{\widehat\tau}
    =\varphi-u\widetilde S(b,\varphi)+v\widetilde S(a,\varphi).
\]
Choosing  the constant $-1/c$
in \eqref{eq:moutard} gives the corresponding $\widehat{\tau}$, 
\[
    \widehat\tau=-\frac1c+\widetilde S(a,b)=-\frac1\tau.
\]
Hence,  $ \widehat\tau $ is nowhere zero, and the resulting potential is
\[
    \widetilde V-2\Delta_g\log|\widehat\tau|
    =\widetilde V+2\Delta_g\log|\tau|=V.
\]
The preceding injectivity argument applies to  $~\widehat{}:\varphi\to \widehat{\varphi}$ gives
$\dim\widetilde E\le\dim E$. Therefore, we obtain $\dim\widetilde E=\dim E $.
\end{proof}

The following constant-multiplicity lemma is an immediate consequence of the min--max principle. For the reader's convenience, we give the details.

\begin{lemma}\label{lem:negative-index}
Let $I\subset\mathbb R$ be an interval, and let $(V_t)_{t\in I}$
be a $C^0(M)$-continuous family of smooth real-valued potentials.
Set $H_t=-\Delta_g+V_t$. If $\dim\ker H_t$ is constant on $I$, then
the number of negative eigenvalues of $H_t$, counted with
multiplicity, is constant on $I$.
\end{lemma}

\begin{proof}
Write $\lambda_j(t)=\lambda_j(H_t)$ in nondecreasing order.
Fix $t_0\in I$,
for every $f\in H^1(M)$ with $\|f\|_{L^2(M)}=1$, the quadratic form
\[
    Q_t(f)=\int_M\bigl(|\nabla_g f|_g^2+V_t|f|^2\bigr)\,\mathrm{d}A_g
\]
satisfies
\[
    |Q_t(f)-Q_{t_0}(f)|
    =\left|\int_M(V_t-V_{t_0})|f|^2\,\mathrm{d}A_g\right|
    \le\|V_t-V_{t_0}\|_{C^0(M)}.
\]
Taking the supremum over the unit sphere of each $j$-dimensional
subspace $L\subset H^1(M)$ and then the infimum over $L$, the
min--max principle gives
\begin{align*}
    \lambda_j(t)
    =\inf_{\dim L=j}
      \sup_{\substack{f\in L\\ \|f\|_{L^2(M)}=1}}Q_t(f)
    \le&
    \inf_{\dim L=j}
      \sup_{\substack{f\in L\\ \|f\|_{L^2(M)}=1}}Q_{t_0}(f)
      +\|V_t-V_{t_0}\|_{C^0(M)}\\
    =&\lambda_j(t_0)+\|V_t-V_{t_0}\|_{C^0(M)}.
\end{align*}
Interchanging $t$ and $t_0$ yields
\begin{equation}\label{eq:potential-eigenvalue-continuity}
    |\lambda_j(t)-\lambda_j(t_0)|
    \le \|V_t-V_{t_0}\|_{C^0(M)}
    \qquad (j\ge1).
\end{equation}

Let $p$ be the number of negative eigenvalues of $H_{t_0}$ and let
$m=\dim\ker H_{t_0}$. Then $\lambda_{p+m+1}(t_0)>0$ and, if $p\ge1$,
$\lambda_p(t_0)<0$. Write $\lambda_0(t)=-\infty$ for convinence. Then by \eqref{eq:potential-eigenvalue-continuity}, for $t$ close to $t_0$, 
we have $\lambda_p(t)<0<\lambda_{p+m+1}(t)$.
Since
$\dim\ker H_t=m$,  then $\lambda_{p+1}(t)=\cdots= \lambda_{p+m}(t)=0$ if $m>0$.
Therefore, exactly $p$ eigenvalues of $H_t$ are negative for $t$ close to $t_0$.
This count is locally constant on the interval $I$ and therefore
constant on $I$.
\end{proof}

The double Moutard transformation
preserves the negative index $p$, while rescaling near each $y\in X$
yields a $(2d_y-1)$-dimensional negative subspace for the transformed
quadratic form. Since these subspaces have disjoint supports, their
direct sum is negative definite, hence we can get \eqref{eq:bubble-count}.

	\begin{lemma}\label{lem:bubbles}
Let $u,v\in E$ be linearly independent, set $s=S(u,v)$, and write
$s_*=\max_M s$ and $ X:=\{y\in M: s(y)=s_*\} $, then $X$ is finite. For each $y\in X$, let $d_y$ be the
common vanishing order of $u$ and $v$ at $y$. Then
\begin{equation}\label{eq:bubble-count}
    \sum_{y\in X}(2d_y-1)\le p.
\end{equation}
% Given pairwise disjoint open neighborhoods $O_y$ of the points $y\in X$, for all sufficiently large $t>0$ there exist subspaces $Z_{t,y}\subset C_c^\infty(O_y)$ of dimension $2d_y-1$ on which the quadratic form of $-\Delta_g+V_t$ is negative definite.
\end{lemma}

\begin{proof}
Firstly, we claim that $X$ is finite. By Lemma \ref{lem:maxima},  every global maximum is isolated.
The set of global maxima is closed in the compact manifold $M$,
and hence the set $X$ is finite. 

Let $\tau_t=s-s_*-1/t$ for $t>0$ and
$V_t=V-2\Delta_g\log|\tau_t|$,  
denote the quadratic form of $-\Delta_g+V_t$ by $Q_t$.
Since $\tau_t<0$ on $M$, Lemma~\ref{lem:transform} gives
\[
    \dim\ker(-\Delta_g+V_t)=\dim E
    \qquad\text{for every }t>0.
\]
Moreover,
\[
    V_t
    =V-2\Delta_g\log\bigl(1+t(s_*-s)\bigr)
    \longrightarrow V
    \qquad\text{in }C^0(M)\quad\text{as }t\downarrow0.
\]
Thus, setting $V_0=V$, we obtain a continuous family $(V_t)_{t\ge0}$
with constant kernel dimension. By Lemma~\ref{lem:negative-index},
$-\Delta_g+V_t$ has exactly $p$ negative eigenvalues
for every $t>0$.

Fix $y\in X$, write $d=d_y$, and choose isothermal coordinates
$z$ centered at $y$ in which $g=e^{2\omega_g}|\mathrm{d}z|^2$. Set
$$c_y:=-c_0=\frac12\det P>0,\qquad \ell_t=(t c_y)^{-1/(2d)}.$$ By \eqref{eq:max-exp} and its
derivative estimates,
\[
    t|\tau_t(\ell_t z)|=1+t\bigl(s(0)-s(\ell_t z)\bigr)
    \longrightarrow 1+|z|^{2d}
    \qquad\text{in }C^2_{\mathrm{loc}}(\mathbb R^2)
    \quad\text{as }t\to\infty.
\]
Since $\Delta_g=e^{-2\omega_g}\Delta_z$, it follows that
\[
    \ell_t^2(e^{2\omega_g} V_t)(\ell_t z)
    \longrightarrow
    -2\Delta_z\log(1+|z|^{2d})
    =-\frac{8d^2|z|^{2d-2}}{(1+|z|^{2d})^2}
    =:U_d(|z|)
\]
locally uniformly. Fix $f\in C_c^\infty(\mathbb R^2)$ in the chosen isothermal coordinates. For sufficiently large $t$,
the support of $f_{\ell_t}(\cdot):=f(\cdot/\ell_t)$ lies in the coordinate neighborhood, so we
extend $f_{\ell_t}$ by zero to $M$. Since $g=e^{2\omega_g}|\mathrm{d}z|^2$,
we have $\mathrm{d}A_g=e^{2\omega_g}\,\mathrm{d}z$ and
$|\nabla_g f_{\ell_t}|_g^2=e^{-2\omega_g}|\nabla f_{\ell_t}|^2$.
Thus, making the
change of variables, we obtain
\begin{align*}
    Q_t\bigl(f_{\ell_t}\bigr)
    =&\int_M \bigl(|\nabla_g f_{\ell_t}|_g^2+V_t f_{\ell_t}^2\bigr)
      \,\mathrm{d}A_g
    =\int_{\mathbb R^2}
      \left(|\nabla f|^2+\ell_t^2(e^{2\omega_g} V_t)(\ell_t z)f^2\right)\mathrm{d}z\\
    &\longrightarrow
    \int_{\mathbb R^2}\left(|\nabla f|^2+U_d f^2\right)\mathrm{d}z
    =:q_d(f).   
\end{align*}

Let $Y(\theta)\in \mathcal{H}_{d-1}$, where $\mathcal{H}_{d-1}$ denotes the space of real trigonometric polynomials of degree at most $d-1$. Then $\mathrm{dim}\mathcal{H}_{d-1}= 2d-1$, and by writing $Y(\theta)=a_0+\sum_{k=1}^{d-1}
(a_k\cos(k\theta)+b_k\sin(k\theta))$, it is elementary to see that
\begin{align}\label{Sec 3 Lem 3 equ-a}
    \int_0^{2\pi}|Y'|^2\,\mathrm{d}\theta
    \le (d-1)^2\int_0^{2\pi}Y^2\,\mathrm{d}\theta.
\end{align}
Fix $\chi\in C_c^\infty((-3,3))$ such that $\chi=1$ on $[-1,1]$
and $|\chi'|\le1$. Write $\xi=\log r$, set $w(\xi)=1/\cosh(d\xi)$,
and for $Y\in\mathcal H_{d-1}$, we define
\[
    f_Y(r,\theta)
    =\chi(\log r)\frac{2r^d}{1+r^{2d}}Y(\theta)
    =\chi(\xi)w(\xi)Y(\theta),
\]
using $|\nabla f_Y|^2=|\partial_rf_Y|^2+r^{-2}|\partial_{\theta}f_Y|^2$, we obtain
\begin{align*}
    |\nabla f_Y|^2+U_df_Y^2
    =\left[r^{-2}|(\chi w)'|^2+U_d(\chi w)^2 \right]Y^2(\theta)+r^{-2}\chi^2w^2 (Y^{'}(\theta))^2
\end{align*}
Since $r^2U_d(r)=-2d^2w^2$ and $w''=d^2(1-2w^2)w$, we have
\[
    |(\chi w)'|^2+r^2U_d(r)(\chi w)^2
    =\bigl(|\chi'|^2-d^2\chi^2\bigr)w^2+(\chi^2ww')'.
\]
The last term integrates to zero. Integrating along the radial direction and using \eqref{Sec 3 Lem 3 equ-a},  we obtain
\begin{align}\label{Sec 3 Lem 3 equ-b}
q_d(f_Y)
    &=\|Y\|_2^2\int_{\mathbb R}|\chi'|^2w^2\,\mathrm{d}\xi
     +\bigl(\|Y'\|_2^2-d^2\|Y\|_2^2\bigr)
      \int_{\mathbb R}\chi^2w^2\,\mathrm{d}\xi\nonumber\\
    &\le \|Y\|_2^2\int_{\mathbb R}
      \bigl(|\chi'|^2-(2d-1)\chi^2\bigr)w^2\,\mathrm{d}\xi\nonumber\\
    &\le \|Y\|_2^2
      \left(\int_{|\xi|>1}w^2\,\mathrm{d}\xi
            -\int_{|\xi|\le1}w^2\,\mathrm{d}\xi\right)\nonumber\\
    &=\frac{2}{d}(1-2\tanh d)\|Y\|_2^2<0.
\end{align}
Set $A=\{z\in\mathbb R^2:e^{-3}\le |z|\le e^3\}$. Each $f_Y$ is
supported in $A$, and separation of variables gives
\[
    \|f_Y\|_{L^2(A)}^2
    =\|Y\|_2^2\int_{-3}^{3}\chi(\xi)^2w(\xi)^2
      e^{2\xi}\,\mathrm{d}\xi.
\]
Since $\ell_t^2(e^{2\omega_g}V_t)(\ell_t\,\cdot)$ converges uniformly
to $U_d$ on $A$, the change of variables above yields
\[
    \bigl|Q_t((f_Y)_{\ell_t})-q_d(f_Y)\bigr|
    \le
    \bigl\|\ell_t^2(e^{2\omega_g}V_t)(\ell_t\,\cdot)-U_d
    \bigr\|_{L^\infty(A)}\|f_Y\|_{L^2(A)}^2
    =o(1)\|Y\|_2^2,
\]
where $o(1)\to0$ as $t\to\infty$, independently of $Y$.
Combining this with the preceding estimate \eqref{Sec 3 Lem 3 equ-b}, we obtain
\[
    Q_t((f_Y)_{\ell_t})
    \le \left[\frac{2}{d}(1-2\tanh d)+o(1)\right]\|Y\|_2^2.
\]
Thus, for all sufficiently large $t$, $Q_t$ is negative definite on
\[
    Z_{t,y}:=\bigl\{(f_Y)_{\ell_t}:Y\in\mathcal H_{d-1}\bigr\}.
\]
Since $f_Y(1,\theta)=Y(\theta)$, the map
$Y\mapsto(f_Y)_{\ell_t}$ is injective, and hence
$\dim Z_{t,y}=\dim\mathcal H_{d-1}=2d_y-1$.

Choose pairwise disjoint coordinate neighborhoods $O_y$ of the points
$y\in X$. In the coordinates centered at $y$, every element of
$Z_{t,y}$ is supported in $\ell_t A$. Since $\ell_t\to0$ at each $y$
and $X$ is finite, we have $Z_{t,y}\subset C_c^\infty(O_y)$ and
$Q_t$ is negative definite on $Z_{t,y}$ simultaneously for all
$y\in X$ when $t$ is sufficiently large. By disjointness of the
supports, for any $F_y\in Z_{t,y}$ with $F_y\ne0$ for at least one $y\in X$, there holds
\[
    Q_t\left(\sum_{y\in X}F_y\right)
    =\sum_{y\in X}Q_t(F_y)<0.
\]
Since $-\Delta_g+V_t$ has exactly $p$ negative eigenvalues, the
min--max principle gives
\[
    \sum_{y\in X}(2d_y-1)
    =\dim\left(\bigoplus_{y\in X}Z_{t,y}\right)\le p.
\]
The proof is complete.
\end{proof}

\section{Descent at a higher-order maximum}\label{sec:descent}

Let $u,v\in E$ be linearly independent and choose a global maximum
$x_{\ast}$ of $s=S(u,v)$. Normalize all stream functions at $x_{\ast}$, so that
$\max_M s=s(x_{\ast})=0$, and write $X=\{s=0\}$. We shall show that if
$m\ge p+2$ and the maximum at $x_{\ast}$ has order at least two, then a
change of potential reduces $p$ and increases $m-p$. We first consider the following continuous path of double Moutard transformation.

\subsection{Selection of the global maxima}
To preserve both $m$ and $p$ under the transformation, motivated by
Lemmas~\ref{lem:transform} and \ref{lem:negative-index}, we consider
the following continuous family of double Moutard transformations.

For $a,b\in E$ and $|t|$ small, taking $c=t^{-1}$ in
\eqref{eq:moutard}, up to a constant: $\tau_t=1+tS(a,b)$, and
\begin{align}\label{Sec 4 f-1}
    f_t=f+\frac{t}{\tau_t}\bigl(aS(b,f)-bS(a,f)\bigr),
    \qquad V_t=V-2\Delta_g\log|\tau_t|.
\end{align}
Here $\tau_t>0$ on $M$. Let $S_t$ denote the stream function for
$H_t:=-\Delta_g+V_t$, normalized at $x_{\ast}$. For $f,h\in E$, we claim that
\begin{equation}\label{eq:streams-change}
    S_t(f_t,h_t)=S(f,h)-\frac{t\left[S(a,f)S(b,h)-S(b,f)S(a,h) \right]}{\tau_t}.
\end{equation}
Indeed, \eqref{eq:moutard-dh} gives
\begin{align}\label{Sec 4 equ-a}
    \mathrm{d}(\tau_t f_t)
    =\tau_t\,\mathrm{d}f+
      t\bigl(S(b,f)\,\mathrm{d}a-S(a,f)\,\mathrm{d}b\bigr).
\end{align}
By the definition of $S_t(f_t,h_t)$, we obtain
\[
    \tau_t\,\mathrm{d}S_t(f_t,h_t)
    =\ast\bigl(f_t\,\mathrm{d}(\tau_t h_t)
              -h_t\,\mathrm{d}(\tau_t f_t)\bigr).
\]
Applying \eqref{Sec 4 equ-a} to $f$ and $h$, dividing by $\tau_t$,
and collecting the terms containing $\mathrm{d}a$ and $\mathrm{d}b$,
we obtain
\begin{align*}
    \mathrm{d}S_t(f_t,h_t)
    &=\ast\bigl(f_t\,\mathrm{d}h-h_t\,\mathrm{d}f\bigr)\\
    &\quad+\frac{t}{\tau_t}\ast\left[
      \bigl(f_tS(b,h)-h_tS(b,f)\bigr)\mathrm{d}a
   -
      \bigl(f_tS(a,h)-h_tS(a,f)\bigr)\mathrm{d}b\right].
\end{align*}
 The first
term in the expression for $\mathrm{d}S_t(f_t,h_t)$ expands as
\begin{align}\label{Sec 4 equ-b}
    f_t\,\mathrm{d}h-h_t\,\mathrm{d}f
    &=f\,\mathrm{d}h-h\,\mathrm{d}f\nonumber\\
    &\quad+\frac{t}{\tau_t}\Big[
       \bigl(aS(b,f)-bS(a,f)\bigr)\,\mathrm{d}h
       -\bigl(aS(b,h)-bS(a,h)\bigr)\,\mathrm{d}f
       \Big].
\end{align}
Write $D=S(a,f)S(b,h)-S(b,f)S(a,h)$. Substituting the expressions
for $f_t$ and $h_t$ into the second term
\begin{align}\label{Sec 4 equ-c}
    f_tS(b,h)-h_tS(b,f)
    =&fS(b,h)-hS(b,f)\nonumber\\
    +&\frac{t}{\tau_t}\Big[
       \bigl(aS(b,f)-bS(a,f)\bigr)S(b,h)
       -\bigl(aS(b,h)-bS(a,h)\bigr)S(b,f)\Big]\nonumber\\
    =&fS(b,h)-hS(b,f)-\frac{tb}{\tau_t}D.
\end{align}
Similarly, we also have
\begin{align}\label{Sec 4 equ-d}
    f_tS(a,h)-h_tS(a,f)
    =fS(a,h)-hS(a,f)-\frac{ta}{\tau_t}D.
\end{align}

Combining \eqref{Sec 4 equ-b}, \eqref{Sec 4 equ-c} and \eqref{Sec 4 equ-d}, we find
\begin{align*}
    \mathrm{d}S_t(f_t,h_t)
    &=\mathrm{d}S(f,h)\\
    &\quad+\frac{t}{\tau_t}\Big[
       S(b,f)\ast(a\,\mathrm{d}h-h\,\mathrm{d}a)
       -S(a,f)\ast(b\,\mathrm{d}h-h\,\mathrm{d}b)\\
    &\hspace{63pt}
       -S(b,h)\ast(a\,\mathrm{d}f-f\,\mathrm{d}a)
       +S(a,h)\ast(b\,\mathrm{d}f-f\,\mathrm{d}b)
       \Big]\\
    &\quad+\frac{t^2D}{\tau_t^2}
       \ast(a\,\mathrm{d}b-b\,\mathrm{d}a).
\end{align*}
By \eqref{Sec 3 def S}, the expression in square brackets is
\begin{align*}
    &S(b,f)\,\mathrm{d}S(a,h)-S(a,f)\,\mathrm{d}S(b,h)-S(b,h)\,\mathrm{d}S(a,f)+S(a,h)\,\mathrm{d}S(b,f)\\
    &=\mathrm{d}\bigl(S(b,f)S(a,h)\bigr)
      -\mathrm{d}\bigl(S(a,f)S(b,h)\bigr)
     =-\mathrm{d}D.
\end{align*}
Moreover,
$\ast(a\,\mathrm{d}b-b\,\mathrm{d}a)=\mathrm{d}S(a,b)$ and
$\mathrm{d}\tau_t=t\,\mathrm{d}S(a,b)$. It follows that
\begin{align*}
    \mathrm{d}S_t(f_t,h_t)
    &=\mathrm{d}S(f,h)-\frac{t}{\tau_t}\,\mathrm{d}D
      +\frac{t^2D}{\tau_t^2}\,\mathrm{d}S(a,b)\\
    &=\mathrm{d}S(f,h)-\frac{t}{\tau_t}\,\mathrm{d}D
      +\frac{tD}{\tau_t^2}\,\mathrm{d}\tau_t\\
    &=\mathrm{d}\!\left(S(f,h)-\frac{tD}{\tau_t}\right).
\end{align*}
Both sides of \eqref{eq:streams-change} vanish at $x_{\ast}$. Therefore, we complete the proof of \eqref{eq:streams-change}.

For the fixed pair $u,v$, write $s_t=S_t(u_t,v_t)$ and put
\[
    \mathcal D_z(a,b)=
    \begin{vmatrix}
        S(u,a)(z)&S(u,b)(z)\\
        S(v,a)(z)&S(v,b)(z)
    \end{vmatrix}.
\]
By skew-symmetry, \eqref{eq:streams-change} gives
\begin{equation}\label{eq:stream-first-variation}
    s_t(z)=s(z)-\frac{t}{\tau_t(z)}\mathcal D_z(a,b),
    \qquad
    \left.\partial_t s_t(z)\right|_{t=0}=-\mathcal D_z(a,b).
\end{equation}

For $Z\subseteq X$, set
\[
    E_Z=\{f\in E:S(u,f)(y)=S(v,f)(y)=0
                  \text{ for every }y\in Z\}.
\]
Since $S(u,v)=0$ on $Z$ and $S$ is skew-symmetric,
$\operatorname{span}\{u,v\}\subseteq E_Z$. As $u,v$ are linearly
independent, $\dim E_Z\ge2$. The number of independent conditions
imposed at the points of $Z$ is therefore $
    \operatorname{codim}_E E_Z=m-\dim E_Z\le m-2$.

\begin{lemma}\label{lem:selection}
After finitely many nonsingular transformations of the form
\eqref{eq:moutard}, the set $K$ of all global maxima of the
transformed stream satisfies $K\subseteq X$ and
\begin{equation}\label{eq:rows}
    x_{\ast}\in K,\qquad \operatorname{codim}_E E_K\le |K|-1,
\end{equation}
where $E$ and $E_K$ are taken for the transformed potential and pair.
The transformations preserve $p$, $m$, and the leading terms of
$u$ and $v$ at $x_{\ast}$.
\end{lemma}

\begin{proof}
\textbf{Step 1:}
If $\operatorname{codim}_E E_X\le |X|-1$, take $K=X$.
Suppose therefore that $\operatorname{codim}_E E_X>|X|-1$, and let
$\mathcal V$ be the subspace of pairs $(a,b)\in E\times E$ satisfying
\begin{align}\label{Sec 4 f-2}
    S(u,a)(y)=S(v,b)(y),\qquad
    S(u,b)(y)=-S(v,a)(y)\qquad y\in X.
\end{align}
By definition, $E_X\times E_X\subset\mathcal V$. Since all streams
vanish at $x_{\ast}$, the defining equations of $\mathcal V$ have total
rank at most $2(|X|-1)$. Thus
\[
    \operatorname{codim}_{E\times E}\mathcal V
    \le 2(|X|-1)
    <2\operatorname{codim}_E E_X
    =\operatorname{codim}_{E\times E}(E_X\times E_X).
\]
Hence, $E_X\times E_X\subsetneq\mathcal V $. Choose
$(a,b)\in\mathcal V\setminus(E_X\times E_X)$. Then
\begin{equation}\label{eq:separate}
    \mathcal D_y(a,b)=S(u,a)(y)^2+S(v,a)(y)^2
    \qquad y\in X.
\end{equation}
Set $X'=\{y\in X:\mathcal D_y(a,b)=0\}$. By \eqref{eq:separate}
and the defining equations of $\mathcal V$, a point $y\in X$
belongs to $X'$ if and only if all four mixed streams \eqref{Sec 4 f-2} vanish there.
In particular $x_{\ast}\in X'$. If $X'=X$, then $a,b\in E_X$, contrary
to our choice. Hence $x_{\ast}\in X'\subsetneq X$.

We claim that $X'$ is the full set of global maxima of $s_t$ for
all sufficiently small $t>0$. Consider the following three cases.
\begin{itemize}
\item If $y\in X\setminus X'$, then $\mathcal D_y(a,b)>0$.
By continuity, $\mathcal D_z(a,b)>0$ on a fixed neighborhood $O_y$
of $y$. Since $s\le0$ and $\tau_t>0$, \eqref{eq:stream-first-variation}
gives
\[
    s_t(z)\le -\frac{t}{\tau_t(z)}\mathcal D_z(a,b)<0,
    \qquad z\in O_y.
\]

\item If $y\in X'$, take isothermal coordinates centered at $y$,
and write $r=|z|$ and $d=d_y$. Since $u=O(r^d)$ and
$\mathrm{d}u=O(r^{d-1})$, we have
\begin{align}\label{Sec 4 euq-f}
    \mathrm{d}S(a,u)
    =\ast(a\,\mathrm{d}u-u\,\mathrm{d}a)=O(r^{d-1}).
\end{align}
As $S(a,u)(y)=0$, integration along radial segments gives
$S(a,u)(z)=O(r^d)$. The same estimate holds for the other three
mixed streams, so $|\mathcal D_z(a,b)|\le C r^{2d}$ near $y$.
By \eqref{eq:max-exp}, there is a fixed neighborhood $O_y$ on which
$s(z)\le -c r^{2d}$ for some $c>0$. Since $\tau_t\to1$ uniformly,
we may take $\tau_t\ge1/2$. Hence $s_t(y)=0$ and
\[
    s_t(z)\le(-c+2Ct)r^{2d}<0, 
    \qquad z\in O_y\setminus\{y\}
\]
for all sufficiently small $t>0$.

\item On the compact set $M\setminus\bigcup_{y\in X}O_y$, the
function $s$ is bounded above by a strictly negative constant.
Since $s_t\to s$ uniformly, $s_t<0$ there for small $t>0$.
\end{itemize}
As $X$ is finite by Lemma~\ref{lem:maxima}, a single choice of $t$
suffices for all these estimates. Thus $\max_M s_t=0$ and
$\{s_t=0\}=X'$.

At $x_{\ast}$, let $d=d_{x_{\ast}}$ and take $f=u$ or $v$. As above,
$S(a,f),S(b,f)=O(r^d)$. By \eqref{eq:moutard-dh},
\[
    \mathrm{d}\bigl(aS(b,f)-bS(a,f)\bigr)
    =S(b,f)\,\mathrm{d}a-S(a,f)\,\mathrm{d}b
      -f\,\mathrm{d}S(a,b)=O(r^d).
\]
Combining with  $aS(b,f)(x_{\ast}) -bS(a,f)(x_{\ast})=0$, so we have
$aS(b,f)-bS(a,f)= O(r^{d+1})$. The transformation formula \eqref{Sec 4 f-1} gives
$f_t-f=O(r^{d+1})$. Thus the degree-$d$ leading terms of $u$ and
$v$ at $x_{\ast}$ are preserved. By Lemmas~\ref{lem:transform}
and~\ref{lem:negative-index}, we can get  $m$ and $p$ are preserved
along this nonsingular path.

\textbf{Step 2:}
Starting with $(V_0,u_0,v_0,X_0)=(V,u,v,X)$, we iterate the
preceding step. At stage $j$, let $E_j=\ker(-\Delta_g+V_j)$, let
$S_j$ be the stream normalized at $x_{\ast}$, and let $X_j$ be the set of
global maxima of $S_j(u_j,v_j)$. For $Z\subseteq X_j$, write
\[
    E_Z^{(j)}
    =\{f\in E_j:S_j(u_j,f)(y)=S_j(v_j,f)(y)=0
                  \text{ for every }y\in Z\}.
\]
Thus $E_Z^{(0)}=E_Z$ for $Z\subseteq X$.

Stop if $\operatorname{codim}_{E_j}E_{X_j}^{(j)}\le |X_j|-1$.
Otherwise, choose $(a_j,b_j)\in E_j\times E_j$ by the preceding
construction and take $t_j>0$ sufficiently small. Denote the
transformed potential by $V_{j+1}$ and set
\[
\begin{aligned}
    T_jh
    =h+t_j\frac{a_jS_j(b_j,h)-b_jS_j(a_j,h)}
                  {1+t_jS_j(a_j,b_j)},\quad (u_{j+1},v_{j+1})=(T_j(u_j),T_j(v_j))
\end{aligned}
\]
By \eqref{eq:separate} and the one-step argument, the new global
maximum set is
\[
    X_{j+1}
    =\{y\in X_j:S_j(u_j,a_j)(y)=S_j(v_j,a_j)(y)=0\},
    \qquad x_{\ast}\in X_{j+1}\subsetneq X_j.
\]
Consequently, $1\le |X_j|\le |X|-j$.
If $j=|X|-1$, then $X_j=\{x_{\ast}\}$ and
$E_{X_j}^{(j)}=E_j$ by normalization at $x_{\ast}$. Hence the iteration
stops at some $N_0\le |X|-1$.
Take $K=X_{N_0}$. Then $x_\ast\in K\subseteq X$ and
$\operatorname{codim}_{E_{N_0}}E_K^{(N_0)}\le |K|-1$ by the stopping
criterion. If $N_0=0$, then $K=X$ and set $T=I$.
Otherwise $K\subsetneq X$ and let
$T=T_{N_0-1}\circ\cdots\circ T_0:E_0\to E_{N_0}$
be the composition of the kernel isomorphisms of Lemma~\ref{lem:transform}.
For $y\in K=X_{N_0}$, we have $y\in X_{j+1}$ for every $0\le j<{N_0}$. Then,
by the definition of $(a_j,b_j)$, we can get
\[
\begin{aligned}
    S_j(u_j,b_j)(y)&=S_j(v_j,a_j)(y)=0,\\
    S_j(v_j,b_j)(y)&=S_j(u_j,a_j)(y)=0.
\end{aligned}
\]
Applying
\eqref{eq:streams-change} successively gives
\begin{align*}
    S_{j+1}(u_{j+1},f_{j+1})(y)=S_{j}(u_{j},f_{j})(y),\qquad S_{j+1}(v_{j+1},f_{j+1})(y)=S_{j}(v_{j},f_{j})(y).
\end{align*}
Then for every $f\in E_0$, we obtain
\[
    S_{N_0}(u_{N_0},Tf)(y)=S(u,f)(y),\qquad
    S_{N_0}(v_{N_0},Tf)(y)=S(v,f)(y).
\]
Hence $f\in E_K^{(0)}$ if and only if $Tf\in E_K^{({N_0})}$.
Since $T$ is bijective, it restricts to an isomorphism from
$E_K^{(0)}$ onto $E_K^{({N_0})}$. Therefore,
\[
    \operatorname{codim}_{E_0}E_K^{(0)}
    =\operatorname{codim}_{E_{N_0}}E_K^{(N_0)}
    \le |K|-1.
\]The proof is complete.
\end{proof}

 Let $K$ be chosen as in Lemma~\ref{lem:selection}. We study the
singular limit of the double Moutard transformation and show that,
when $d_{x_{\ast}}\ge2$, the negative index $p$ decreases while $m-p$ increases.

\begin{proposition}\label{prop:descent}
Suppose that $m\ge p+2$ and that $S(u,v)$ has a global maximum $x_{\ast}$
of order $d_{x_{\ast}}\ge2$. There is a smooth real potential $\widehat V$
on $(M,g)$ such that $-\Delta_g+\widehat V$ has $\widehat p$
negative eigenvalues and nullity $\widehat m$ at zero, with
\[
    \widehat p<p,\qquad \widehat m-\widehat p\ge m-p+1.
\]
\end{proposition}

\begin{proof}
Apply Lemma~\ref{lem:selection} at $x_{\ast}$, keeping the same notation
for the transformed potential and pair. Now $K=\{s=0\}$,
$\dim E_K\ge m-|K|+1$, and $d_{x_{\ast}}\ge2$. Set
\[
    V_\infty=V-2\Delta_g\log|s|\qquad\text{on }M\setminus K.
\]
In the isothermal coordinates of \eqref{eq:max-exp} near $y\in K$, write
$|s|=c_y r^{2d_y}(1+a)$. Taylor's formula for the smooth function
$s$ bounds each derivative of its remainder of order $j\le2$ by
$O(r^{2d_y+1-j})$. Dividing by $r^{2d_y}$ gives
$|\nabla^j a|=O(r^{1-j})$ for $0\le j\le2$.
On the punctured disk $\Delta\log r=0$, so
$$\Delta_g\log|s|=e^{-2\omega_g}\Delta\log(1+a)=O(r^{-1}).$$
Consequently the potential $V_\infty$ extends as an
$L^q(M)$ function for $1<q<2$.
By the truncation estimate preceding Proposition~\ref{prop:transfer},  its quadratic form 
\begin{align*}
    Q_{\infty}(f)=\int_M\bigl(|\nabla_g f|_g^2+V_{\infty}|f|^2\bigr)\,\mathrm{d}A_g
\end{align*}
on $H^1(M)$ is closed and
bounded below. Write $\widehat p$ and $\widehat m$ for its negative
index and nullity.

For $f\in E_K$, by \eqref{Sec 4 euq-f}, we get the streams $S(u,f), S(v,f)= O(r^{d_y})$ near  $y\in K$. Moreover, using Lemma \ref{lem:maxima}, we also have $|s|\geq c r^{2d_y}$ near $y$. Hence
\[
    f_\infty=f+\frac{uS(v,f)-vS(u,f)}{s}
\]
is bounded near $K$. The local calculation in the proof of Lemma~\ref{lem:transform} gives
$\Delta_g f_\infty=V_\infty f_\infty$ away from $K$.
Since $f_\infty$ is bounded near $K$ and $V_\infty\in L^q(M)$,
we have $V_\infty f_\infty\in L^q(M)$.
For each $y\in K$, choose a small coordinate disk $D$ with
$D\cap K=\{y\}$, and let $h\in W^{2,q}(D)$ solve
$\Delta_g h=V_\infty f_\infty$ in $D$ with zero boundary values.
Since $q>1$, $h$ is bounded, so $f_\infty-h$ is a bounded
harmonic function on $D\setminus\{y\}$, hence the singularity
$y$ is removable. By the standard elliptic estimates in
\cite{Adams,GT},
$f_\infty\in W^{2,q}(M)\subset C^0(M)\cap H^1(M)$. Thus, we obtain the map
\begin{align}\label{Sec 4 map}
  f\to f_{\infty}: \qquad E_K\to \ker(-\Delta_g+ V_{\infty})
\end{align}
We claim that the kernel of the map \eqref{Sec 4 map} is exactly $\operatorname{span}\{u,v\}$, namely,   \begin{align}\label{Sec 4 equ-l}
    \{f\in E_K: f_{\infty}=0\}=\operatorname{span}\{u,v\}.
\end{align}
Since $s=0$ on $K$, both $u,v$ belong to $E_K$,
and substitution
gives $$u_\infty=u+\frac{uS(v,u)}{s}=0,\quad v_\infty=v-\frac{vS(u,v)}{s}=0,$$
where we have used $S(u,u)=0 , S(u,v)=-S(v,u)$.
Conversely, if $f_\infty=0$, then \eqref{eq:inject} with $\tau=s$ shows that
\begin{align*}
   s\,\mathrm{d}S(u,f)=S(u,f)\,\mathrm{d}s,
    \qquad
    s\,\mathrm{d}S(v,f)=S(v,f)\,\mathrm{d}s.  
\end{align*}
$S(u,f)/s$ and $S(v,f)/s$ are constant on the connected set
$M\setminus K$. Since $$f=\left(\frac{S(u,f) }{s}\right)v-\left(\frac{S(v,f) }{s}\right)u,$$  this linear identity extends to $K$ by continuity. Thus, we complete the proof of  \eqref{Sec 4 equ-l}. It follows that
\begin{equation}\label{eq:outer-kernel}
    \widehat m\ge\dim E_K-2\ge m-|K|-1.
\end{equation}

For $t>0$, let
\[
    V_t=V-2\Delta_g\log(t^{-1}-s)
\]
and denote the quadratic form of $-\Delta_g+V_t$ by $Q_t$.
As in Lemma~\ref{lem:bubbles}, this operator has exactly $p$
negative eigenvalues for every $t>0$. Moreover, since $s<0$ on
$M\setminus K$, we have $V_t\to V_\infty$ uniformly on compact
subsets of $M\setminus K$ as $t\to\infty$.

Choose disjoint isothermal coordinate disks centered at the points
of $K$. In each disk, with $r=|z|$, let
\begin{align*}
    \chi_\delta=\begin{cases}
        \displaystyle 0 \qquad &\mathrm{for}\quad r\le\delta^2\\
         \displaystyle \frac{1}{|\log\delta|}\log\left(\frac{r}{\delta^2}\right) \qquad &\mathrm{for}\quad \delta^2<r<\delta\\
         \displaystyle 1 \qquad &\mathrm{for}\quad \delta\leq r
    \end{cases}
\end{align*}
 Then
\[
    \int_M |\nabla_g\chi_\delta|_g^2\,\mathrm{d}A_g
    =\frac{2\pi |K|}{|\log\delta|}\longrightarrow 0.
\]
Thus $\chi_\delta\phi\to\phi$ in $H^1(M)$ for every
$\phi\in C^\infty(M)$. Smoothing $\chi_\delta\phi$ away from $K$
and using the density of $C^\infty(M)$ in $H^1(M)$ gives
\[
    \overline{C_c^\infty(M\setminus K)}^{\,H^1(M)}=H^1(M).
\]

If $\widehat p=0$, take $W=\{0\}$. Otherwise, let
$\phi_1,\ldots,\phi_{\widehat p}$ be an $L^2(M)$-orthonormal basis
of negative eigenfunctions for $-\Delta_g+V_\infty$, with
corresponding eigenvalues $\lambda_i<0$. Since
$Q_\infty(\phi_i,\phi_j)=\lambda_i\delta_{ij}$ and
$Q_\infty$ is continuous on $H^1(M)\times H^1(M)$ by
\eqref{eq:real-form}, we may choose
$\psi_i\in C_c^\infty(M\setminus K)$ sufficiently close to $\phi_i$
in $H^1(M)$ that the matrix $(Q_\infty(\psi_i,\psi_j))_{i,j}$
is negative definite. Hence
$W=\operatorname{span}\{\psi_1,\ldots,\psi_{\widehat p}\}$
has dimension $\widehat p$. On the compact set
$C=\bigcup_i\operatorname{supp}\psi_i\subset M\setminus K$,
we have $V_t\to V_\infty$ uniformly. Therefore,
\[
    Q_t(\psi_i,\psi_j)-Q_\infty(\psi_i,\psi_j)
    =\int_C(V_t-V_\infty)\psi_i\psi_j\,\mathrm{d}A_g
    \longrightarrow 0.
\]
It follows that $(Q_t(\psi_i,\psi_j))_{i,j}$ remains negative
definite for all sufficiently large $t$, as required.

Choose pairwise disjoint neighborhoods $O_y$ of the points $y\in K$
with $O_y\cap C=\varnothing$. By Lemma~\ref{lem:bubbles}, for all
sufficiently large $t$ and each $y\in K$, there is a subspace
$Z_{t,y}\subset C_c^\infty(O_y)$ of dimension $2d_y-1$ on which
$Q_t$ is negative definite. Since the supports are disjoint,
\[
    Q_t\!\left(w+\sum_{y\in K}z_y\right)
    =Q_t(w)+\sum_{y\in K}Q_t(z_y),
    \qquad w\in W,\quad z_y\in Z_{t,y}.
\]
Thus $Q_t$ is negative definite on
$W\oplus\bigoplus_{y\in K}Z_{t,y}$. The min--max principle gives
\begin{equation}\label{eq:outer-index}
    \widehat p+\sum_{y\in K}(2d_y-1)\le p.
\end{equation}

Since $d_y\ge1$ for every $y\in K$ and $d_{x_{\ast}}\ge2$,
\[
    \sum_{y\in K}(2d_y-1)\ge (|K|-1)+3=|K|+2.
\]
Therefore $\widehat p\le p-|K|-2<p$. Combining this estimate with
\eqref{eq:outer-kernel}, we obtain
\[
    \widehat m-\widehat p
    \ge (m-|K|-1)-(p-|K|-2)=m-p+1.
\]
In particular $\widehat m>0$, since $m\ge p+2$.
Proposition~\ref{prop:transfer}(2), applied to $V_\infty$ at zero,
now gives a smooth real potential $\widehat V$ on $(M,g)$ with
negative index $\widehat p$ and nullity $\widehat m$.
\end{proof}

For $H=-\Delta_g+V$ on $(M,g)$, write
$$p(H)=\#\{j:\lambda_j(H)<0\},\quad m(H)=\dim\ker H.$$ Define
\begin{align}
    \mathfrak C=\{H=-\Delta_g+V:V\in C^\infty(M,\mathbb R),
                      p(H)\ge3,\ m(H)\ge p(H)+2\}.
\end{align}
We shall prove in Section~\ref{sec:upper} that $\mathfrak C=\emptyset$.
\begin{cor}\label{cor:simple-global-maxima}
If $\mathfrak C\ne\emptyset$, let $H_{\min}\in\mathfrak C$ satisfy
$p(H_{\min})=\min_{H\in\mathfrak C}p(H)$. Then, for every linearly
independent pair $u,v\in\ker H_{\min}$, every global maximum of $S(u,v)$
has order one.
\end{cor}

\begin{proof}
For a smooth real Schr\"odinger operator $H$ on the two-sphere with
zero as an eigenvalue, write $p=p(H)$ and $m=m(H)$. In our indexing, Nadirashvili's estimate
\cite[Theorem~1(A), p.~225]{Nadirashvili} gives $m_k(H)\le 2k-1$.
For $k\ge3$, M. Hoffmann-Ostenhof, T. Hoffmann-Ostenhof and
Nadirashvili's result \cite[Theorem~1 and Remark~(ii)]{HHOsphere}
give $m_k(H)\le 2k-3$ for smooth real potentials.
Since zero first occurs at index $p+1$, we have $m=m_{p+1}(H)$ and hence
\begin{equation}\label{eq:low-level-bounds}
    p=0\ \Longrightarrow\ m=1,\qquad
    p=1\ \Longrightarrow\ m\le3,\qquad
    p=2\ \Longrightarrow\ m\le3.
\end{equation}
If $S(u,v)$ had a global maximum of order at least two,
Proposition~\ref{prop:descent}, applied to $H_{\min}$, would give
a smooth real potential $\widehat V$ on $(M,g)$ with
\[
    \widehat p<p(H_{\min}),\qquad
    \widehat m-\widehat p\ge m(H_{\min})-p(H_{\min})+1\ge3.
\]
By \eqref{eq:low-level-bounds}, $\widehat p\ge3$. Since
$\widehat m\ge\widehat p+2$, the operator $-\Delta_g+\widehat V$
belongs to $\mathfrak C$, contrary to the minimality of $p(H_{\min})$.
\end{proof}

	\section{The upper bound by an integral identity}\label{sec:upper}

Suppose, to the contrary, that $\mathfrak C\ne\emptyset$. Choose
$H_{\min}=-\Delta_g+V\in\mathfrak C$ with minimal negative
index, i.e., $p(H_{\min})=\min_{H\in \mathfrak C}p(H)$. Write
$E=\ker H_{\min}$, $p=p(H_{\min})$, and $m=\dim E$.
By Corollary~\ref{cor:simple-global-maxima}, every global maximum
of $S(u,v)$ has order one for every linearly independent pair
$u,v\in E$. Using this property, we now employ a global argument
to derive a contradiction for $m\geq p+2$.

\textbf{Differential forms:}
We briefly recall the notation for differential forms used below. A smooth
$j$-form $\beta_1$ assigns to each point an alternating $j$-linear form on the
tangent space. For smooth $j$-form $\beta_1$ and an arbitrary differential form $\gamma_1$, we denote their
wedge product by $\beta_1\wedge\gamma_1$. The exterior derivative is denoted by
$\mathrm d$ and satisfies
\[
\mathrm d(\beta_1\wedge\gamma_1)
=\mathrm d\beta_1\wedge\gamma_1
+(-1)^{j}\beta_1\wedge\mathrm d\gamma_1,
\qquad
\mathrm d^2=0.
\]
 For a vector field $X$, the contraction operator $\iota_X$ maps a $j$-form
$\beta_1$ to the $(j-1)$-form defined by
\[
(\iota_X\beta_1)(X_1,\ldots,X_{j-1})
=\beta_1(X,X_1,\ldots,X_{j-1}).
\]
In particular, for a $1$-form $\beta$, $\iota_X\beta=\beta(X)$. We will also
use the identity
\begin{align}\label{Sec 5 idot}
\iota_X(\beta_1\wedge\gamma_1)
=(\iota_X\beta_1)\wedge\gamma_1
+(-1)^{j}\beta_1\wedge(\iota_X\gamma_1).
\end{align}

\subsection{An integral identity on the space of orthonormal pairs}

Fix a real $L^2(M,\mathrm{d}A_g)$-orthonormal basis
$\phi_1,\ldots,\phi_m$ of $E$. We parametrize ordered orthonormal
pairs of eigenfunctions by
\begin{align}\label{Sec 5 Def fm}
    \mathcal F_m
    =\{q=(\mathbf a,\mathbf b)\in\mathbb R^m\times\mathbb R^m:
          |\mathbf a|=|\mathbf b|=1,\ \mathbf a\cdot\mathbf b=0\},
\end{align}
where the dot denotes the Euclidean inner product. For
$\mathbf a=(a_1,\ldots,a_m)$ and $\mathbf b=(b_1,\ldots,b_m)$,
the eigenfunctions corresponding to $q$ are
\begin{equation}\label{eq:frame-eigenfunctions}
    u_q(x)=\sum_{j=1}^m a_j\phi_j(x),\qquad
    v_q(x)=\sum_{j=1}^m b_j\phi_j(x),\qquad x\in M.
\end{equation}
Thus $q$ is the parameter, and $x$ is the point at which the
eigenfunctions are evaluated. Then $q\in\mathcal F_m$ if and only if
$\|u_q\|_{L^2(dA_g)}=\|v_q\|_{L^2(dA_g)}=1$
and $\int_M u_qv_q\,\mathrm{d}A_g=0$.

The definition of $\mathcal{F}_m$ involves three algebraic constraints; see \eqref{Sec 5 Def fm}. Hence $\mathcal F_m$ is a compact smooth manifold without boundary of dimension $2m-3$. Differentiating these three equations along a curve through $q$ gives the tangent space
\[
    \mathcal{T}_q\mathcal F_m=
    \left\{(\mathbf h,\mathbf k)\in\mathbb R^m\times\mathbb R^m:
    \begin{aligned}
       &\mathbf a\cdot\mathbf h=0,\qquad
        \mathbf b\cdot\mathbf k=0,\\
       &\mathbf h\cdot\mathbf b+\mathbf a\cdot\mathbf k=0
    \end{aligned}\right\}.
\]
The rotation of the coefficient pair is given by
\begin{equation}\label{eq:frame-rotation}
    R_\theta q
    =(\mathbf a\cos\theta-\mathbf b\sin\theta,\
      \mathbf a\sin\theta+\mathbf b\cos\theta).
\end{equation}
Its infinitesimal generator is
\begin{align}\label{Sec 5 def T}
    T_q
    =\left.\frac{\mathrm d}{\mathrm d\theta}R_\theta q
      \right|_{\theta=0}=(-\mathbf b,\mathbf a)\in \mathcal{T}_q\mathcal F_m.
\end{align}
Write $T$ for this vector field. By
\eqref{eq:frame-eigenfunctions}, the corresponding eigenfunctions
satisfy
\[
    u_{R_\theta q}=u_q\cos\theta-v_q\sin\theta,\qquad
    v_{R_\theta q}=u_q\sin\theta+v_q\cos\theta.
\]
then $S( u_{R_\theta q},v_{R_\theta q})=S(u_q,v_q) $.
 The coordinate functions $a_j,b_j$ define the forms
\begin{equation}\label{eq:frame-forms}
    \vartheta=\frac12\sum_{j=1}^m
       (a_j\,\mathrm{d}b_j-b_j\,\mathrm{d}a_j),
    \qquad
    \omega=\mathrm{d}\vartheta
          =\sum_{j=1}^m\mathrm{d}a_j
                         \wedge\mathrm{d}b_j .
\end{equation}
Both forms are invariant under $R_\theta$. At each $q$ we have
$\vartheta_q(T_q)=1$, and for every
$(\mathbf h,\mathbf k)\in \mathcal{T}_q\mathcal F_m$,
\[
    \omega_q(T_q,(\mathbf h,\mathbf k))
    =-\mathbf b\cdot\mathbf k-\mathbf a\cdot\mathbf h=0.
\]
Set $W_q=\operatorname{span}\{\mathbf a,\mathbf b\}^{\perp}
\subset\mathbb R^m$. For $(\mathbf h,\mathbf k)\in \mathcal{T}_q\mathcal F_m $, we have
\[
    \vartheta_q(\mathbf h,\mathbf k)
      =\tfrac12(\mathbf a\cdot\mathbf k-\mathbf b\cdot\mathbf h)
      =\mathbf a\cdot\mathbf k=-\mathbf b\cdot\mathbf h.
\]
Consequently $\ker\vartheta_q=W_q\times W_q\subset \mathcal{T}_q\mathcal F_m$. For every nonzero
$(\mathbf h,\mathbf k)\in W_q\times W_q$, the vector
$(-\mathbf k,\mathbf h)$ also belongs to $W_q\times W_q$, and
\[
    \omega_q\big((\mathbf h,\mathbf k),(-\mathbf k,\mathbf h)\big)
    =|\mathbf h|^2+|\mathbf k|^2>0.
\]
Thus $\omega_q$ is nondegenerate on the $(2m-4)$-dimensional
space $\ker\vartheta_q$. Since $\vartheta_q(T_q)=1$,
the form $\vartheta\wedge\omega^{m-2}$ is nowhere zero. Orienting
$\mathcal F_m$ by this form, we obtain
\begin{equation}\label{eq:frame-volume}
    \int_{\mathcal F_m}\vartheta\wedge\omega^{m-2}>0.
\end{equation}

The first lemma on the existence of zeros is a direct consequence of
Stokes' theorem and is independent of the preceding arguments in this paper.

\begin{lemma}\label{lem:frame-stokes}
Let $1\le p\le m-2$, and let $\ell_1,\ldots,\ell_p$ be positive
integers. Suppose that a smooth map
$F=(F_1,\ldots,F_p):\mathcal F_m\to\mathbb C^p$ satisfies
\[
    F_\nu(R_\theta q)=e^{i\ell_\nu\theta}F_\nu(q)
    \qquad
    \theta\in\mathbb R,\ 1\le\nu\le p.
\]
Then there exists $q_*\in\mathcal F_m$ such that
$F_1(q_*)=\cdots=F_p(q_*)=0$.
\end{lemma}

\begin{proof}
Suppose, to the contrary, that $F$ is nowhere zero. Then
$G=F/|F|$ is a smooth map from $\mathcal F_m$ to the unit sphere
$\mathbb S^{2p-1}\subset\mathbb C^p$. On this sphere, define the
real one-form
\[
    \gamma=
    \frac{\operatorname{Im}\sum_{\nu=1}^p
                 \overline{Z_\nu}\,\mathrm dZ_\nu}
         {\sum_{\nu=1}^p\ell_\nu|Z_\nu|^2}.
\]
Its denominator is positive since $\ell_\nu\ge1$ and
$\sum_{\nu=1}^p|Z_\nu|^2=1$. Set $\zeta=G^*\gamma$ and
$\kappa=\mathrm d\zeta$. Consider the $(2p-1)$-form
\[
    B=(\vartheta-\zeta)\wedge
      \sum_{j=0}^{p-1}\omega^{p-1-j}\wedge\kappa^j.
\]
We claim that
\begin{equation}\label{eq:global-primitive}
    \mathrm dB=\omega^p,\qquad \iota_TB=0.
\end{equation}
To achieve this , we need to prove it by two steps.

\textbf{Step 1}:
We first verify that
\begin{align}\label{Sec 5 lem-1 equ-a}
    \zeta(T)=1\qquad\mathrm{and}\qquad\iota_T\kappa=0.
\end{align}
Let $\rho_\theta$ denote the rotation on the target sphere given by
\[
    \rho_\theta(Z_1,\ldots,Z_p)
    =(e^{i\ell_1\theta}Z_1,\ldots,e^{i\ell_p\theta}Z_p).
\]
Since $|F(R_\theta q)|=|F(q)|$, the rotation law for $F$ gives
\begin{equation}\label{eq:normalized-frame-rotation}
    G\circ R_\theta=\rho_\theta\circ G.
\end{equation}
Differentiating at $\theta=0$ and using the definition of $T$
in \eqref{Sec 5 def T}, we obtain
\begin{equation}\label{eq:frame-generators}
    (\mathrm dG)_q(T_q)=W_{G(q)},\qquad
    W_Z:=\left.\frac{\mathrm d}{\mathrm d\theta}\right|_{\theta=0}
             \rho_\theta(Z)
        =(i\ell_1Z_1,\ldots,i\ell_pZ_p).
\end{equation}
By the definition of $\gamma$,
\[
    \gamma_Z(W_Z)
    =\frac{\operatorname{Im}\sum_{\nu=1}^p
                 \overline{Z_\nu}(i\ell_\nu Z_\nu)}
           {\sum_{\nu=1}^p\ell_\nu|Z_\nu|^2}
    =1.
\]
It follows from \eqref{eq:frame-generators} that
\[
    \zeta_q(T_q)
    =\gamma_{G(q)}\bigl((\mathrm dG)_q(T_q)\bigr)
    =\gamma_{G(q)}(W_{G(q)})=1.
\]
For each fixed $\theta$, the factors $\overline{Z_\nu}$ and
$\mathrm dZ_\nu$ acquire opposite phases under $\rho_\theta$.
Thus $\rho_\theta^*\gamma=\gamma$, and
\eqref{eq:normalized-frame-rotation} gives
\[
    R_\theta^*\zeta
    =(G\circ R_\theta)^*\gamma
    =G^*(\rho_\theta^*\gamma)=\zeta.
\]
Differentiating at $\theta=0$ yields $\mathcal L_T\zeta=0$, where $\mathcal{L}$ denotes the Lie derivative.
Since $\zeta(T)\equiv1$, Cartan's formula now gives
\[
    \iota_T\kappa
    =\mathcal L_T\zeta-\mathrm d\bigl(\zeta(T)\bigr)=0.
\]
Therefore, we complete the proof of \eqref{Sec 5 lem-1 equ-a}.

\textbf{Step 2}: To prove the first identity in \eqref{eq:global-primitive}, observe that
\[
    \kappa^p=G^*\bigl((\mathrm d\gamma)^p\bigr)=0,
\]
because $(\mathrm d\gamma)^p$ has degree $2p$ on
$\mathbb S^{2p-1}$. Since $\omega$ and $\kappa$ are closed
two-forms and $\mathrm d(\vartheta-\zeta)=\omega-\kappa$, we obtain
\[
    \mathrm dB
    =(\omega-\kappa)\wedge
      \sum_{j=0}^{p-1}\omega^{p-1-j}\wedge\kappa^j
    =\omega^p-\kappa^p=\omega^p.
\]
For the second identity, recall that $\vartheta(T)=1$ and
$\iota_T\omega=0$. Hence
\[
    \iota_T(\vartheta-\zeta)=\vartheta(T)-\zeta(T)=0,
    \qquad \iota_T\omega=\iota_T\kappa=0.
\]
The contraction rule \eqref{Sec 5 idot} therefore gives
$\iota_TB=0$, proving \eqref{eq:global-primitive}.

Finally, since $p\le m-2$, set
\[
    \Sigma=B\wedge\omega^{m-p-2},
    \qquad \omega^0=1.
\]
By \eqref{eq:global-primitive}, $\mathrm d\omega=0$, and
$\iota_T\omega=0$, we have
\[
    \mathrm d\Sigma=\omega^{m-2},
    \qquad \iota_T\Sigma=0.
\]
Consequently $\iota_T(\omega\wedge\Sigma)=0$. The form
$\omega\wedge\Sigma$ has degree $2m-3=\dim\mathcal F_m$.
At each $q\in\mathcal F_m$, complete the nonzero vector $T_q$
to a basis of $\mathcal{T}_q\mathcal F_m$. The form vanishes on this basis
and hence $\omega\wedge\Sigma=0$. It follows that
\[
    \mathrm d(\vartheta\wedge\Sigma)
    =\omega\wedge\Sigma-\vartheta\wedge\mathrm d\Sigma
    =-\vartheta\wedge\omega^{m-2}.
\]
Since $\mathcal F_m$ is compact and has no boundary, Stokes' formula
gives
\[
    0=\int_{\mathcal F_m}\mathrm d(\vartheta\wedge\Sigma)
     =-\int_{\mathcal F_m}\vartheta\wedge\omega^{m-2},
\]
contrary to \eqref{eq:frame-volume}.
\end{proof}

\subsection{Construction of map by the global maxima}

Normalize all stream functions on $E$ at the same point $x_{\ast}\in M$ as in
\eqref{Sec 3 def S}. For $q\in\mathcal F_m$, let $u_q,v_q$ be the
eigenfunctions defined in \eqref{eq:frame-eigenfunctions}, and set
\[
    w_q=u_q+iv_q,\qquad s_q=S(u_q,v_q),\qquad
    X(q)=\{x\in M:s_q(x)=\max_M s_q\}.
\]
By \eqref{eq:frame-rotation} and the bilinearity and skew-symmetry
of $S$,
\[
    w_{R_\theta q}=e^{i\theta}w_q,\qquad s_{R_\theta q}=s_q.
\]
Every global maximum has order one. Hence \eqref{eq:bubble-count}
gives
\begin{equation}\label{eq:number-maxima}
    1\le |X(q)|= \sum_{y\in X(q)}(2d_y-1) \le p.
\end{equation}
At each $y\in X(q)$, the maximum has order one. Thus, in the
isothermal coordinates $\Tilde{z}=(\Tilde{z}_1,\Tilde{z}_2)$ of Lemma~\ref{lem:maxima},
the linear term of $(u_q,v_q)$ is $P\cdot \Tilde{z}$, so its Jacobian matrix
at $y$ is precisely $P$, with $\det P>0$.

Identify $M$ with the Riemann sphere by an orientation-preserving
diffeomorphism. In a complex coordinate $z=z_1+iz_2$, write
$\partial_z=(\partial_{z_1}-i\partial_{z_2})/2$.
Since $w_q=u_q+iv_q$, the chain rule gives
\begin{equation}\label{eq:max-complex-derivative}
    |(w_q)_z(y)|^2-|(w_q)_{\bar z}(y)|^2
    =(\det P)\,
      \det\frac{\partial(\Tilde{z}_1,\Tilde{z}_2)}{\partial(z_1,z_2)}(y)>0,
\end{equation}
where the coordinate change preserves orientation. In particular,
$(w_q)_z(y)\ne0$ at every $y\in X(q)$.

For $t>0$, let
\[
    \mathrm d\mu_{q,t}
    =\frac{e^{t s_q}}{\int_M e^{t s_q}\,\mathrm dA_g}\,\mathrm dA_g .
\]
Then, we define the map
\begin{equation}\label{eq:maximum-moment-map}
    F_t=(F_{t,1},\ldots,F_{t,p}):
    \mathcal F_m\longrightarrow\mathbb C^p
\end{equation}
as follows. Fix an affine coordinate $z$ on the Riemann sphere and set
\[
    F_{t,j}(q)=\int_M z^{j-1}(w_q)_z\,\mathrm d\mu_{q,t},
    \qquad 1\le j\le3.
\]
For $4\le j\le p$, set
\begin{equation}\label{eq:vandermonde-maxima}
    F_{t,j}(q)
    =\int_{M^j}
        \Phi_{q,j}(x_1,\ldots,x_j)\,
       \mathrm{d}\mu_{q,t}^{\otimes j}(x_1,\cdots,x_j),
\end{equation}
where $z_a=z(x_a)$,
\begin{align*}
    \Phi_{q,j}(x_1,\ldots,x_j)&=\prod_{1\le a<b\le j}
       \Bigl[(z_a-z_b)^2(w_q)_z(x_a)(w_q)_z(x_b)\Bigr],\\
       \mathrm{d}\mu_{q,t}^{\otimes j}(x_1,\cdots,x_j)&=\mathrm d\mu_{q,t}(x_1)\cdots\mathrm d\mu_{q,t}(x_j).
\end{align*}

We verify that the integrands are globally defined. Under the fixed
identification of $M$ with the Riemann sphere, the two coordinates
$z$ and $\xi=1/z$ cover $M$. On their overlap, the chain rule gives
\[
    z^{j-1}(w_q)_z=-\xi^{3-j}(w_q)_\xi,\qquad 1\le j\le3.
\]
Since $3-j\ge0$, each of the first three integrands extends smoothly
across $\xi=0$.
For $j\ge4$, each pair factor in \eqref{eq:vandermonde-maxima}
is invariant under a M\"obius change of coordinate. Indeed, if
$\xi=f(z)=(Az+B)/(Cz+D)$ with $AD-BC\ne0$, then
\[
(f(z_a)-f(z_b))^2 =f'(z_a)f'(z_b)(z_a-z_b)^2,
    \qquad
    (w_q)_{\xi}(x_a)=\frac{(w_q)_z(x_a)}{f'(z_a)}.
\]
The factors $f'(z_a)$ and $f'(z_b)$ therefore cancel. Thus each pair
factor is globally smooth, and their product is smooth on $M^j$.

Finally, under $q\mapsto R_\theta q$, the measure $\mu_{q,t}$ is
unchanged and $(w_q)_z$ acquires the factor $e^{i\theta}$.
Since each of the $j(j-1)/2$ pair factors contains two derivatives,
\begin{equation}\label{eq:many-phase}
    F_{t,j}(R_\theta q)=
    \begin{cases}
        e^{i\theta}F_{t,j}(q), & 1\le j\le3,\\
        e^{ij(j-1)\theta}F_{t,j}(q), & 4\le j\le p.
    \end{cases}
\end{equation}

\begin{lemma}\label{lem:maximum-moments}
For all sufficiently large $t$, the map $F_t$ is nowhere zero
on $\mathcal F_m$.
\end{lemma}

\begin{proof}
Suppose that there exist $t_\ell\to\infty$ and $q_\ell\in\mathcal F_m$
such that $F_{t_\ell}(q_\ell)=0$. Write
$\mu_\ell=\mu_{q_\ell,t_\ell}$. By compactness of $\mathcal F_m$
and $M$, after passing to a subsequence we may assume that
$q_\ell\to q_0$ and $\mu_\ell\rightharpoonup\mu$ for a probability
measure $\mu$ on $M$. Since $E$ is finite dimensional and $S$ is
bilinear, we have
\[
    w_{q_\ell}\longrightarrow w_{q_0}\quad\text{in }C^1(M),
    \qquad
    s_{q_\ell}\longrightarrow s_{q_0}\quad\text{in }C^0(M).
\]

We first show that $\operatorname{supp}\mu\subseteq X(q_0)$.
Let $s_*=\max_M s_{q_0}$. Given $\delta>0$, choose an open ball
$B\subset\{s_{q_0}>s_*-\delta\}$ and set
$U_\delta=\{s_{q_0}<s_*-4\delta\}$. For all sufficiently large
$\ell$, uniform convergence gives $s_{q_\ell}>s_*-2\delta$ on $B$
and $s_{q_\ell}<s_*-3\delta$ on $U_\delta$. Hence
\[
    \mu_\ell(U_\delta)
    \le
    \frac{\operatorname{area}_g(M)e^{t_\ell(s_*-3\delta)}}
         {\operatorname{area}_g(B)e^{t_\ell(s_*-2\delta)}}
    =
    \frac{\operatorname{area}_g(M)}{\operatorname{area}_g(B)}
    e^{-\delta t_\ell}.
\]
Since $U_\delta$ is open, weak convergence gives
$\mu(U_\delta)\le\liminf_{\ell\to\infty}\mu_\ell(U_\delta)=0$.
As $\delta\downarrow0$, these sets increase to $M\setminus X(q_0)$.
Thus, by \eqref{eq:number-maxima}, there exist distinct points $y_1,\ldots,y_n\in X(q_0) $ such that
\[
    \mu=\sum_{a=1}^n c_a\delta_{y_a},
    \qquad
    c_a>0,\quad \sum_{a=1}^n c_a=1,\quad 1\le n\le p.
\]

 The affine coordinate $z$ used to define $F_t$
was fixed before these points were obtained, so one of $y_1,\ldots,y_n\in X(q_0) $ may be
the point at infinity of $z$. Choose $z_*\in\mathbb C$ such that
the point with $z$-coordinate $z_*$ lies outside
$\{y_1,\ldots,y_n\}$. Then $\xi=1/(z-z_*)$ is finite at every
$y_a$, and we write $\xi_a=\xi(y_a)$.

Suppose first that $n\le3$. Passing to the limit in the first three
components of $F_{t_\ell}(q_\ell)=0$ gives
\[
    \int_M \mathring{p}(z)(w_{q_0})_z\,\mathrm{d}\mu=0,
    \qquad \mathrm{for}\quad\deg \mathring{p}\le2.
\]
To use this identity for interpolation in $\xi$, observe that
for $j=0,1,2$,
\[
    \int_M \xi^j(w_{q_0})_\xi\,\mathrm{d}\mu
    =-\int_M (z-z_*)^{2-j}(w_{q_0})_z\,\mathrm{d}\mu=0.
\]
Here the first equality is the chain rule, and the second follows
because $(z-z_*)^{2-j}$ has degree at most two. Thus
\[
    \sum_{a=1}^n c_a \mathring{p}(\xi_a)(w_{q_0})_\xi(y_a)=0,
     \qquad \mathrm{for}\quad\deg \mathring{p}\le2.
\]
Choose
\begin{equation}\label{eq:three-values}
    \mathring{p}(\xi)=\prod_{a=2}^n(\xi-\xi_a),
\end{equation}
with $\mathring{p}=1$ if $n=1$. Since $\deg \mathring{p}=n-1\le2$, we obtain
\[
    0=c_1\mathring{p}(\xi_1)(w_{q_0})_\xi(y_1).
\]
This is impossible. Indeed, $c_1>0$, the distinctness of the
$\xi_a$ gives $\mathring{p}(\xi_1)\ne0$, and
\eqref{eq:max-complex-derivative} gives $(w_{q_0})_\xi(y_1)\ne0$.

Suppose now that $4\le n\le p$. We use the $n$-th component of
$F_{t_\ell}(q_\ell)$ in \eqref{eq:vandermonde-maxima}, where $n$
is the number of points in $\operatorname{supp}\mu$. Since
$(q,x_1,\ldots,x_n)\mapsto\Phi_{q,n}(x_1,\ldots,x_n)$ is continuous
on the compact set $\mathcal F_m\times M^n$, we have
\[
    \|\Phi_{q_\ell,n}-\Phi_{q_0,n}\|_{C^0(M^n)}\longrightarrow0.
\]
Moreover, $\mu_\ell^{\otimes n}\rightharpoonup\mu^{\otimes n}$,
since weak convergence holds on products of continuous functions
and their linear span is dense in $C(M^n)$. As $\mu_\ell$ is a
probability measure, we may pass to the limit and obtain
\begin{align*}
    0
    &=\lim_{\ell\to\infty}F_{t_\ell,n}(q_\ell)
      =\int_{M^n}\Phi_{q_0,n}\,\mathrm{d}\mu^{\otimes n}\\
    &=\sum_{i_1,\ldots,i_n=1}^n
      c_{i_1}\cdots c_{i_n}\,
      \Phi_{q_0,n}(y_{i_1},\ldots,y_{i_n}).
\end{align*}

To evaluate this sum, use the coordinate $\xi=1/(z-z_*)$ chosen
above, in which all $\xi_a=\xi(y_a)$ are finite. By the M\"obius
invariance of each pair factor in \eqref{eq:vandermonde-maxima},
\[
    \Phi_{q_0,n}(y_{i_1},\ldots,y_{i_n})
    =\prod_{1\le a<b\le n}
      \Bigl[(\xi_{i_a}-\xi_{i_b})^2
      (w_{q_0})_\xi(y_{i_a})(w_{q_0})_\xi(y_{i_b})\Bigr].
\]
If $i_a=i_b$ for some $a<b$, the factor
$(\xi_{i_a}-\xi_{i_b})^2$ vanishes. The remaining index tuples
are precisely the $n!$ permutations of $(1,\ldots,n)$. Since
$\Phi_{q_0,n}$ is symmetric, the sum reduces to
\[
    0=n!\left(\prod_{a=1}^n c_a\right)
        \left(\prod_{1\le a<b\le n}(\xi_a-\xi_b)^2\right)
        \prod_{a=1}^n\bigl((w_{q_0})_\xi(y_a)\bigr)^{n-1}.
\]
The right-hand side is nonzero, since each $c_a>0$, the $\xi_a$
are pairwise distinct, and \eqref{eq:max-complex-derivative} gives
$(w_{q_0})_\xi(y_a)\ne0$ at every $y_a\in X(q_0)$.
This contradiction completes the proof.
\end{proof}

\begin{proof}[\bf Proof of Theorem \ref{Thm 2}]
If $\mathfrak C\neq\emptyset$, then by Corollary~\ref{cor:simple-global-maxima}, there exists an operator $H_{\min}\in\mathfrak C$ such that
$
p(H_{\min})=\min_{H\in\mathfrak C}p(H).$
Lemma~\ref{lem:maximum-moments} gives a
nowhere-zero map $F_t$ for all sufficiently large $t$. On the other
hand, since $p\le m-2$, the rotation laws \eqref{eq:many-phase}
and Lemma~\ref{lem:frame-stokes} imply that $F_t$ has
a zero for every $t>0$. This contradiction proves that
$\mathfrak C=\emptyset$. 
Hence every smooth real spherical
Schr\"odinger operator $H$ with $p(H)\ge3$ satisfies $m(H)\le p(H)+1$.\end{proof}

\begin{proof}[\bf Proof of  Theorem \ref{thm:main} (Inequality part)]
Now fix $\lambda_k(\Omega)=\lambda_k^{(\alpha,\beta)}(\Omega)$ and let $p$ be the number of
eigenvalues strictly below it, counted with multiplicity.
Proposition~\ref{prop:transfer}(1) gives a spherical operator with
negative index $p$ and nullity $m_k(\Omega)$. The spherical bound
and \eqref{eq:low-level-bounds} give
\[
    m_k(\Omega)\le\max\{3,p+1\}.
\]
Since $p+1$ is the first index at which $\lambda_k(\Omega)$ occurs,
we have $p+1\le k$, see also \eqref{Intro f-1}. It follows that $m_k(\Omega)\le k$ for every
$k\ge3$. Especially for $k=2$, then $m_2(\Omega)\le3$.
\end{proof}

\section{A matrix on a planar graph}\label{sec:graph}
We now turn to the sharpness of the upper bound $m_k(\Omega)\leq k$.
To this end, for each fixed $k\geq3$, we need to construct a planar
domain $\Omega_k$ such that $m_k(\Omega_k)=k$. In the context of
spectral prescription, Colin de Verdi\`ere \cite{CdV,CdVArnold}
developed an approach based on graph theory, in which each vertex
is replaced by a room and each edge by a thin tube, so that the
spectrum of the resulting domain approximates the prescribed
spectral data of the graph. However, for our purposes, we need to
construct a special planar graph together with a corresponding
matrix having the desired spectral properties.

Let $G=(\mathcal V_G,\mathcal E)$ be a finite graph with vertex set
$\mathcal V_G=\{a_1,\ldots,a_N\}$. We write the edge $e=\{i,j\}\in\mathcal E$ if
$a_i$ and $a_j$ are joined by an edge. We denote by $\mathcal L(G)$
the space of real symmetric $N\times N$ matrices whose $(i,j)$-entry
vanishes whenever $i\ne j$ and $e=\{i,j\}\notin\mathcal E$.
The diagonal entries may be chosen arbitrarily.
A matrix $\mathcal{M}\in\mathcal L(G)$ is called \emph{well signed} if
$\mathcal{M}_{ij}<0$ for every edge $e=\{i,j\}\in\mathcal E$.

For example, if $G$ is the path
$\underset{a_1}{\bullet}\!-\!\underset{a_2}{\bullet}\!-\!
 \underset{a_3}{\bullet}$, then every matrix in $\mathcal L(G)$
has the form
\[
    \mathcal{M}=
    \begin{pmatrix}
        a & b & 0 \\
        b & c & d \\
        0 & d & f
    \end{pmatrix}.
\]
The two zero entries reflect the fact that $a_1$ and $a_3$ are
not adjacent. Such a matrix is well signed precisely when $b<0$
and $d<0$, while the diagonal entries $a,c,f$ may be chosen freely.
We write $n_-(\mathcal{M})$ for the number of negative eigenvalues of $\mathcal{M}$,
counted with multiplicity.

For each $k\ge3$, we shall construct a connected planar graph $G$
and a well-signed matrix $\mathcal{M}\in\mathcal L(G)$ with
$n_-(\mathcal{M})=k-1$ and $\dim\ker \mathcal{M}=k$. The construction also gives
control over perturbations of $\mathcal{M}$ on its zero eigenspace.
We write $E=\ker \mathcal{M}$ and denote by $\Sym(E)$ the space of symmetric
bilinear forms on $E$.

\smallskip
\noindent\textbf{The discrete Dirichlet problem.}
Let $G$ be a finite connected graph embedded in a closed disk.
Label its boundary vertices $a_1,\ldots,a_{k-1}$ in cyclic order and
its interior vertices $a_{k},\ldots,a_N$, as in the graph below.
Assign to each edge $e=\{i,j\}\in\mathcal E$ a positive weight
$c_{ij}=c_{ji}$.
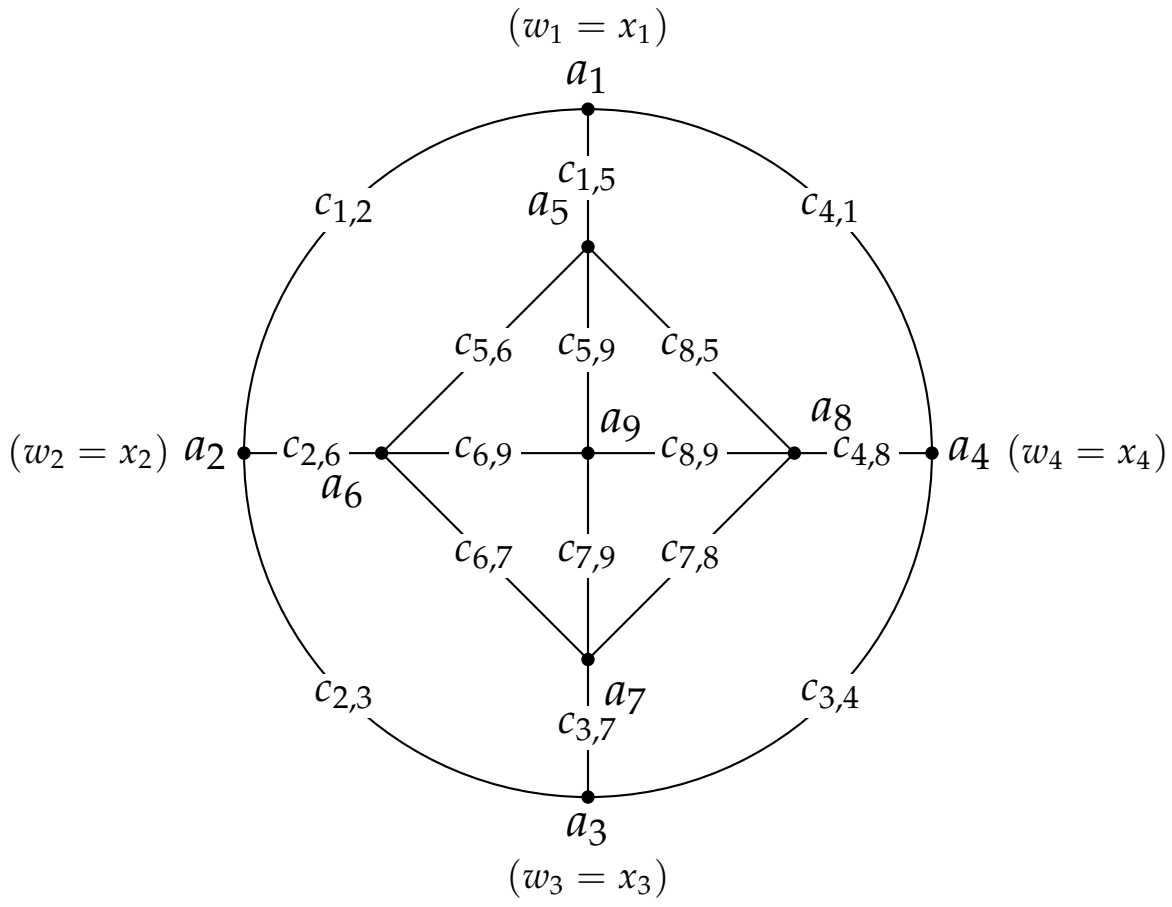
\[
\begin{tikzpicture}[
    scale=0.70, font=\small,
    weight/.style={fill=white,inner sep=1pt,font=\footnotesize}]
  \def\N{4}
  \def\outeredges{1} % 1: weighted outer cycle; 0: disk outline only
  \pgfmathtruncatemacro{\centerindex}{2*\N+1}

  \ifnum\outeredges=1\relax
    \draw (0,0) circle (3.25cm);
  \else
    \draw[gray,densely dotted] (0,0) circle (3.25cm);
  \fi
  \coordinate (a\centerindex) at (0,0);
  \fill (a\centerindex) circle (1.8pt);
  \node[above right,inner sep=2pt] at (a\centerindex) {$a_{\centerindex}$};

  \foreach \i in {1,...,\N} {
    \pgfmathsetmacro{\ang}{mod(90+360*(\i-1)/\N,360)}
    \pgfmathtruncatemacro{\k}{\N+\i}
    \coordinate (a\i) at (\ang:3.25cm);
    \fill (a\i) circle (1.8pt);
    \path (a\i) ++(\ang:4pt)
      node[anchor={\ang+180},inner sep=0pt] (name\i) {$a_{\i}$};
    \path (name\i.\ang) ++(\ang:3pt)
      node[anchor={\ang+180},inner sep=0pt,font=\scriptsize]
      {$(w_{\i}=x_{\i})$};
    \coordinate (a\k) at (\ang:1.95cm);
    \fill (a\k) circle (1.8pt);
    \pgfmathsetmacro{\labelang}{\ang+45}
    \path (a\k) ++(\labelang:8pt)
      node[anchor={\labelang+180},inner sep=0pt] {$a_{\k}$};
  }

  \foreach \i in {1,...,\N} {
    \pgfmathtruncatemacro{\j}{mod(\i,\N)+1}
    \pgfmathtruncatemacro{\k}{\N+\i}
    \pgfmathtruncatemacro{\lidx}{\N+\j}
    \ifnum\outeredges=1\relax
      \pgfmathsetmacro{\midang}{90+360*(\i-0.5)/\N}
      \node[weight] at (\midang:3.25cm) {$c_{\i,\j}$};
    \fi
    \draw (a\i) -- node[midway,weight] {$c_{\i,\k}$} (a\k);
    \draw (a\k) -- node[midway,weight] {$c_{\k,\lidx}$} (a\lidx);
    \draw (a\k) -- node[midway,weight] {$c_{\k,\centerindex}$}
                   (a\centerindex);
  }
\end{tikzpicture}
\]
For a function $w=(w_1,\ldots,w_N)$ on $\mathcal V_G$, let
\[
    \mathcal E_c(w)
    =\sum_{\{i,j\}\in\mathcal E}c_{ij}(w_i-w_j)^2
\]
be its Dirichlet energy. Given boundary values
$x=(x_1,\ldots,x_{k-1})^T\in\mathbb R^{k-1}$, the response matrix $R(c)$
is characterized by
\begin{equation}\label{Sec 6 equ-a}
    x^T R(c)x
    =
    \min\bigl\{\mathcal E_c(w):
        w_i=x_i \text{ for } 1\le i\le k-1\bigr\}.
\end{equation}
By the discrete Dirichlet principle
\cite[Theorem~0.41]{KLW}, the minimizer is the unique
\emph{harmonic extension} $w^x$ of $x$ to $\mathcal V_G$. More explicitly,
if $L_c$ denotes the weighted graph Laplacian
\[
    (L_cw)_i
    =\sum_{j:\,\{i,j\}\in\mathcal E}c_{ij}(w_i-w_j),
\]
then $w^x$ solves the discrete Dirichlet problem
\[
    \begin{cases}
        (L_cw^x)_i=0, & k\le i\le N,\\
        w_i^x=x_i,    & 1\le i\le k-1.
    \end{cases}
\]
Indeed, the interior equations follow by differentiating
$\mathcal E_c(w)$ with respect to the free variables
$w_{k},\ldots,w_N$.

The response matrix $R(c)$ is the discrete Dirichlet-to-Neumann
matrix: it maps the boundary values $x$ to the boundary currents
\[
    (R(c)x)_i=(L_cw^x)_i,\qquad 1\le i\le k-1.
\]
The harmonic extension of $\mathbf 1=(1,\ldots,1)^T\in\mathbb R^{k-1}$
is constant on $G$, so $R(c)\mathbf 1=0$. Moreover, by \eqref{Sec 6 equ-a} and the connectedness of $G$,
we have $x^T R(c)x=0$ if and only if $x_1=\cdots=x_{k-1}$.
Thus $R(c)$ restricts to a positive definite symmetric bilinear
form on $H_0=\mathbf 1^\perp$. Identifying $R(c)$ with this
restriction, we obtain the map
\begin{equation}\label{Sec 6 equ-c}
    R\colon (0,\infty)^{|\mathcal E|}
    \longrightarrow \Sym(H_0),\qquad c\longmapsto R(c),
\end{equation}
where $|\mathcal E|$ denotes the number of edges of $G$.

\begin{lemma}\label{lem:network-response}
For every $k\ge3$, there exist a finite connected graph
$G=(\mathcal V_G,\mathcal E)$ embedded in a closed disk, with boundary vertices
$a_1,\ldots,a_{k-1}$ in cyclic order, and positive edge weights
$c^0=(c_{ij}^0)_{\{i,j\}\in\mathcal E}$ such that the map
\[
    R:(0,\infty)^{|\mathcal E|}\longrightarrow\Sym(H_0)
\]
defined by \eqref{Sec 6 equ-c} has surjective differential at $c^0$.
Thus
\begin{equation}\label{eq:response-full}
    \bigl\{DR(c^0)[\dot c]:\dot c\in\R^{|\mathcal E|}\bigr\}
    =\Sym(H_0),
    \qquad H_0=\mathbf 1^\perp\subset\R^{k-1}.
\end{equation}
In other words, for every $Q\in\Sym(H_0)$ there exists an edge-weight
variation $\dot c=(\dot c_e)_{e\in\mathcal E}$ such that
\[
    \left.\frac{\mathrm d}{\mathrm dt}\right|_{t=0}
    x^T R(c^0+t\dot c)y=Q(x,y)
    \qquad\text{for all }x,y\in H_0.
\]
\end{lemma}

\begin{proof}
For $k=3$, take the graph consisting of a single edge $\{1,2\}$.
Writing $c=c_{12}$, we have $x^TR(c)x=c(x_1-x_2)^2$. Hence
\[
    R(c)=c\begin{pmatrix}
        1&-1\\
        -1&1
    \end{pmatrix},
    \qquad
    DR(c)[\dot c]=\dot c\begin{pmatrix}
        1&-1\\
        -1&1
    \end{pmatrix}.
\]
Since $H_0=\operatorname{span}\{(1,-1)^T\}$, the space $\Sym(H_0)$
is one-dimensional. Taking $\dot c=1$, we obtain
$
    (1,-1)\,DR(c)[1]\,(1,-1)^T=4.
$
Hence the restriction of $DR(c)[1]$ to $H_0\times H_0$ is nonzero
and thus spans $\Sym(H_0)$. This proves \eqref{eq:response-full} for $k=3$.

Suppose now that $k\ge4$. Let $G$ be the well-connected graphs used by Curtis
and Morrow in \cite[Section~4.6 and Section 2.5 P.20-P.24]{CM}. It is embedded in a disk,
has $k-1$ boundary vertices $a_1,\ldots,a_{k-1}$ in cyclic order,
and has exactly $(k-1)(k-2)/2$ edges.

We first identify the matrix $R(c)$ defined by \eqref{Sec 6 equ-a}
with the response matrix $\Lambda_c$ of Curtis and Morrow.
Fix positive edge weights $c$ and boundary data $x\in\R^{k-1}$.
As shown above, the minimizer $w^x$ in \eqref{Sec 6 equ-a} is
the unique harmonic extension of $x$, characterized by
\[
    (L_cw^x)_i=0 \quad k\le i\le N,
    \qquad
    w_i^x=x_i \quad 1\le i\le k-1.
\]
Thus $w^x$ is also the harmonic extension used in
\cite[Section~3.4]{CM}. By the definition of the boundary bilinear
form $L(c)$ in \cite[Section~4.6]{CM}, we have
\[
    x^T\Lambda_c x
    =L(c)(x,x)
    =\sum_{\{i,j\}\in\mathcal E}c_{ij}(w_i^x-w_j^x)^2
    =x^T R(c)x.
\]
Since both matrices are symmetric, polarization gives
$\Lambda_c=R(c)$ and hence
\[
    L(c)(x,y)=y^T R(c)x
    \qquad\text{for all }x,y\in\R^{k-1}.
\]

By \cite[Theorem~4.6]{CM}, $DL(c)$ is injective at every positive
edge-weight vector $c$. To apply this result on $H_0$, recall that
$R(c)\mathbf 1=0$ for every $c$. Differentiating this identity gives
$DR(c)[\dot c]\mathbf 1=0$ for every variation $\dot c$.
Suppose that $DR(c)[\dot c]$ vanishes on $H_0\times H_0$, i.e., $x^T DR(c)[\dot c]\,y=0$ for all $x,y\in H_0$ .
By symmetry and the decomposition $\R^{k-1}=H_0\oplus\R\mathbf 1$,
it vanishes on $\R^{k-1}\times\R^{k-1}$. Consequently,
\[
    DL(c)[\dot c](x,y)=y^TDR(c)[\dot c]x=0
    \qquad\text{for all }x,y\in\R^{k-1}.
\]
The injectivity of $DL(c)$ gives $\dot c=0$. Thus the map
\[
    \dot c\longmapsto DR(c)[\dot c]\big|_{H_0\times H_0}
\]
is injective. Since
\[
    |\mathcal E|=\frac{(k-1)(k-2)}2=\dim\Sym(H_0),
\]
it is also surjective. Taking $c_e^0=1$ for every edge
$e\in\mathcal E$ proves \eqref{eq:response-full}.
\end{proof}

% In the preamble: \usepackage{amsmath,amssymb,tikz}
\begin{lemma}\label{net:graph}
For every $k\ge3$, there exist a finite connected planar graph $G$
and a well-signed matrix $\mathcal{M}\in\mathcal L(G)$ such that
\begin{align}\label{Sec 6 Lem 2 formula}
    n_-(\mathcal{M})=k-1,\qquad \dim\ker \mathcal{M}=k.
\end{align}
Moreover, if $E=\ker \mathcal{M}$, the restriction map
$R_E: \mathcal L(G)\to\operatorname{Sym}(E)$ is surjective.
\end{lemma}

\begin{proof}
Take the weighted disk graph $G_0=(\mathcal V_{G_0},\mathcal E_0)$ of
Lemma~\ref{lem:network-response}, with weights $c_e=c_e^0$ and boundary vertices at $a_1,\ldots,a_{k-1}$ on the ball
in counterclockwise order. Write $i^+=i+1$ for $i<k-1$ and $(k-1)^+=1$. Let $F_i$ be the boundary face incident to the open arc from
$a_i$ to $a_{i^+}$. Then, we also  number the
remaining interior faces by $F_{k},\cdots $, we set $\mathcal{F}=\{F_1,\cdots,F_{k-1},F_{k},\cdots\}$.

Assign $x_i$ to each vertex $a_i$, $y_j$ to each face $F_j$. For an edge
$e$, we give a  direction from $a_i$ to $a_j$ with $i<j$, put
\[
    \delta_e x=x_i-x_j,\qquad
    \delta_e^\perp y=y_{r(e)}-y_{\ell(e)},
\]
where $r(e)$ and $\ell(e)$ index the faces to its right and left.
For illustrate more clearly, we give the picture of example for $k=5$:  the edge $a_1\to a_2$ has $F_1$ on its right
and $F_5$ on its left, so $\delta_e^\perp y=y_1-y_5$. 

% Figure 1. Oriented edges and faces.
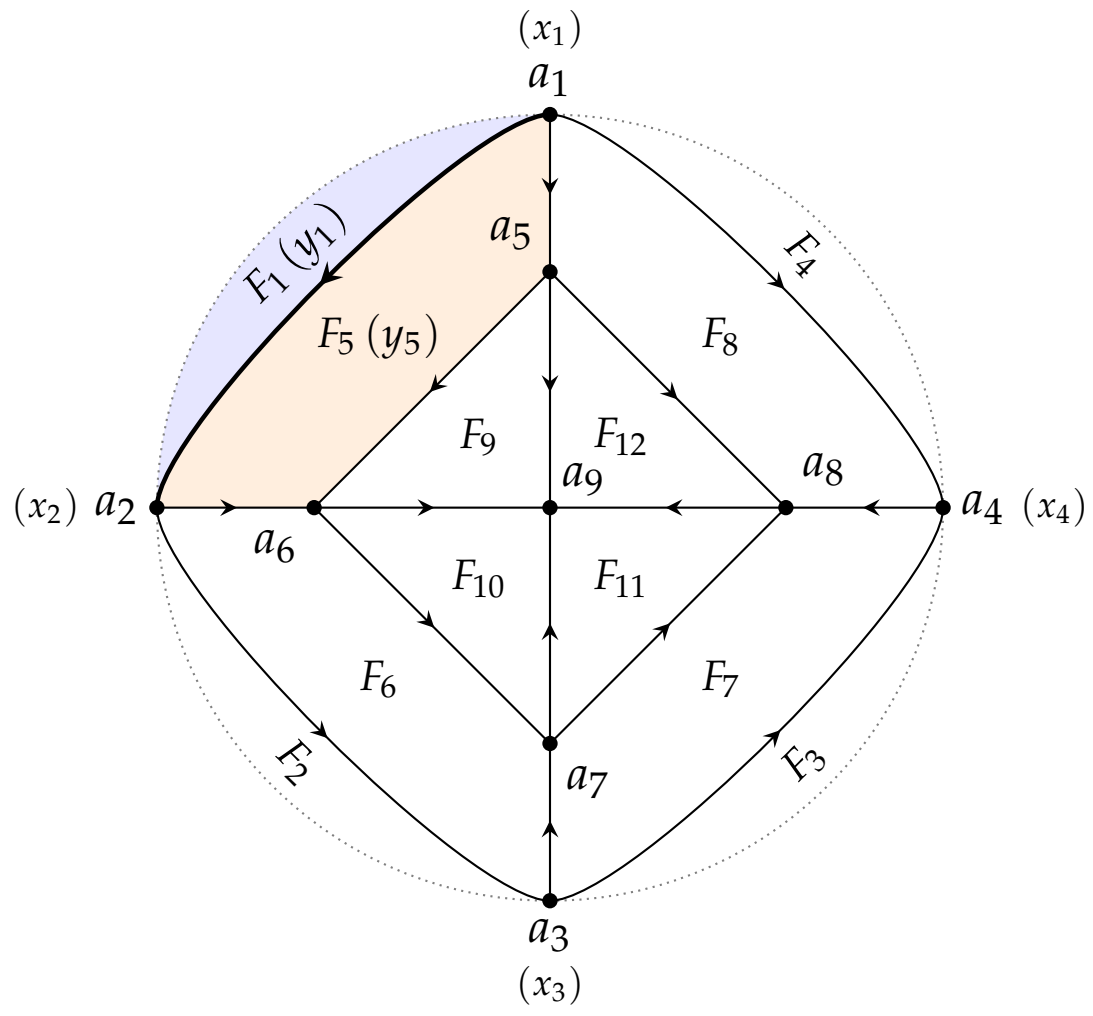
\usetikzlibrary{decorations.markings}
\begin{center}
\begin{tikzpicture}[scale=.80,font=\small,>=stealth,
    forward/.style={postaction={decorate},
      decoration={markings,mark=at position .51 with {\arrow{>}}}},
    backward/.style={postaction={decorate},
      decoration={markings,mark=at position .51 with {\arrow{<}}}},
    face/.style={font=\footnotesize,inner sep=1pt}]
  \def\innerR{1.95}
  \def\outerHandle{.66}
  \coordinate (a9) at (0,0);
  \foreach \i/\k/\t in {1/5/90,2/6/180,3/7/270,4/8/0}{
    \coordinate (a\i) at (\t:3.25);
    \coordinate (a\k) at (\t:\innerR);
    \path (a\i) ++({\t+90}:\outerHandle) coordinate (out\i);
    \path (a\i) ++({\t-90}:\outerHandle) coordinate (in\i);
  }
  \fill[blue!10] (a1) arc (90:180:3.25)
    .. controls (in2) and (out1) .. (a1) -- cycle;
  \fill[orange!13] (a1) .. controls (out1) and (in2) .. (a2)
    -- (a6) -- (a5) -- cycle;
  \draw[gray,densely dotted] (0,0) circle[radius=3.25];

  \foreach \i/\j/\k/\lidx in {1/2/5/6,2/3/6/7,3/4/7/8,4/1/8/5}{
    \ifnum\i=4\relax
      \tikzset{cycle arrow/.style={backward}}
    \else
      \tikzset{cycle arrow/.style={forward}}
    \fi
    \draw[cycle arrow,line width={\i==1 ? .85pt : .40pt}]
      (a\i) .. controls (out\i) and (in\j) .. (a\j);
    \draw[forward] (a\i) -- (a\k);
    \draw[cycle arrow] (a\k) -- (a\lidx);
    \draw[forward] (a\k) -- (a9);
  }
  \foreach \i/\k/\m/\t/\tilt in
    {1/5/9/135/45,2/6/10/225/-45,3/7/11/315/45,4/8/12/45/-45}{
    \ifnum\i=1\relax
      \node[face,rotate=\tilt] at (\t:2.97) {$F_1\,(y_1)$};
      \node[face] at (\t:2.00) {$F_5\,(y_5)$};
    \else
      \node[face,rotate=\tilt] at (\t:2.97) {$F_{\i}$};
      \node[face] at (\t:2.00) {$F_{\k}$};
    \fi
    \node[face] at (\t:.83) {$F_{\m}$};
  }
  \foreach \i in {1,...,9} \fill (a\i) circle[radius=1.8pt];
  \foreach \i/\t in {1/90,2/180,3/270,4/0}{
    \path (a\i) ++(\t:4pt)
      node[anchor={\t+180},inner sep=0pt] (name\i) {$a_{\i}$};
    \path (name\i.\t) ++(\t:3pt)
      node[anchor={\t+180},inner sep=0pt,font=\scriptsize] {$(x_{\i})$};
  }
  \foreach \k/\t in {5/135,6/225,7/315,8/45}
    \path (a\k) ++(\t:7pt)
      node[anchor={\t+180},inner sep=0pt] {$a_{\k}$};
  \node[above right,inner sep=2pt] at (a9) {$a_9$};
\end{tikzpicture}
\end{center}

\textbf{Step 1: Construction of $\mathcal{M}$} We introduce an
additional variable $z_i$ to each boundary, see Figure \ref{1}. For $u=(x,y,z)$, define the symmetric matrix $\mathcal{M}$ by
\begin{align*}
    u^T\mathcal{M}u=Q_0(x,y)-\sum_{i=1}^{k-1}(z_i+x_i+y_i)^2,
\end{align*}
where
\begin{align*}
    Q_0(x,y)=\sum_{e\in\mathcal E_0}
      \left(\sqrt{c_e}\,\delta_e x+
            \frac{\delta_e^\perp y}{\sqrt{c_e}}\right)^2 .
\end{align*}
The invertible linear change of variables
\[
    (x,y,z)\longmapsto(x,y,w),
    \qquad w_i=z_i+x_i+y_i,\quad 1\le i\le k-1,
\]
transforms the quadratic form into
$
    u^T\mathcal{M}u=Q_0(x,y)-|w|^2.
$
then we can see
$
    n_-(\mathcal{M})=k-1.
$
To identify the kernel, write $u=(x,y,z)$ and recall that
$\nabla_u(u^T\mathcal{M}u)=2\mathcal{M}u$. Since $Q_0$ is independent of $z$, then
\[
    0=\frac{\partial(u^T\mathcal{M}u)}{\partial z_i}
      =-2(z_i+x_i+y_i),
    \qquad 1\le i\le k-1.
\]
Hence $z_i=-x_i-y_i$. Under these identities, the derivatives
of the subtracted square terms vanish, so the remaining kernel
equations reduce to
\[
    \nabla Q_0(x,y)=0.
\]
As $Q_0$ is positive semidefinite, this is equivalent to
$Q_0(x,y)=0$. Consequently,
\begin{align*}
    \ker \mathcal{M}
    =&
    \bigl\{(x,y,z):
    Q_0(x,y)=0,\ 
    z_i=-x_i-y_i\ \text{for }1\le i\le k-1\bigr\}\\
    \simeq& \{(x,y):Q_0(x,y)=0\}:=\mathcal N
\end{align*}

Expanding the face differences and collecting the coefficient of
each $y_j$ gives
\[
    \sum_e\delta_e x\,\delta_e^\perp y
     =\sum_j y_j\left(\sum_{r(e)=j}\delta_e x
                      -\sum_{\ell(e)=j}\delta_e x\right).
\]
There are two cases to consider.
\begin{itemize}
\item For an interior face $F_j$, the coefficient of $y_j$
is zero, since the signed differences telescope along its
closed boundary walk, with each occurrence of an edge counted
with its induced sign. For instance, along
$\partial F_5$ in the first diagram, their sum is
\[
    (x_2-x_1)+(x_6-x_2)+(x_5-x_6)+(x_1-x_5)=0.
\]

\item For a boundary face $F_i$, there is only one boundary arc.
It runs from $a_i$ to $a_{i^+}$, so the graph edges run back from
$a_{i^+}$ to $a_i$. Thus their signed differences telescope to
$x_i-x_{i^+}$, which is the coefficient of $y_i$.
\end{itemize}
Therefore, we obtain
\begin{equation}\label{eq:net-summation}
    \sum_{e\in\mathcal E_0}\delta_e x\,\delta_e^\perp y
    =\sum_{i=1}^{k-1}(x_i-x_{i^+})y_i .
\end{equation}

By \eqref{eq:net-summation}, the terms $2x_i y_i$ in $Q_0$ cancel
upon subtraction of $(z_i+x_i+y_i)^2$. Expanding the squares gives
\[
\begin{aligned}
u^T\mathcal{M}u
={}&\sum_{e\in\mathcal E_0}
    \left(c_e(\delta_e x)^2
          +c_e^{-1}(\delta_e^\perp y)^2\right)\\
 &-\sum_{i=1}^{k-1}(x_i^2+y_i^2+z_i^2)
  -2\sum_{i=1}^{k-1}(x_i z_i+z_i y_i+y_i x_{i^+}).
\end{aligned}
\]

\textbf{Step 2: Construction of $G$.}
We take the disk embedding of $G_0$ to meet the boundary
circle only at $a_1,\ldots,a_{k-1}$.
Distinct open boundary edges correspond to distinct boundary
faces. Indeed, suppose that two distinct open boundary edges
$e$ and $e'$ belong to the same boundary face $F$.
Choose $P\in e$ and $Q\in e'$. Since $F$ is connected,
there exists a simple arc $\gamma$ joining $P$ to $Q$,
with its interior lying in $F$ and in the interior of the disk.
In particular, $\gamma\cap G_0=\varnothing$.
By the Jordan curve theorem, $\gamma$ separates the disk into
two connected components, each containing at least one
boundary vertex of $G_0$. Since $G_0$ is connected,
it must intersect $\gamma$, a contradiction
(see the figure below).

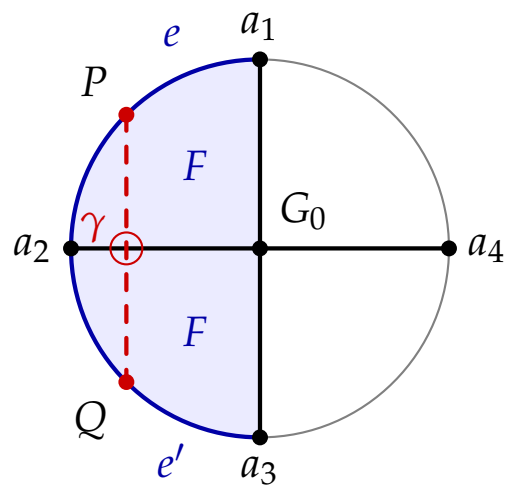
\begin{center}
\begin{tikzpicture}[scale=1,font=\footnotesize,
    line cap=round,line join=round]
    \def\rad{1.25}
    \coordinate (a1) at (90:\rad);
    \coordinate (a2) at (180:\rad);
    \coordinate (a3) at (270:\rad);
    \coordinate (a4) at (0:\rad);
    \coordinate (a5) at (0,0);
    \coordinate (P) at (135:\rad);
    \coordinate (Q) at (225:\rad);

    % The two regions assumed to belong to the same face F
    \fill[blue!8]
        (a5)--(a1) arc (90:270:\rad)--cycle;

    % The boundary circle
    \draw[gray] (0,0) circle[radius=\rad];

    % The two open boundary edges
    \draw[blue!65!black,thick]
        (a1) arc (90:180:\rad);
    \draw[blue!65!black,thick]
        (a2) arc (180:270:\rad);
    \node[blue!65!black] at (112.5:1.52) {$e$};
    \node[blue!65!black] at (247.5:1.52) {$e'$};

    % The connected graph G_0
    \foreach \i in {1,...,4}
        \draw[black,thick] (a\i)--(a5);
    \foreach \i in {1,...,5}
        \fill (a\i) circle[radius=1.5pt];

    % Boundary vertices
    \node[above] at (a1) {$a_1$};
    \node[left]  at (a2) {$a_2$};
    \node[below] at (a3) {$a_3$};
    \node[right] at (a4) {$a_4$};
    \node[above right] at (a5) {$G_0$};

    % Face labels under the contradiction assumption
    \node[blue!65!black] at (-.43,.55) {$F$};
    \node[blue!65!black] at (-.43,-.55) {$F$};

    % The hypothetical separating arc
    \draw[red!80!black,dashed,thick]
        (P)--(Q) node[pos=.42,left] {$\gamma$};
    \fill[red!80!black] (P) circle[radius=1.5pt];
    \fill[red!80!black] (Q) circle[radius=1.5pt];
    \node[above left] at (P) {$P$};
    \node[below left] at (Q) {$Q$};

    % The unavoidable intersection
    \draw[red!80!black] (P |- a5) circle[radius=3pt];
\end{tikzpicture}
\end{center}

We construct $G$ in three steps. In the figure, $k=5$, and
the orange circle is the boundary of the disk containing $G_0$.
The black graph is drawn inside this circle, its disk dual is
drawn outside it, and the joining paths lie on the circle.

\begin{enumerate}
\item \textbf{Draw the original graph in black.}
Place the boundary vertices $a_1,\ldots,a_{k-1}$ on the circle in
counterclockwise order and draw the remaining vertices and edges
of $G_0$ inside the disk. Give each black edge $e$ weight $-c_e$.
In the figure, $a_1,\ldots,a_4$ lie on the circle and
$a_5,\ldots,a_9$ lie inside it. The black edge $a_1a_2$
has weight $-c_{12}$.

\item \textbf{Draw the dual graph in blue.}
Introduce one vertex $b_j$ corresponding to each face $F_j$ of
$G_0$, including the boundary faces. For every original edge $e$
separating two distinct faces $F_j$ and $F_\ell$, join $b_j$ and
$b_\ell$ by a blue edge and give it weight $-c_e^{-1}$.
To draw this disk dual outside the circle, first draw it
separately in a copy of the disk, with $b_i$ on the open
boundary arc from $a_i$ to $a_{i^+}$ for $1\le i\le k-1$.
The usual face-by-face construction gives such an embedding,
with all other vertices and all edge interiors inside the disk.
After a small perturbation in the interior, we may assume that
the drawing avoids the center of the circle.
Inversion in the boundary circle then fixes the boundary
vertices $b_i$ and places all other dual vertices and all
dual edge interiors outside the disk, without introducing
any crossings.

In the figure, $b_1,\ldots,b_4$ lie on the circle,
$b_5,\ldots,b_8$ form the inner blue cycle, and
$b_9,\ldots,b_{12}$ form the outer blue cycle. Both cycles
are drawn using circular arcs. The radial blue edges form
the paths $b_i-b_{i+4}-b_{i+8}$ for $1\le i\le4$.
For example, the black edge $a_1a_2$ separates the two shaded
faces $F_1$ and $F_5$, so its dual is $b_1b_5$, with weight
$-c_{12}^{-1}$. Similarly, $b_5b_9$ is dual to $a_5a_6$ and
has weight $-c_{56}^{-1}$.

\item \textbf{Draw the joining paths in orange.}
For each boundary face $F_i$, draw an orange arc along the
boundary circle from $a_i$ to $a_{i^+}$ through $b_i$.
Insert a vertex $s_i$ on this arc between $a_i$ and $b_i$.
The resulting path is
\[
    a_i\!-\!s_i\!-\!b_i\!-\!a_{i^+},
\]
and each of its three edges has weight $-1$. For example, the
orange arc from $a_1$ to $a_2$ runs along the boundary arc
incident to the shaded face $F_1$ and consists of the three
edges $a_1s_1$, $s_1b_1$, and $b_1a_2$. The black and blue edge interiors lie on opposite sides of
the circle, while the orange paths lie on the circle.
By construction, all edges meet only at their prescribed
common endpoints. Hence this gives a planar embedding of $G$.

\end{enumerate}
% In the preamble: \usepackage{tikz,float}
% In the text: Figure~\ref{1}

\begin{figure}[H]
\centering
\renewcommand{\thefigure}{\ensuremath{1}}
\begin{tikzpicture}[scale=.76,font=\footnotesize,
    line cap=round,line join=round,
    primal/.style={draw=black,line width=.45pt},
    dual/.style={draw=blue!55!black,line width=.60pt},
    join edge/.style={draw=orange!80!black,line width=.75pt},
    bvertex/.style={
      circle,
      draw=blue!55!black,
      fill=white,
      inner sep=1.65pt
    },
    svertex/.style={
      rectangle,
      fill=orange!80!black,
      inner sep=1.5pt
    },
    face name/.style={
      font=\small\bfseries\boldmath,
      inner sep=.6pt
    },
    edge label/.style={
      font=\scriptsize,
      fill=white,
      inner sep=1.4pt
    }]

  \def\boundaryR{3.25}
  \def\innerR{1.52}
  \def\dualR{4.20}
  \def\outerDualR{5.20}
  \def\outerHandle{.40}

  % The original graph stays inside the circle.
  \coordinate (a9) at (0,0);
  \foreach \i/\j/\t in
    {1/5/90,2/6/180,3/7/270,4/8/0}{
    \coordinate (a\i) at (\t:\boundaryR);
    \coordinate (a\j) at (\t:\innerR);
    \path (a\i) ++({\t+90}:\outerHandle)
      coordinate (out\i);
    \path (a\i) ++({\t-90}:\outerHandle)
      coordinate (in\i);
  }

  % The shaded faces belong to the original disk graph G_0.
  \fill[blue!10]
    (a1)
    arc[start angle=90,end angle=180,radius=\boundaryR]
    .. controls (in2) and (out1) .. (a1) -- cycle;
  \fill[orange!13]
    (a1) .. controls (out1) and (in2) .. (a2)
    -- (a6) -- (a5) -- cycle;

  \foreach \i/\j/\r/\s in
    {1/2/5/6,2/3/6/7,3/4/7/8,4/1/8/5}{
    \draw[primal]
      (a\i) .. controls (out\i) and (in\j) .. (a\j);
    \draw[primal] (a\i) -- (a\r) -- (a9);
    \draw[primal] (a\r) -- (a\s);
  }
  \draw[primal,line width=.85pt]
    (a1) .. controls (out1) and (in2) .. (a2);

  % Boundary joining paths and the exterior disk dual.
  \foreach \i/\j/\r/\s/\t in
    {1/2/5/9/90,2/3/6/10/180,
     3/4/7/11/270,4/1/8/12/360}{
    \coordinate (s\i) at ({\t+22.5}:\boundaryR);
    \coordinate (b\i) at ({\t+45}:\boundaryR);
    \coordinate (b\r) at ({\t+45}:\dualR);
    \coordinate (b\s) at ({\t+45}:\outerDualR);
    \draw[join edge]
      (a\i)
      arc[start angle=\t,
          end angle={\t+90},
          radius=\boundaryR];
    \draw[dual] (b\i) -- (b\r) -- (b\s);
  }

  % Each circle consists of four edges,
  % meeting at the marked vertices.
  \draw[dual] (0,0) circle[radius=\dualR];
  \draw[dual] (0,0) circle[radius=\outerDualR];
  \draw[dual,line width=.95pt] (b1) -- (b5);

  % Vertex markers.
  \foreach \i in {1,...,9}
    \fill (a\i) circle[radius=1.7pt];
  \foreach \i in {1,...,12}
    \node[bvertex] at (b\i) {};
  \foreach \i in {1,...,4}
    \node[svertex] at (s\i) {};

  % Labels for the original vertices.
  \node[above=3pt] at (a1) {$a_1\,(x_1)$};
  \node at (-2.31,-.32) {$a_2\,(x_2)$};
  \node[below=3pt] at (a3) {$a_3$};
  \node[right=3pt] at (a4) {$a_4$};

  \foreach \i/\t in {5/12,6/102,7/192,8/282}
    \path (a\i) ++(\t:.30)
      node[inner sep=.5pt] {$a_{\i}$};
  \node[below right=2pt] at (a9) {$a_9$};

  % The four s_i labels lie between the orange and blue circles.
  \foreach \i/\ang/\tilt in
    {1/112.5/22.5,2/202.5/-67.5,3/292.5/22.5,4/22.5/-67.5}{
    \node[rotate=\tilt] at (\ang:3.72) {$s_{\i}\,(z_{\i})$};
  }

  % The four boundary b_i labels lie inside the orange circle.
  \node[rotate=45] at (135:2.91) {$b_1\,(y_1)$};
  \foreach \i/\ang in {2/225,3/315,4/45}
    \node at (\ang:2.87) {$b_{\i}$};

  % Exterior dual labels, offset from the circles and rays.
  \node[rotate=35] at (146:4.71) {$b_5\,(y_5)$};
  \node[anchor=north east] at (220:4.34) {$b_6$};
  \node[anchor=north west] at (320:4.34) {$b_7$};
  \node[anchor=south west] at (40:4.34) {$b_8$};
  \node[anchor=south east] at (135:5.38) {$b_9$};
  \node[anchor=north east] at (225:5.38) {$b_{10}$};
  \node[anchor=north west] at (315:5.38) {$b_{11}$};
  \node[anchor=south west] at (45:5.38) {$b_{12}$};

  % Face labels refer to the original disk graph G_0.
  \node[face name,text=blue!45!black,rotate=30]
    at (119:2.88) {$F_1$};
  \node[face name,text=orange!40!black]
    at (-1.23,1.48) {$F_5$};
  \node[font=\scriptsize,text=black!60]
    at (-.57,.28) {$F_9$};

  % Selected edge weights.
  \path (a1) .. controls (out1) and (in2) ..
    node[pos=.29,sloped,below=2pt,font=\scriptsize]
      {$-c_{12}$} (a2);
  \path (b1) --
    node[midway,sloped,above=3pt,
         font=\scriptsize,text=blue!55!black]
      {$-c_{12}^{-1}$} (b5);
  \path (b5) --
    node[midway,sloped,above=3pt,
         font=\scriptsize,text=blue!55!black]
      {$-c_{56}^{-1}$} (b9);
  \node[text=orange!80!black,font=\scriptsize]
    at (163:3.65) {$-1$};

  \begin{scope}[yshift=-5.80cm,font=\footnotesize]
    \draw[primal] (-4.0,0) -- (-3.65,0)
      node[right,black] {primal $(-c_e)$};
    \draw[dual] (-1.25,0) -- (-.90,0)
      node[right,black] {dual $(-c_e^{-1})$};
    \draw[join edge] (1.65,0) -- (2.0,0)
      node[right,black] {joining $(-1)$};
  \end{scope}
\end{tikzpicture}
\caption{}
\label{1}
\end{figure}

The disk dual is connected, since any two faces can be joined
by a path in the open disk avoiding the vertices of $G_0$
and crossing its edges transversely. Since $G_0$ is connected
and the orange paths join it to its disk dual, $G$ is connected.

We read the coefficients of $u^T\mathcal{M}u$ from the weighted
graph $G$ as follows. The vertices $a_i$, $b_j$, and $s_\ell$
carry the variables $x_i$, $y_j$, and $z_\ell$, respectively.
Parallel edges, if any, are merged into a single edge by
adding their weights.
For two distinct vertices, the coefficient of the product of
their variables is twice the corresponding edge weight and
is zero if there is no edge between them.

\begin{enumerate}
\item \emph{Terms $x_i x_j$.}
For $x_i^2$, sum $c_e$ over all black edges $e$ incident to $a_i$
and subtract $1$ if $a_i$ is a boundary vertex. For example, the
black edges incident to $a_1$ are $a_1a_2$, $a_1a_4$, and $a_1a_5$,
so the coefficient of $x_1^2$ is $c_{12}+c_{14}+c_{15}-1$.
For $i\ne j$, each black edge $a_i a_j$ contributes
$-2c_{ij}x_ix_j$. In particular, $a_1a_2$ contributes
$-2c_{12}x_1x_2$.

\item \emph{Terms $y_j y_\ell$.}
For $y_j^2$, sum $c_e^{-1}$ over the primal edges $e$ whose dual
edges are incident to $b_j$, and subtract $1$ if $F_j$ is a
boundary face. For $j\ne\ell$, each blue edge $b_jb_\ell$ dual
to $e$ contributes $-2c_e^{-1}y_jy_\ell$. For example, $b_1$
represents the shaded boundary face $F_1$, and the only blue edge
incident to it is $b_1b_5$. This edge is dual to $a_1a_2$ and joins
$b_1$ to the vertex $b_5$ representing $F_5$. Hence the coefficient of
$y_1^2$ is $c_{12}^{-1}-1$, and the coefficient of $y_1y_5$
is $-2c_{12}^{-1}$.

\item \emph{Terms $z_i z_j$.}
The diagonal correction gives the coefficient $-1$ for each
$z_i^2$, in particular for $z_1^2$. For $i\ne j$, the coefficient
of $z_i z_j$ is zero since $s_i$ and $s_j$ are not adjacent.

\item \emph{Terms $x_i y_j$.}
The orange edge $b_j a_{j^+}$ contributes $-2x_{j^+}y_j$
for each $1\le j\le k-1$. For example, $b_1a_2$ contributes
$-2x_2y_1$. The coefficient of $x_1y_1$ is zero since $s_1$
lies between $a_1$ and $b_1$ on the orange path and there is
no edge $a_1b_1$.

\item \emph{Terms $x_i z_j$.}
The orange edge $a_j s_j$ contributes $-2x_jz_j$ for each
$1\le j\le k-1$. For example, $a_1s_1$ contributes $-2x_1z_1$.

\item \emph{Terms $y_i z_j$.}
The orange edge $s_j b_j$ contributes $-2y_jz_j$ for each
$1\le j\le k-1$. For example, $s_1b_1$ contributes $-2y_1z_1$.
\end{enumerate}

Since all edge weights are strictly negative, we conclude that
$G$ is a finite connected planar graph and that
$\mathcal{M}\in\mathcal L(G)$ is well signed.

\smallskip
\textbf{Step 3: Proof of $\dim\ker \mathcal{M}=k$.} 
By the previous discussion, we know that
$\ker \mathcal{M}\simeq\mathcal N=\{(x,y):Q_0(x,y)=0\}$. We prove it by two steps:
\begin{itemize}
  
\item $\dim\mathcal N\geq k$. There is one variable $x_i$ for each vertex in $\mathcal V_{G_0}$ and one variable
$y_j$ for each face in $\mathcal F$. For each edge $e\in\mathcal E_0$,
the condition $Q_0(x,y)=0$ imposes one homogeneous linear equation
\[
    c_e\delta_e x+\delta_e^\perp y=0.
\]
Hence, the number of independent linear constraints is at most
$|\mathcal E_0|$, and therefore
\[
    \dim\mathcal N\ge |\mathcal V_{G_0}|+|\mathcal F|-|\mathcal E_0|.
\]

To compute the right-hand side, adjoin the $k-1$ boundary arcs joining
consecutive boundary vertices. The resulting decomposition of the disk
has $|\mathcal V_{G_0}|$ vertices, $|\mathcal E_0|+k-1$ edges, and $|\mathcal F|$
faces. Euler's formula for the disk gives
\[
    |\mathcal V_{G_0}|-\bigl(|\mathcal E_0|+k-1\bigr)+|\mathcal F|=1.
\]
It follows that
\[
    \dim\mathcal N\ge |\mathcal V_{G_0}|+|\mathcal F|-|\mathcal E_0|=k.
\]
  \item $\dim\mathcal N\leq k$. Consider the linear map
\begin{align}\label{Sec 6 lem 2 equ-a}
    \Phi:\mathcal N\longrightarrow\mathbb R^{k-1}\times\mathbb R,\qquad
    \Phi(x,y)=(x_\partial,y_1),\qquad x_\partial=(x_1,\ldots,x_{k-1})^T. 
\end{align}
To prove that $\Phi$ is injective, suppose that $(x,y)\in\mathcal N$ satisfies
$x_\partial=0$ and $y_1=0$. Since $x_\partial=0$, the identity
\eqref{eq:net-summation} from Step~2 gives
\[
    \sum_{e\in\mathcal E_0}\delta_e x\,\delta_e^\perp y
    =\sum_{i=1}^{k-1}(x_i-x_{i^+})y_i=0.
\]
Thus the mixed term in the expansion of $Q_0(x,y)$ vanishes.
Since $(x,y)\in\mathcal N$, we obtain
\[
    0=Q_0(x,y)
    =\sum_{e\in\mathcal E_0}c_e(\delta_e x)^2
     +\sum_{e\in\mathcal E_0}c_e^{-1}(\delta_e^\perp y)^2.
\]
As $c_e>0$ for every edge $e$, each summand vanishes. Hence $x$
takes the same value at the two endpoints of every black edge.
Since $G_0$ is connected and $x$ vanishes at the boundary vertices,
we have $x=0$. Similarly, $y$ takes the same value at the two
endpoints of every blue edge. Since the disk dual is connected,
all face values equal $y_1=0$. Hence $y=0$, proving that $\Phi$
is injective. In fact, $\Phi$ is an isomorphism from $\mathcal N\to \R^{k-1}\times\R$.
\end{itemize}
\smallskip
\smallskip
\smallskip
\noindent\textbf{Step 4. Surjectivity of the restriction map.}
By \textbf{Step~3} and \eqref{Sec 6 lem 2 equ-a}, the boundary values
$x_\partial=(x_1,\ldots,x_{k-1})^T$ together with $y_1$ give linear
coordinates $(x_{\partial},y_1)$ on $E=\ker \mathcal{M}$. We say that a quadratic form $q$ on
$E$ is realized by $D\in\mathcal L(G)$ if
\[
    u^TDu=q(u)\qquad\text{for every }u\in E.
\]
Since $\mathcal L(G)$ is a vector space, it suffices to realize
the monomials $x_ix_j$ for $1\le i\le j\le k-1$, $x_jy_1$ for
$1\le j\le k-1$, and $y_1^2$.

\smallskip
\noindent\emph{(i) The monomials $x_ix_j$.}
For $1\le i\le k-1$, let $D\in\mathcal L(G)$ be the matrix whose
only nonzero entry is $D_{a_i,a_i}=1$. Then
$u^TDu=x_i^2$ for every $u=(x,y,z)\in E$. 

Fix  $1\le i<j\le k-1$, we need to achieve $x_ix_j$. By Lemma~\ref{lem:network-response}, there
exists a variation $\dot c$ of the black edge weights such that
\begin{align}\label{Sec 6 equ-d}
    \xi^T DR(c^0)[\dot c]\,\xi=(\xi_i-\xi_j)^2
    \qquad\text{for all }\xi\in H_0.
\end{align}
Since $R(c)\mathbf 1=0$ for every $c$, both sides are unchanged
by adding a constant to $\xi$. Thus the identity holds on $\R^{k-1}$.

Now let $u=(x,y,z)\in E$. For every $w$ with
$w_\partial=x_\partial$, \eqref{eq:net-summation} gives
\[
    \mathcal E_{c^0}(w)-\mathcal E_{c^0}(x)
    =Q_0(w,y)-Q_0(x,y)=Q_0(w,y)\ge0.
\]
Thus $x$ attains the minimum in \eqref{Sec 6 equ-a}. By the definition of $R(c)$ in \eqref{Sec 6 equ-a}, we obtain
\[
    x_\partial^T R(c^0+t\dot c)x_\partial
    -\sum_{e\in\mathcal E_0}(c_e^0+t\dot c_e)(\delta_e x)^2\leq 0
\]
for all sufficiently small $|t|$, with equality at $t=0$.
The left-hand side has a maximum at $t=0$. Its derivative vanishes there, so \eqref{Sec 6 equ-d} gives
\[
    (x_i-x_j)^2
    =x_\partial^T DR(c^0)[\dot c]\,x_\partial
    =\sum_{e\in\mathcal E_0}\dot c_e(\delta_e x)^2.
\]
Each term on the right is supported on the endpoints of a black
edge, so the sum is realized by a matrix in $\mathcal L(G)$.
Since $x_i^2$ and $x_j^2$ are already realized, the identity
$2x_ix_j=x_i^2+x_j^2-(x_i-x_j)^2$ gives $x_ix_j$ as well.

\smallskip
\noindent\emph{(ii) The monomials $x_jy_1$.}
For each $1\le i\le k-1$, let $D^{(i)}$ be the symmetric matrix
indexed by the vertices of $G$ whose only nonzero entries are
\[
    D^{(i)}_{a_{i^+},b_i}
    =D^{(i)}_{b_i,a_{i^+}}=\frac12.
\]
Since $a_{i^+}$ and $b_i$ are joined by an orange edge,
$D^{(i)}\in\mathcal L(G)$. For every $u=(x,y,z)\in E$, we have
\[
    u^T D^{(i)}u
    =\frac12 x_{i^+}y_i+\frac12 y_i x_{i^+}
    =x_{i^+}y_i.
\]
To express this form in the
coordinates $(x_\partial,y_1)$, observe that there exists cofficents $a_{ij}$ such that $y_i=y_1+\sum_{j=1}^{k-1}a_{ij}x_j $, then
\[
    x_{i^+}y_1=x_{i^+}y_i-\sum_{j=1}^{k-1}a_{ij}x_{i^{+}} x_j .
\]
The second term on the right is a quadratic form in $x_\partial$
and is therefore realized by part~(i). Thus $x_{i^+}y_1$ is
also realized. Since $i^+$ runs through all boundary indices,
this gives every $x_jy_1$.

\smallskip
\noindent\emph{(iii) The monomial $y_1^2$.}
Taking the diagonal entry at $b_1$ equal to $1$ and all other
entries equal to zero realizes $y_1^2$. 
We have thus realized every quadratic form on $E$. By polarization,
the restriction map
\[
    \mathcal L(G)\longrightarrow\Sym(E),
    \qquad D\longmapsto D|_{E\times E},
\]
is surjective.
\end{proof}

\section{The sharpness of the multiplicity on the planar domains} \label{sec:domain}

In this section, for each boundary condition, we construct a smooth
bounded connected planar domain whose Laplacian has an eigenvalue
of multiplicity $\dim\ker\mathcal M$, with exactly $n_-(\mathcal M)$
eigenvalues below it, counted with multiplicity. The matrices constructed
in Section~\ref{sec:graph} satisfy $\dim\ker\mathcal M=k$ and
$n_-(\mathcal M)=k-1$, so the resulting eigenvalue first occurs
at index $k$. In the Dirichlet case, the domain is obtained by
removing finitely many circular holes from a disk.
In Subsection~\ref{sec:robin-realization}, using the same matrix family,
we construct domains consisting of rooms connected by thin tubes
for both the Neumann problem and the Robin problem with any fixed
parameter $\sigma\in\mathbb R$.

Let $\mathcal{M}$ be a well-signed matrix on a finite connected planar graph
$G$ with $N$ vertices. Assume that $E=\ker \mathcal{M}\ne\{0\},\ m:=\dim E$ and that
the restriction map
\[
    R_E:\mathcal L(G)\longrightarrow\Sym(E),
    \qquad R_E(D)=D|_{E\times E},
\]
is surjective. Choose a basis $q_1,\ldots,q_{d_m}\ (d_m=\dim\Sym(E)=\frac{m(m+1)}{2})$ of $\Sym(E)$ and,
for each $j$, a matrix $D_j\in\mathcal L(G)$ such that $R_E(D_j)=q_j$.
Extending $L(q_j)=D_j$ linearly gives a right inverse
$L:\Sym(E)\to\mathcal L(G)$ of $R_E$, which means that $R_E\circ L=I_{\Sym(E)}$ but $L\circ R_E$ may not be identity map.

We introduce a parameter $t\in\Sym(E)$ and set $\mathcal{M}(t)=\mathcal{M}+L(t)$. Since $E=\ker \mathcal{M}$ and $R_E(L(t))=t$, we have
\[
    u^T \mathcal{M}(t)v = u^T L(t)v = t(u,v)
    \qquad\text{for all }u,v\in E.
\]
On a sufficiently small closed ball
$\mathcal B=\overline{B_{\tilde r}}\subset\Sym(E)$ centered at zero, the matrices $\mathcal M(t)$
remain well signed. Write
\begin{equation}\label{net:family}
 \begin{aligned}
 x^T\mathcal{M}(t)x
 &=\frac{1}{2}\sum_{i,j}c_{ij}(t)(x_i-x_j)^2
   +\sum_{i=1}^N b_i(t)x_i^2,\\
 c_{ij}(t)&=-\mathcal{M}(t)_{ij},
 \qquad b_i(t)=\sum_{j=1}^N \mathcal{M}(t)_{ij}.
 \end{aligned}
\end{equation}
Since $\mathcal M$ is well signed, the weight corresponding to each edge
$e=\{i,j\}$ satisfies $c_{ij}>0$.  Fix \begin{align}\label{Sec 7 equ-a}
    \kappa\ge 1+\max_{t\in\mathcal B,\,1\le i\le N}|b_i(t)|,
\end{align}
then $\mathcal{M}(t)+\kappa I\ge I$ for every $t\in\mathcal B$.

\subsection{From the graph to a smooth potential}

We follow the construction with fixed rooms and tubes in
\cite[Sections~3--4]{CLZ}, adding a zeroth-order term to
represent the vertex term in \eqref{net:family}.
Let $\mathbb B\subset\R^2$ be the unit disk. A planar embedding
of $G$ gives disjoint Lipschitz rooms $R_i$, all of area $v>0$,
joined by fixed Lipschitz tubes $T_e$, $e=\{i,j\}$, with disjoint
end interfaces $\Gamma_{ei}$ and $\Gamma_{ej}$. We choose the
configuration so that
\[
 D=\operatorname{int}\left(\bigcup_i\overline{R_i}
       \cup\bigcup_e\overline{T_e}\right)\Subset\mathbb B
\]
is connected and Lipschitz.
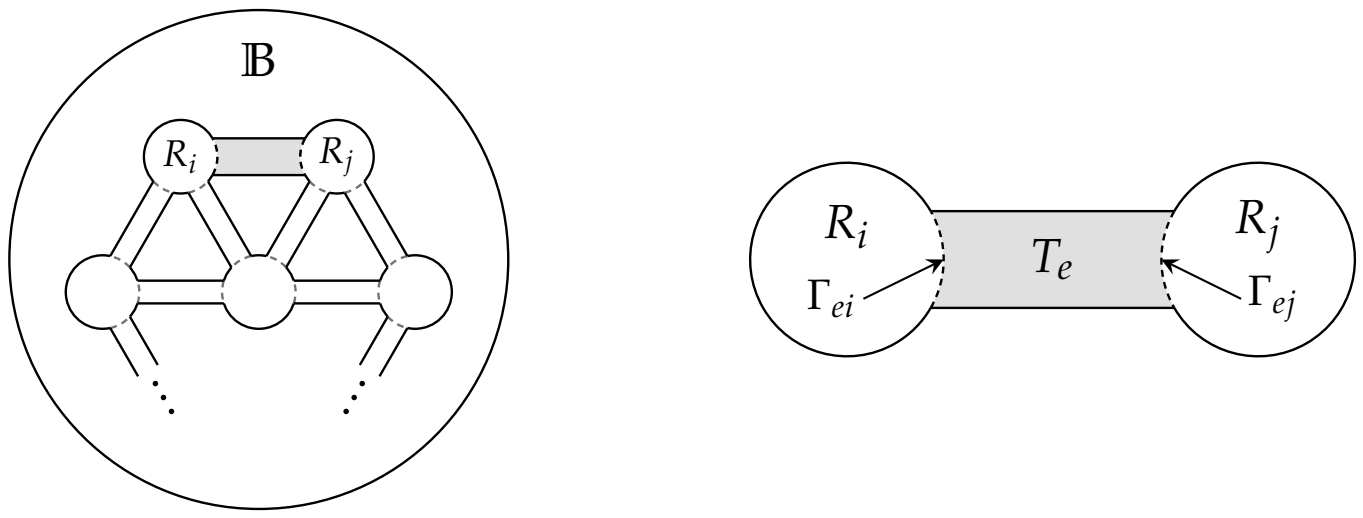
\begin{figure}[htbp]
\centering
\begin{tikzpicture}[
    x=1cm, y=1cm, font=\small, line width=.45pt,
    interface/.style={dash pattern=on 1.3pt off 1pt}
]
    % For one horizontal edge the disk centers are (-a,0) and (a,0).
    \def\a{.68}
    \def\r{.32}
    \pgfmathsetmacro{\w}{\r*sin(30)}
    \pgfmathsetmacro{\s}{\a-\r*cos(30)}

    % Draw the strip outside the two disks. Its ends are circular arcs.
    \newcommand{\drawtube}[2]{%
        \fill[#1] (-\s,\w) -- (\s,\w) arc (150:210:\r)
                    -- (-\s,-\w) arc (-30:30:\r) -- cycle;
        \draw (-\s,\w) -- (\s,\w)  (-\s,-\w) -- (\s,-\w);
        \draw[interface,#2] (-\s,-\w) arc (-30:30:\r);
        \draw[interface,#2] ( \s, \w) arc (150:210:\r);
    }

    % Left panel: a fragment of a triangular network of disk rooms.
    \draw (0,0) circle[radius=1.65];
    \node at (0,1.32) {$\mathbb B$};
    \begin{scope}[yshift=.68cm,scale=.76]
        \pgfmathsetmacro{\depth}{sqrt(3)*\a}
        \coordinate (A) at (-\a,0);
        \coordinate (B) at ( \a,0);
        \coordinate (C) at (-2*\a,-\depth);
        \coordinate (D) at (0,-\depth);
        \coordinate (E) at ( 2*\a,-\depth);

        % Each port is specified by its direction and angular half-width.
        \newcommand{\netroom}[2]{%
            \begin{scope}[shift={(#1)}]
                \draw (0,0) circle[radius=\r];
                \foreach \direction/\aperture in {#2}{
                    \pgfmathsetmacro{\lo}{\direction-\aperture}
                    \pgfmathsetmacro{\hi}{\direction+\aperture}
                    \draw[white,line width=1.1pt]
                        (\lo:\r) arc (\lo:\hi:\r);
                    \draw[interface,black!55]
                        (\lo:\r) arc (\lo:\hi:\r);
                }
            \end{scope}
        }
        \netroom{A}{0/30,240/18,300/18}
        \netroom{B}{180/30,240/18,300/18}
        \netroom{C}{0/18,60/18,300/18}
        \netroom{D}{0/18,60/18,120/18,180/18}
        \netroom{E}{120/18,180/18,240/18}

        \drawtube{black!12}{black}
        \node[font=\scriptsize] at (A) {$R_i$};
        \node[font=\scriptsize] at (B) {$R_j$};

        \pgfmathsetmacro{\offset}{\r*cos(18)}
        \pgfmathsetmacro{\halfwidth}{\r*sin(18)}
        \pgfmathsetmacro{\farend}{2*\a-\offset}
        \foreach \site/\direction in {A/240,A/300,B/240,B/300,C/0,D/0}{
            \begin{scope}[shift={(\site)},rotate=\direction]
                \foreach \side in {-\halfwidth,\halfwidth}
                    \draw (\offset,\side) -- (\farend,\side);
            \end{scope}
        }

        % Open tube walls and three dots indicate the omitted network.
        \foreach \site/\direction in {C/300,E/240}{
            \begin{scope}[shift={(\site)},rotate=\direction]
                \foreach \side in {-\halfwidth,\halfwidth}
                    \draw (\offset,\side) -- (.79,\side);
                \foreach \distance in {.92,1.06,1.20}
                    \fill (\distance,0) circle[radius=.024];
            \end{scope}
        }
    \end{scope}

    % Right panel: the selected edge in detail.
    \begin{scope}[xshift=5.25cm, scale=2]
        \drawtube{black!12}{black}
        \draw (-\a,0) ++(30:\r)  arc (30:330:\r);
        \draw ( \a,0) ++(210:\r) arc (210:510:\r);
        % Midpoints of the two circular interface arcs.
        \coordinate (Gi) at ({-\a+\r},0);
        \coordinate (Gj) at ({ \a-\r},0);
        \node at (-\a,.11) {$R_i$};
        \node at (0,0)     {$T_e$};
        \node at ( \a,.11) {$R_j$};

        \node[font=\footnotesize,inner sep=1pt] (Li)
            at (-.73,-.13) {$\Gamma_{ei}$};
        \node[font=\footnotesize,inner sep=1pt] (Lj)
            at ( .73,-.13) {$\Gamma_{ej}$};
        \draw[->,>=stealth,thin] (Li.east) -- (Gi);
        \draw[->,>=stealth,thin] (Lj.west) -- (Gj);
    \end{scope}
\end{tikzpicture}
\end{figure}

Let $h_e$ minimize $\int_{T_e}|\nabla h|^2$ subject to $h|_{\Gamma_{ei}}=0$ and $h|_{\Gamma_{ej}}=1$. By our previous work, such a minimizer $h_e$ exists and is unique; see \cite[p.~12]{CLZ}. Moreover, $0\leq h_e\leq 1$. Put
\[
 \kappa_e=\int_{T_e}|\nabla h_e|^2>0.
\]
For $0<\tau\le1$, define
\begin{equation}\label{eq:planar-density}
 \widehat\rho_{\tau,t}=
 \begin{cases}
  1,&\text{in }R_i,\\
  \frac{\tau vc_e(t)}{\kappa_e},&\text{in }T_e,\\
  \tau^2,&\text{in }\mathbb B\setminus\overline D,
 \end{cases}
 \qquad b_t=\sum_{i=1}^N b_i(t)\mathbf1_{R_i}.
\end{equation}
The corresponding forms on $H^1_0(\mathbb B)$ are
\[
 \begin{aligned}
 \widehat a^0_{\tau,t}(f,g)
   &=\int_{\mathbb B}\widehat\rho_{\tau,t}\nabla f\cdot\nabla g,\\
 \widehat b_{\tau,t}(f,g)
   &=\int_{\mathbb B}\widehat\rho_{\tau,t}fg,\\
 \widehat a_{\tau,t}(f,g)
   &=\widehat a^0_{\tau,t}(f,g)
     +\tau\int_{\mathbb B}b_t\widehat\rho_{\tau,t}fg.
 \end{aligned}
\]
A nonzero $f\in H^1_0(\mathbb B)$ is a critical point of the Rayleigh
quotient $\widehat a_{\tau,t}(f,f)/\widehat b_{\tau,t}(f,f)$ with
critical value $\nu$ if and only if
\[
    \widehat a_{\tau,t}(f,g)=\nu\,\widehat b_{\tau,t}(f,g)
    \qquad\text{for all }g\in H^1_0(\mathbb B).
\]
Equivalently,
\[
    -\operatorname{div}(\widehat\rho_{\tau,t}\nabla f)
    +\tau b_t\widehat\rho_{\tau,t}f
    =\nu\widehat\rho_{\tau,t}f
\]
in the sense of distributions. Under the unitary change of variables
$u=\sqrt{\widehat\rho_{\tau,t}}\,f$, this equation becomes
$\mathcal L_{\tau,t}u=\nu u$ in $L^2(\mathbb B)$, where
\begin{align}\label{Sec 7 equ-c}
    \mathcal L_{\tau,t}u
    =-\frac{1}{\sqrt{\widehat\rho_{\tau,t}}}
      \operatorname{div}\!\left(
          \widehat\rho_{\tau,t}
          \nabla\!\left(\frac{u}{\sqrt{\widehat\rho_{\tau,t}}}\right)
      \right)
      +\tau b_tu.
\end{align}
Its domain consists of all $u=\sqrt{\widehat\rho_{\tau,t}}\,f$ with
$f\in H^1_0(\mathbb B)$ and
$\operatorname{div}(\widehat\rho_{\tau,t}\nabla f)\in L^2(\mathbb B)$.
Define the fixed isometry $J_0:\mathbb R^N\to L^2(\mathbb B)$ and
the operator $T(t)$ by
\begin{align}\label{Sec 7 Lem 1 equ-s}
    J_0x=\frac{1}{\sqrt v}\sum_{i=1}^N x_i\mathbf 1_{R_i},
    \qquad
    T(t)=J_0(\mathcal{M}(t)+\kappa I)^{-1}J_0^*.
\end{align}
\begin{lemma}\label{lem:rooms}
For every $\delta_0>0$, there exist $\tau>0$ and real potentials
$V_t\in C^\infty(\overline{\mathbb B})$, continuous in
$t\in\mathcal B$ in the uniform norm, such that, for $A_t=\tau^{-1}(-\Delta_{\mathbb B}+V_t)$,
\[
 A_t+\kappa\ge I,\qquad
 \sup_{t\in\mathcal B}
 \|(A_t+\kappa)^{-1}-T(t)\|<\delta_0.
\]
where $\kappa$ is defined by \eqref{Sec 7 equ-a}.
\end{lemma}

\begin{proof}
Define $J:\R^N\to H^1(D)$ by
\[
 (Jx)|_{R_i}=\frac{x_i}{\sqrt v},\qquad
 (Jx)|_{T_e}=\frac{(1-h_e)x_i+h_ex_j}{\sqrt v},
 \qquad e=\{i,j\}.
\]
Using a fixed extension operator and a cutoff, extend $J$ to
$\widehat J:\R^N\to H^1_0(\mathbb B)$. The calculation in
\cite[Lemma 4.3 (4.8)]{CLZ} gives, uniformly in $t$,
\begin{equation}\label{eq:planar-trial-forms}
 \begin{aligned}
 \widehat b_{\tau,t}(\widehat Jx,\widehat Jy)
   &=x\cdot y+O(\tau)|x||y|,\\
 \widehat a_{\tau,t}(\widehat Jx,\widehat Jy)
   &=\tau x^T\mathcal{M}(t)y+O(\tau^2)|x||y|.
 \end{aligned}
\end{equation}
Using the Poincar\'e inequality and \eqref{eq:planar-density}, we can get
\[
 \tau^2\|f\|_{L^2(\mathbb B)}^2
 \le C\tau^2\|\nabla f\|_{L^2(\mathbb B)}^2
 \le C\widehat a^0_{\tau,t}(f,f),
\]
they give $\widehat b_{\tau,t}(f,f)\le
C\widehat a^0_{\tau,t}(f,f)$ whenever
$\int_{R_i}f=0$ for every $i=1,\cdots, N$. The added multiplication operator
has norm at most $\tau\sup_t\|b_t\|_\infty$ in the weighted
space. Applying the min--max principle for functional $\widehat a_{\tau,t}(f,f)/\widehat b_{\tau,t}(f,f) $ and
\eqref{eq:planar-trial-forms} imply
\begin{equation}\label{eq:planar-gap}
 \begin{gathered}
 -C\tau\le\nu_1(\tau,t)\le\cdots\le\nu_N(\tau,t)\le C\tau,\quad
 \nu_{N+1}(\tau,t)\ge\widehat\gamma>0.
 \end{gathered}
\end{equation}
for all sufficiently small $\tau$, where $\nu_j(\tau,t)$ are
the eigenvalues of $\mathcal L_{\tau,t}$.

Let $\widehat P$ be the $\widehat b$-orthogonal projection onto the
spectral subspace corresponding to $\nu_1,\ldots,\nu_N$. We suppress
the subscripts $\tau,t$ on the forms and the density.

For $x\in\mathbb R^N$, set $z=(I-\widehat P)\widehat Jx$. The
gradient estimate in \cite[(4.10)]{CLZ}, together with
\[
    \left|\tau\int_{\mathbb B}b_t\widehat\rho\,(\widehat Jx)g\right|
    \le C\tau|x|\,\widehat b(g,g)^{1/2},
\]
gives
\[
    |\widehat a(\widehat Jx,g)|
    \le C\tau|x|
       \bigl(\widehat a^0(g,g)+\widehat b(g,g)\bigr)^{1/2}.
\]
Since $z$ belongs to the high spectral subspace,
\eqref{eq:planar-gap} gives
$\widehat b(z,z)\le\widehat\gamma^{-1}\widehat a(z,z)$.
As $\widehat a^0$ and $\widehat a$ differ by at most
$\tau\|b_t\|_\infty\widehat b$ on the diagonal, taking $g=z$ and
using spectral orthogonality yields
\[
    \widehat a(z,z)=\widehat a(\widehat Jx,z)
    \le C\tau|x|\,\widehat a(z,z)^{1/2}.
\]
Consequently,
\[
    \widehat a(z,z)+\widehat b(z,z)\le C\tau^2|x|^2.
\]

By spectral orthogonality and \eqref{eq:planar-trial-forms},

\begin{align}
    \widehat b(\widehat P\widehat Jx,\widehat P\widehat Jx)
    &=\widehat b(\widehat Jx,\widehat Jx)-\widehat b(z,z)
      =|x|^2+O(\tau)|x|^2,\label{Sec 7 Lem equ-a}\\
    \widehat a(\widehat P\widehat Jx,\widehat P\widehat Jx)
    &=\widehat a(\widehat Jx,\widehat Jx)-\widehat a(z,z)
      =\tau x^T\mathcal{M}(t)x+O(\tau^2)|x|^2.\label{Sec 7 Lem equ-b}
\end{align}

Define
\[
 \begin{aligned}
 (\widehat G_{\tau,t})_{ij}=\widehat b_{\tau,t}
  (\widehat P_{\tau,t}\widehat Je_i,
   \widehat P_{\tau,t}\widehat Je_j),\quad
 U_{\tau,t}=\sqrt{\widehat\rho_{\tau,t}}\,
       \widehat P_{\tau,t}\widehat J\widehat G_{\tau,t}^{-1/2},\quad 
 \widehat B_{\tau,t}=U_{\tau,t}^*\mathcal L_{\tau,t}U_{\tau,t},
 \end{aligned}
\]
where $U_{\tau,t}^*$ is the adjoint operator of $U_{\tau,t} $, then $U_{\tau,t}^*U_{\tau,t}=I $. The preceding energy estimate \eqref{Sec 7 Lem equ-b} gives
\begin{align}\label{Sec 7 Lem equ-d}
    \widehat B_{\tau,t}
    =\widehat G_{\tau,t}^{-1/2}
       \bigl[\tau \mathcal{M}(t)+O(\tau^2)\bigr]
       \widehat G_{\tau,t}^{-1/2}
    =\tau \mathcal{M}(t)+O(\tau^2),
\end{align}

Finally, we also claim that
\begin{equation}\label{eq:planar-low-matrix}
    \|U_{\tau,t}-J_0\|\le C\sqrt\tau .
\end{equation}
On each room, $\sqrt{\widehat\rho}\,\widehat Jx=J_0x$.
The density is $O(\tau)$ on the tubes and $\tau^2$ on the exterior,
so $\|\sqrt{\widehat\rho}\,\widehat J-J_0\|\le C\sqrt\tau$. The bound on $\widehat b(z,z)$ also gives
$\|\sqrt{\widehat\rho}(I-\widehat P)\widehat J\|\le C\tau$,
and \eqref{Sec 7 Lem equ-a} gives
$\widehat G_{\tau,t}^{-1/2}=I+O(\tau)$.
The definition of $U_{\tau,t}$ and the triangle inequality now give
\eqref{eq:planar-low-matrix}.

Since $\tau^{-1}\mathcal L_{\tau,t}+\kappa\ge I$,
the spectral decomposition gives
\[
 \begin{aligned}
 (\tau^{-1}\mathcal L_{\tau,t}+\kappa)^{-1}
 &=U_{\tau,t}(\tau^{-1}\widehat B_{\tau,t}+\kappa I)^{-1}
       U_{\tau,t}^*+O(\tau).
 \end{aligned}
\]
It follows from \eqref{Sec 7 Lem equ-d} and \eqref{eq:planar-low-matrix} that
\begin{equation}\label{net:rooms}
 \sup_{t\in\mathcal B}
 \|(\tau^{-1}\mathcal L_{\tau,t}+\kappa)^{-1}-T(t)\|
 \le C\sqrt\tau.
\end{equation}

Choose $\tau>0$ so that the right-hand side of \eqref{net:rooms}
is less than $\delta_0/2$, and fix this $\tau$. Extend
$\widehat\rho_{\tau,t}$ by $\tau^2$ and $b_t$ by zero outside
$\mathbb B$. For a nonnegative mollifier $\psi_\eta$ of integral
one, put
\[
    \widehat\rho_{\tau,t,\eta}:=\psi_\eta*\widehat\rho_{\tau,t}, \qquad\varphi_t=(\widehat\rho_{\tau,t,\eta})^{1/2},
    \qquad b_{t,\eta}=\psi_\eta*b_t .
\]
We suppress the dependence of $\varphi_t$ on $\eta$. For this fixed
$\tau$, we have $\|b_{t,\eta}\|_{L^\infty(\mathbb B)}\leq
\|b_t\|_{L^\infty(\mathbb B)}$ and
$0<c_{\tau}\leq\widehat\rho_{\tau,t,\eta}\leq C_{\tau}$ for some
positive constants $c_\tau,C_\tau$ independent of $t$ and $\eta$.
Moreover, $\kappa+b_{t,\eta}\ge1$ by \eqref{Sec 7 equ-a}.

We follow the quadratic form calculation in
\cite[(4.24)--(4.25)]{CLZ}. For $f\in H^1_0(\mathbb B)$, put
$u=\varphi_t f$. Direct computation gives
\[
   \int_{\mathbb B}\widehat\rho_{\tau,t,\eta}|\nabla f|^2
    =\int_{\mathbb B}\varphi_t^2|\nabla f|^2
    =\int_{\mathbb B}|\nabla u-u\nabla\log\varphi_t|^2
    =\int_{\mathbb B}\left(
        |\nabla u|^2+\frac{\Delta\varphi_t}{\varphi_t}u^2
      \right).
\]
Using $\int_{\mathbb B}\widehat\rho_{\tau,t,\eta}f^2
=\int_{\mathbb B}u^2$, for $f\ne0$,  we can obtain 
\begin{align*}
   \frac{\int_{\mathbb B}\widehat\rho_{\tau,t,\eta}|\nabla f|^2
   +\tau\int_{\mathbb B}b_{t,\eta}\widehat\rho_{\tau,t,\eta}f^2}
   {\int_{\mathbb B}\widehat\rho_{\tau,t,\eta}f^2}
   =\frac{\int_{\mathbb B}\left(
        |\nabla u|^2+\left(\frac{\Delta\varphi_t}{\varphi_t}
        +\tau b_{t,\eta}\right)u^2
      \right)}{\int_{\mathbb B}u^2}.
\end{align*}
Thus,  the two eigenvalue problems are equivalent under the map
$f\mapsto\varphi_t f$.  Set
$V_t=\Delta\varphi_t/\varphi_t+\tau b_{t,\eta}$, suppressing its
dependence on $\eta$. Then $f$ is a critical point of the weighted
Rayleigh quotient with critical value $\nu$ if and only if
\begin{align*}
   \mathcal L_{\tau,t,\eta}u=(-\Delta_{\mathbb B}+V_t)u=\nu u,
\end{align*}
where $\mathcal L_{\tau,t,\eta}$ is defined as in
\eqref{Sec 7 equ-c}, replacing $\widehat\rho_{\tau,t}$ by
$\widehat\rho_{\tau,t,\eta}$ and $b_t$ by $b_{t,\eta}$.

Throughout this argument, $\tau>0$ is fixed. In the setting of
\cite[Lemma~4.5]{CLZ}, the authors prove convergence of each
fixed-index eigenvalue and of the energy matrix on the selected
spectral subspace. However, as $\eta\downarrow0$,   we need to show that
$(\tau^{-1}\mathcal L_{\tau,t,\eta}+\kappa)^{-1}$ converges to
$(\tau^{-1}\mathcal L_{\tau,t}+\kappa)^{-1}$    uniformly in $t\in\mathcal B$. 

\smallskip
\textbf{Step 1. Uniform estimate.}
Let $\eta_n\downarrow0$, $t_n\to t$ in $\mathcal B$, and
$F_n\rightharpoonup F$ in $L^2(\mathbb B)$. Write
$(\widehat\rho_n,b_n)$ for either
$(\widehat\rho_{\tau,t_n},b_{t_n})$ or
$(\widehat\rho_{\tau,t_n,\eta_n},b_{t_n,\eta_n})$.
In either case, the weighted solution $f_n\in H^1_0(\mathbb B)$
satisfies
\[
    \tau^{-1}\int_{\mathbb B}\widehat\rho_n\nabla f_n\cdot\nabla g
    +\int_{\mathbb B}(\kappa+b_n)\widehat\rho_n f_ng
    =\int_{\mathbb B}F_n\sqrt{\widehat\rho_n}\,g,
    \qquad g\in H^1_0(\mathbb B).
\]
 Taking $g=f_n$ and using
$c_\tau\le\widehat\rho_n\le C_\tau$ and $\kappa+b_n\ge1$, we obtain
\[
    c_\tau\left(\tau^{-1}\|\nabla f_n\|_{L^2}^2
                      +\|f_n\|_{L^2}^2\right)
    \le C_\tau^{1/2}\|F_n\|_{L^2}\|f_n\|_{L^2}.
\]
Since $(F_n)$ is bounded in $L^2$, the sequence $(f_n)$ is bounded
in $H^1_0(\mathbb B)$.

\textbf{Step 2. Uniqueness of  limit.}
By Rellich compactness, after passing to a subsequence,
\[
    f_n\rightharpoonup f \quad\text{in }H^1_0(\mathbb B),
    \qquad f_n\to f \quad\text{in }L^2(\mathbb B).
\]
For either choice of coefficients, 
\[
    \widehat\rho_n\to\widehat\rho_{\tau,t},\qquad b_n\to b_t
    \quad\text{in measure on }\mathbb B.
\]
The coefficients are uniformly bounded. Hence, for each fixed
$g\in H^1_0(\mathbb B)$,
\[
\begin{aligned}
    \widehat\rho_n\nabla g
        \to\widehat\rho_{\tau,t}\nabla g,\quad
    (\kappa+b_n)\widehat\rho_n g
        \to(\kappa+b_t)\widehat\rho_{\tau,t}g,\quad
    \sqrt{\widehat\rho_n}\,g
        \to\sqrt{\widehat\rho_{\tau,t}}\,g
\end{aligned}
\]
in $L^2(\mathbb{B})$. 
Letting $n\to\infty$ yields
\[
    \tau^{-1}\int_{\mathbb B}\widehat\rho_{\tau,t}\nabla f\cdot\nabla g
    +\int_{\mathbb B}(\kappa+b_t)\widehat\rho_{\tau,t}fg
    =\int_{\mathbb B}F\sqrt{\widehat\rho_{\tau,t}}\,g,
    \qquad g\in H^1_0(\mathbb B).
\]
This equation has a unique solution by coercivity. Thus every
convergent subsequence has the same limit, and $f_n\to f$ in
$L^2$ along the whole sequence. Moreover,
\[
    \|\sqrt{\widehat\rho_n}\,f_n
             -\sqrt{\widehat\rho_{\tau,t}}\,f\|_{L^2}
    \le C_\tau^{1/2}\|f_n-f\|_{L^2}
       +\|(\sqrt{\widehat\rho_n}
             -\sqrt{\widehat\rho_{\tau,t}})f\|_{L^2}
    \longrightarrow0 .
\]
Consequently, the original and mollified resolvent images of $F_n$
both converge in $L^2$ to
$(\tau^{-1}\mathcal L_{\tau,t}+\kappa)^{-1}F$.

\smallskip
\textbf{Step 3. Uniform convergence in operator norm.}
If the asserted convergence failed, there would exist
$\varepsilon>0$, $\eta_n\downarrow0$, $t_n\in\mathcal B$, and
$F_n\in L^2(\mathbb B)$ with $\|F_n\|_{L^2}=1$ such that
\[
    \left\|
       (\tau^{-1}\mathcal L_{\tau,t_n,\eta_n}+\kappa)^{-1}F_n
       -(\tau^{-1}\mathcal L_{\tau,t_n}+\kappa)^{-1}F_n
    \right\|_{L^2}\ge\varepsilon .
\]
By compactness of $\mathcal B$ and weak sequential compactness
of the $L^2$ unit ball, we may pass to a subsequence for which
$t_n\to t$ and $F_n\rightharpoonup F$. Step~2 shows that the
left-hand side tends to zero, a contradiction. Recalling that
$\mathcal L_{\tau,t,\eta}=-\Delta_{\mathbb B}+V_t$, we conclude
that
\begin{align}\label{eq:planar-mollified-resolvent}
    \lim_{\eta\downarrow0}\sup_{t\in\mathcal B}
    \left\|
        \bigl(\tau^{-1}(-\Delta_{\mathbb B}+V_t)+\kappa\bigr)^{-1}
        -\bigl(\tau^{-1}\mathcal L_{\tau,t}+\kappa\bigr)^{-1}
    \right\|=0 .
\end{align}

Choose $\eta>0$ so that the  error in \eqref{eq:planar-mollified-resolvent} is less than
$\delta_0/2$, uniformly in $t$. Combining  with \eqref{net:rooms} yields
\[
    \sup_{t\in\mathcal B}
    \bigl\|(A_t+\kappa)^{-1}-T(t)\bigr\|<\delta_0 .
\]The proof is complete.
\end{proof}

The approximation of Schr\"odinger operators by Dirichlet Laplacians
on planar domains is a classical topic in homogenization theory;
see \cite{DM,KP}. The basic idea is to remove a large number of small
interior holes whose capacities are chosen to approximate the
prescribed potential.

\begin{lemma}\label{lem:holes}
Let 
$t\mapsto W_t\in C(\overline{\mathbb B})$ be continuous in the uniform norm, with 
$W_t\ge1$. There exist domains $\Omega_{\delta,t}$ obtained
from $\mathbb B$ by removing finitely many disjoint closed
interior disks such that 
\begin{align}\label{Sec 7 equ-f}
 \sup_{t\in\mathcal B}
 \|(-\Delta_{\Omega_{\delta,t}}+1)^{-1}
       -(-\Delta_{\mathbb B}+W_t+1)^{-1}\|
 \longrightarrow0\qquad(\delta\downarrow0).
\end{align}
Moreover, for every fixed $\delta>0$, the map
\[
    t\longmapsto(-\Delta_{\Omega_{\delta,t}}+1)^{-1}
\]
is continuous in the operator norm on $L^2(\mathbb B)$, where
the Dirichlet resolvents are extended by zero outside
$\Omega_{\delta,t}$.
\end{lemma}

\begin{proof}
 We give the two-dimensional
construction and the uniformity needed here.
Let $Q_z$ be the squares of the side-$\delta$ grid whose
closures lie in $\mathbb B$, indexed by their centers $z$.
Set
\begin{equation}\label{eq:circular-hole-radius}
 R_\delta=\frac\delta4,\qquad
 r_z(t)=R_\delta
       \exp\!\left(-\frac{2\pi}{\delta^2W_t(z)}\right),
 \qquad
 \Omega_{\delta,t}=\mathbb B\setminus
       \bigcup_z\overline{B(z,r_z(t))}.
\end{equation}

% Preamble: \usepackage{tikz,amssymb,float}
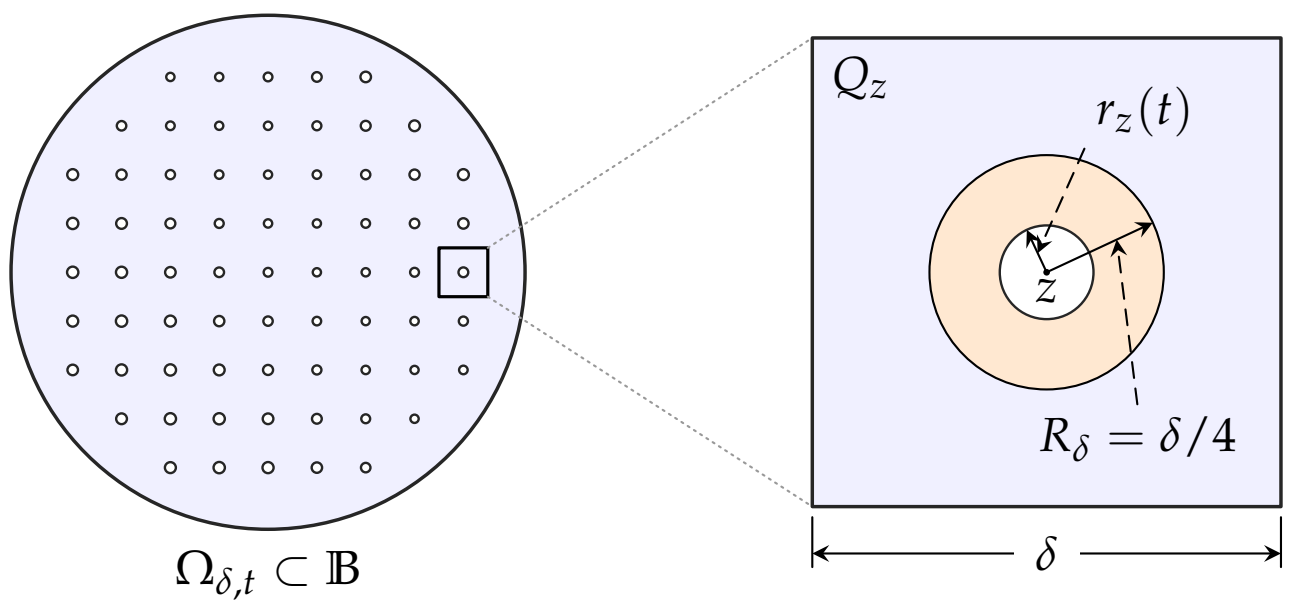
\begin{figure}[H]
\centering
\begin{tikzpicture}[
    font=\small, >=stealth, line cap=round, line join=round,
    disk region/.style={draw=black!85, fill=blue!6, line width=.65pt},
    hole/.style={draw=black!85, fill=white, line width=.45pt},
    guide/.style={draw=black!40, densely dotted, thin}
]
% Left: the unit disk. Keep only cells whose four corners lie in it.
\begin{scope}[scale=1.7]
    \def\mesh{.19}
    \draw[disk region] (0,0) circle[radius=1];
    \foreach \i in {-5,...,5}{
        \foreach \j in {-5,...,5}{
            \pgfmathtruncatemacro{\inside}{
                \mesh^2*((abs(\i)+.5)^2+(abs(\j)+.5)^2)<1}
            \ifnum\inside=1
                \pgfmathsetmacro{\rh}{.019*(1+.16*sin(35*(\i-4)+25*\j))}
                \draw[hole] ({\i*\mesh},{\j*\mesh}) circle[radius=\rh];
            \fi
        }
    }
    % The selected cell has center (4 delta,0).
    \draw[semithick] ({3.5*\mesh},{-.5*\mesh})
        rectangle ({4.5*\mesh},{.5*\mesh});
    \coordinate (upper) at ({4.5*\mesh},{.5*\mesh});
    \coordinate (lower) at ({4.5*\mesh},{-.5*\mesh});
\end{scope}
\node at (0,-2.00) {$\Omega_{\delta,t}\subset\mathbb B$};

% Right: the selected cell. Its side is 3.10 drawing units.
\begin{scope}[xshift=5.15cm]
    \def\halfside{1.55}
    \def\outerradius{.775}  % R_delta = delta/4
    \def\holeradius{.31}   % schematic hole radius
    \draw[guide] (upper) -- (-\halfside,\halfside);
    \draw[guide] (lower) -- (-\halfside,-\halfside);
    \draw[disk region] (-\halfside,-\halfside) rectangle (\halfside,\halfside);
    \node[below right] at (-\halfside,\halfside) {$Q_z$};
    \fill[orange!18] (0,0) circle[radius=\outerradius];
    \draw[thin] (0,0) circle[radius=\outerradius];
    \draw[hole] (0,0) circle[radius=\holeradius];

    \draw[->,thin] (0,0) -- (115:\holeradius)
        coordinate[pos=.45] (innerMark);
    \draw[->,thin] (0,0) -- (25:\outerradius)
        coordinate[pos=.66] (outerMark);

    \node[inner sep=2pt] (innerLabel) at (.65,1.06) {$r_z(t)$};
    \node[inner sep=2pt] (outerLabel) at (.60,-1.10) {$R_\delta=\delta/4$};
    \draw[->,densely dashed,thin] (innerLabel.south west) -- (innerMark);
    \draw[->,densely dashed,thin] (outerLabel.north) -- (outerMark);
    \fill (0,0) circle[radius=.022];
    \node[below,inner sep=1pt] at (0,0) {$z$};

    \draw[thin] (-\halfside,-1.62) -- (-\halfside,-1.95)
                (\halfside,-1.62) -- (\halfside,-1.95);
    \draw[<->,thin]
        (-\halfside,-1.86) -- (\halfside,-1.86)
        node[midway,fill=white,inner xsep=4pt] {$\delta$};
\end{scope}
\end{tikzpicture}
\end{figure}

 Define $\chi_{\delta,t}$ on $\mathbb B$ by setting it equal to zero
on each closed disk $\overline{B(z,r_z(t))}$ and equal to one on
$\mathbb B\setminus\bigcup_z B(z,R_\delta)$. On each intervening
annulus, let $\chi_{\delta,t}$ be the unique solution of
\[
    \begin{cases}
        \Delta\chi_{\delta,t}=0
        & \text{in } B(z,R_\delta)\setminus\overline{B(z,r_z(t))},\\
        \chi_{\delta,t}=0
        & \text{on } \partial B(z,r_z(t)),\\
        \chi_{\delta,t}=1
        & \text{on } \partial B(z,R_\delta).
    \end{cases}
\]
This solution is radial about $z$ and is given by
\[
    \chi_{\delta,t}(x)
    =\frac{\log\bigl(|x-z|/r_z(t)\bigr)}
           {\log\bigl(R_\delta/r_z(t)\bigr)},
    \qquad r_z(t)<|x-z|<R_\delta.
\]
Hence, integrating in polar coordinates and using
\eqref{eq:circular-hole-radius}, we obtain
\begin{equation}\label{eq:circular-hole-capacity}
    \int_{B(z,R_\delta)\setminus\overline{B(z,r_z(t))}}
        |\nabla\chi_{\delta,t}|^2\,\mathrm dx
    =\frac{2\pi}{\log\bigl(R_\delta/r_z(t)\bigr)}
    =\delta^2W_t(z).
\end{equation}
Writing $\log(R_\delta/r_z(t))=2\pi/(\delta^2W_t(z))$, direct integration gives
\[
 \int_{Q_z}|1-\chi_{\delta,t}|^2
 \le C r_z(t)^2+C\frac{R_\delta^2}{[\log(R_\delta/r_z(t))]^2}.
\]
The uniform bounds on $W_t$ therefore imply
$\chi_{\delta,t}\to1$ in $L^2(\mathbb B)$ as $\delta\to 0$ and
$\|\nabla\chi_{\delta,t}\|_{L^2(\mathbb B) }\le C$.
In particular, $\nabla\chi_{\delta,t}\rightharpoonup0$.

Now, we begin to study the convergence as $\delta\downarrow0$ and $t\to t_0$. Let $F_\delta\rightharpoonup F$ in $L^2(\mathbb B)$, and let $u_\delta$ solve the boundary value problem
\begin{align}\label{Sec 7 equ-d}
    -\Delta_{\Omega_{\delta,t}}u_{\delta}+u_{\delta}=F_{\delta} \quad\mathrm{in}\quad\Omega_{\delta,t},\qquad u_{\delta}=0\quad\mathrm{on}\quad\partial\Omega_{\delta,t}
\end{align}

By zero extension, we may view
$u_\delta\in H^1_0(\mathbb{B})$, and they
are bounded in $H^1_0(\mathbb B)$ by a standard energy estimate. Passing to a subsequence,
let $u_\delta\rightharpoonup u$ in $H^1$ and
$u_\delta\to u$ in $L^2$. Following arguments similar to those in \cite[Section 5]{CLZ}, we claim that the limit equation is 
\begin{align*}
    -\Delta_{\mathbb B}u+(W_{t_0}+1)u= F\quad\mathrm{in}\quad\mathbb{B},\qquad u=0\quad\mathrm{on}\quad\partial\mathbb{B}.
\end{align*}

Fix $\phi\in C_c^\infty(\mathbb B)$.
 On $\partial B(z,R_\delta)$, the radial formula for $\chi_{\delta,t}$
and \eqref{eq:circular-hole-capacity} give
\[
    \partial_\nu\chi_{\delta,t}
    =\frac{1}{R_\delta\log(R_\delta/r_z(t))}
    =\frac{\delta^2W_t(z)}{2\pi R_\delta}.
\]
Since $\chi_{\delta,t}$ is harmonic on each annulus and $u_\delta$
has zero trace on its inner circle, integration by parts gives
\[
    \int_{B(z,R_\delta)\setminus\overline{B(z,r_z(t))}}
        \nabla(u_\delta\phi)\cdot\nabla\chi_{\delta,t}\,\mathrm dx
    =\int_{\partial B(z,R_\delta)}
        u_\delta\phi\,\partial_\nu\chi_{\delta,t}\,\mathrm ds
    =\delta^2W_t(z)
        \fint_{\partial B(z,R_\delta)}u_\delta\phi\,\mathrm ds.
\]
As $\nabla\chi_{\delta,t}$ vanishes almost everywhere outside the
annuli, summing over $z$ and expanding $\nabla(u_\delta\phi)$ yields
\begin{equation}\label{eq:hole-flux}
    \int_{\mathbb B}\phi\nabla u_\delta\cdot\nabla\chi_{\delta,t}\,\mathrm dx
    =\sum_z\delta^2W_t(z)
        \fint_{\partial B(z,R_\delta)}u_\delta\phi\,\mathrm ds
    -\int_{\mathbb B}u_\delta\nabla\phi\cdot\nabla\chi_{\delta,t}\,\mathrm dx.
\end{equation}
The last term tends to zero by the strong $L^2$ convergence and
the weak convergence of $\nabla\chi_{\delta,t}$. Applying the trace and Poincar\'e inequalities as in \cite[Lemma 5.2]{CLZ} gives
\[
 \left|\fint_{\partial B(z,R_\delta)}g
               -\fint_{Q_z}g\right|
 \le C\|\nabla g\|_{L^2(Q_z)}.
\]
Consequently,
\begin{equation}\label{eq:hole-average-error}
 \left|\sum_z\delta^2W_t(z)
    \left(\fint_{\partial B(z,R_\delta)}g
                   -\fint_{Q_z}g\right)\right|
 \le C\delta\|\nabla g\|_{L^2(\mathbb B)}.
\end{equation}
Apply \eqref{eq:hole-average-error} to $g=u_\delta\phi$.
The uniform $H^1$ bound gives an $O(\delta)$ error in replacing
the circle averages in the first term of \eqref{eq:hole-flux} by cell
averages. For the second term, write $\chi=\chi_{\delta,t}$.
Since $R_\delta=\delta/4$ and
$\log(R_\delta/r_z(t))\ge c\delta^{-2}$, summing the cell estimate
following \eqref{eq:circular-hole-capacity} over $O(\delta^{-2})$
cells gives $\|1-\chi\|_2=O(\delta^2)$.
Since $\phi$ is smooth and compactly supported, integration by parts
on $\mathbb B$ gives
\[
 \int_{\mathbb B}u_\delta\nabla\phi\cdot\nabla\chi\,\mathrm dx
 =-\int_{\mathbb B}(\chi-1)
       (\nabla u_\delta\cdot\nabla\phi+u_\delta\Delta\phi)\,\mathrm dx
 =O(\delta^2).
\] Hence, by \eqref{eq:hole-flux} and $|Q_z|=\delta^2$,
\[
    \int_{\mathbb B}\phi\,\nabla u_\delta\cdot\nabla\chi_{\delta,t}\,\mathrm dx
    =\sum_z\int_{Q_z}W_t(z)u_\delta\phi\,\mathrm dx+O(\delta)
    \longrightarrow
    \int_{\mathbb B}W_{t_0}u\phi\,\mathrm dx.
\]
Indeed, for small $\delta$ the squares cover $\operatorname{supp}\phi$,
the step functions with value $W_t(z)$ on $Q_z$ converge uniformly
to $W_{t_0}$ there, and $u_\delta\to u$ in $L^2(\mathbb B)$.

Testing the weak equation for \eqref{Sec 7 equ-d} with
$\chi_{\delta,t}\phi\in H^1_0(\Omega_{\delta,t})$ gives
\[
    \int_{\mathbb B}\chi_{\delta,t}\nabla u_\delta\cdot\nabla\phi\,\mathrm dx
    +\int_{\mathbb B}\phi\,\nabla u_\delta\cdot\nabla\chi_{\delta,t}\,\mathrm dx
    +\int_{\mathbb B}u_\delta\chi_{\delta,t}\phi\,\mathrm dx
    =\int_{\mathbb B}F_\delta\chi_{\delta,t}\phi\,\mathrm dx.
\]
Since $\chi_{\delta,t}\to1$ in $L^2(\mathbb B)$ and $\phi$ is smooth
with compact support, both $\chi_{\delta,t}\phi$ and
$\chi_{\delta,t}\nabla\phi$ converge strongly in $L^2$ to $\phi$ and
$\nabla\phi$, respectively. Using the preceding limit for the second
term and the weak convergences $\nabla u_\delta\rightharpoonup\nabla u$
and $F_\delta\rightharpoonup F$, we obtain
\[
    \int_{\mathbb B}
    \bigl(\nabla u\cdot\nabla\phi+(1+W_{t_0})u\phi\bigr)\,\mathrm dx
    =\int_{\mathbb B}F\phi\,\mathrm dx.
\]

 To obtain \eqref{Sec 7 equ-f}, we apply a contradiction argument similar to that in \cite[Proposition~5.5]{CLZ}. Since the details are very similar, we omit them here.

It remains to prove norm-continuity for fixed $\delta$.
Write
\[
    B_{\delta,t}
    =(-\Delta_{\Omega_{\delta,t}}+1)^{-1},
\]
with zero extension to $L^2(\mathbb B)$.
Let $t_n\to t_0$. Since $W_{t_n}\to W_{t_0}$ uniformly,
the finitely many hole radii satisfy
\[
    r_z(t_n)\longrightarrow r_z(t_0)>0.
\]
Let $F_n\rightharpoonup F$ in $L^2(\mathbb B)$ and put
$u_n=B_{\delta,t_n}F_n$.
The energy estimate bounds the zero extensions $u_n$
in $H^1_0(\mathbb B)$.
After passing to a subsequence, we have
\[
    u_n\rightharpoonup u \quad\text{in }H^1_0(\mathbb B),
    \qquad
    u_n\longrightarrow u \quad\text{in }L^2(\mathbb B).
\]
Every compact subset of a limiting hole lies in the
corresponding hole for all sufficiently large $n$.
Hence $u$ vanishes in the limiting holes.
Since $u\in H^1(\mathbb B)$, its traces from the two sides of
each boundary circle agree and are therefore zero.
Together with the zero trace on $\partial\mathbb B$, this gives
\[
    u|_{\Omega_{\delta,t_0}}
    \in H^1_0(\Omega_{\delta,t_0}).
\]

For every $\phi\in C_c^\infty(\Omega_{\delta,t_0})$,
its support lies in $\Omega_{\delta,t_n}$ for all
sufficiently large $n$.
Passing to the limit in the weak equation gives
\[
    \int_{\Omega_{\delta,t_0}}
        \bigl(\nabla u\cdot\nabla\phi+u\phi\bigr)\,\mathrm dx
    =
    \int_{\Omega_{\delta,t_0}}F\phi\,\mathrm dx.
\]
Thus $u=B_{\delta,t_0}F$.
Uniqueness of the limit and Rellich compactness imply
\[
    B_{\delta,t_n}F_n
    \longrightarrow B_{\delta,t_0}F
    \quad\text{in }L^2(\mathbb B).
\]

If $\|B_{\delta,t_n}-B_{\delta,t_0}\|\not\to0$, then,
after passing to a subsequence, there exist unit vectors
$F_n\in L^2(\mathbb B)$ and $\varepsilon>0$ such that
\[
    \|(B_{\delta,t_n}-B_{\delta,t_0})F_n\|_{L^2(\mathbb B)}
    \ge\varepsilon.
\]
We may further assume that $F_n\rightharpoonup F$.
The preceding argument gives
$B_{\delta,t_n}F_n\to B_{\delta,t_0}F$.
Since $B_{\delta,t_0}$ is compact, we also have
$B_{\delta,t_0}F_n\to B_{\delta,t_0}F$,
a contradiction.
This proves norm-continuity.
\end{proof}

\begin{proposition}\label{net:realization}
There exists a disk with finitely many disjoint closed
interior disks removed whose Dirichlet Laplacian has an
eigenvalue of multiplicity $\dim E$, with exactly $n_-(\mathcal{M})$
eigenvalues below it, counted with multiplicity.
\end{proposition}
\begin{proof}
The basic idea is to apply Lemma~\ref{lem:stability} to complete the proof, with the choice
\[
    T=T(0),\qquad \mu=\kappa^{-1},
    \qquad m=\dim E,\qquad p=n_-(\mathcal{M}).
\]

First, we need to verify that the differential of $T(0)$ satisfies the assumption of Lemma~\ref{lem:stability}.
By the definition of $T(t)$ (see \eqref{Sec 7 Lem 1 equ-s}) and the identity $\mathcal{M}(0)=\mathcal{M}$, we have $T(0)=J_0(\mathcal{M}+\kappa I)^{-1}J_0^*$, $J_0^*J_0=I$,
$\ker \mathcal{M}=E$, and $\mathcal{M}+\kappa I\ge I$, we have
\[
    \ker\bigl(T(0)-\kappa^{-1}I\bigr)=J_0E,\qquad
    \#\{j:\lambda_j(T(0))>\kappa^{-1}\}=n_-(\mathcal{M}).
\]
By the surjectivity of  $R_E:\mathcal L(G)\longrightarrow\Sym(E)$ (see Lemma \ref{net:graph}), the right inverse $L$
chosen above satisfies
\[
    \langle L(h)x,y\rangle=h(x,y),
    \qquad h\in\Sym(E),\quad x,y\in E.
\]
here $m=\dim E$, $\dim \Sym(E)=\frac{m(m+1)}{2}=d_m$.  Let $P$ be the orthogonal projection onto
$\ker(T(0)-\kappa^{-1}I)=J_0E$. Since $PJ_0y=J_0y$ for every
$y\in E$, we have, for $h\in\Sym(E)$ and $x,y\in E$,
\[
\begin{aligned}
    \bigl\langle P\,DT(0)[h]J_0x,J_0y\bigr\rangle
    &=\bigl\langle DT(0)[h]J_0x,J_0y\bigr\rangle\\
    &=-\bigl\langle L(h)(\mathcal M+\kappa I)^{-1}x,
                        (\mathcal M+\kappa I)^{-1}y\bigr\rangle\\
    &=-\kappa^{-2}h(x,y).
\end{aligned}
\]
Thus, under the isometric identification $J_0:E\to J_0E$, the map
\[
    h\longmapsto P\,DT(0)[h]\big|_{J_0E}
\]
is multiplication by $-\kappa^{-2}$ on $\Sym(E)$ and hence is an
isomorphism.

Let $\tilde{r}$ be the radius of the parameter ball chosen above the formula \eqref{net:family} and let $r_0$
be given by Lemma~\ref{lem:stability} for this model family. Fix
\[
    r_1\le \tfrac12\min\{r_0,\tilde{r}\},\qquad \mathcal B=\overline B_{r_1}(0),
\]
and let $\delta_1(r_1)>0$ be given by Lemma \ref{lem:stability}.
Since $G$ is planar and $\mathcal{M}(t)$ is well signed on $\mathcal B$,
Lemma~\ref{lem:rooms} with $\delta_0=\delta_1(r_1)/2$
gives $\tau>0$ and potentials $V_t$ such that
\[
    A_t=\tau^{-1}(-\Delta_{\mathbb B}+V_t),\qquad
    A_t+\kappa\ge I,
\]
and 
\begin{align}\label{Sec 7 Lem-3 equ-d}
    \sup_{t\in\mathcal B}\|(A_t+\kappa)^{-1}-T(t)\|<\frac{\delta_1(r_1)}{2}.
\end{align}

Choose $K>\tau\kappa$ such that $V_t+K\ge1$ for all
$t\in\mathcal B$. Applying Lemma~\ref{lem:holes} to $V_t+K$ gives
\begin{align}\label{Sec 7 Lem 3 equ-a}
    \sup_{t\in\mathcal B}
    \left\|(-\Delta_{\Omega_{\delta,t}}+1)^{-1}
          -(-\Delta_{\mathbb B}+V_t+K+1)^{-1}\right\|
    \longrightarrow0
    \qquad(\delta\downarrow0). 
\end{align}

 We claim that 
 \begin{align}\label{Sec 7 Lem 3 equ-b}
     \sup_{t\in\mathcal B}
    \left\|
        \bigl[\tau^{-1}(-\Delta_{\Omega_{\delta,t}}-K)+\kappa\bigr]^{-1}
        -(A_t+\kappa)^{-1}
    \right\| \to 0 \qquad\mathrm{as}\qquad\delta\to 0.
 \end{align}
Put $c=K+1-\tau\kappa>0$ and write
\[
    B_t=(-\Delta_{\mathbb B}+V_t+K+1)^{-1},\qquad
    B_{\delta,t}=(-\Delta_{\Omega_{\delta,t}}+1)^{-1},
\]
where $B_{\delta,t}$ is extended by zero. Since
$B_t^{-1}=\tau(A_t+\kappa)+cI$ and $A_t+\kappa\ge I$, we have
\[
    I-cB_t\ge\frac{\tau}{\tau+c}I.
\]
By \eqref{Sec 7 Lem 3 equ-a},
$\sup_{t\in\mathcal B}\|B_{\delta,t}-B_t\|\to0$. Hence, for all
sufficiently small $\delta$,
\[
    I-cB_{\delta,t}\ge\frac{\tau}{2(\tau+c)}I
    \qquad\text{for every }t\in\mathcal B.
\]
Using $1-c=\tau\kappa-K$, we obtain
\[
\begin{aligned}
    &\bigl[\tau^{-1}(-\Delta_{\Omega_{\delta,t}}-K)+\kappa\bigr]^{-1}
       -(A_t+\kappa)^{-1} \\
    &\quad=\tau B_{\delta,t}(I-cB_{\delta,t})^{-1}
             -\tau B_t(I-cB_t)^{-1} \\
    &\quad=\tau(I-cB_{\delta,t})^{-1}
             (B_{\delta,t}-B_t)(I-cB_t)^{-1}.
\end{aligned}
\]
The preceding lower bounds therefore give
\[
\begin{aligned}
    &\sup_{t\in\mathcal B}
    \left\|
        \bigl[\tau^{-1}(-\Delta_{\Omega_{\delta,t}}-K)+\kappa\bigr]^{-1}
        -(A_t+\kappa)^{-1}
    \right\| \\
    &\qquad\le
    \frac{2(\tau+c)^2}{\tau}
    \sup_{t\in\mathcal B}\|B_{\delta,t}-B_t\|
    \longrightarrow0
    \qquad(\delta\downarrow0).
\end{aligned}
\]
With $\tau$ and $c$ fixed, choose $\delta>0$ sufficiently small that
\[
    \frac{2(\tau+c)^2}{\tau}
    \sup_{t\in\mathcal B}\|B_{\delta,t}-B_t\|
    <\frac{\delta_1(r_1)}{2}.
\]
Put
\[
    \widetilde T(t)
    =\bigl[\tau^{-1}(-\Delta_{\Omega_{\delta,t}}-K)+\kappa\bigr]^{-1}.
\]
By the preceding estimate and \eqref{Sec 7 Lem-3 equ-d},
\[
\begin{aligned}
    \sup_{t\in\mathcal B}\|\widetilde T(t)-T(t)\|
    &\le \sup_{t\in\mathcal B}\|\widetilde T(t)-(A_t+\kappa)^{-1}\|
       +\sup_{t\in\mathcal B}\|(A_t+\kappa)^{-1}-T(t)\|\\
    &<\frac{\delta_1(r_1)}{2}
       +\frac{\delta_1(r_1)}{2}
     =\delta_1(r_1).
\end{aligned}
\]

By Lemma~\ref{lem:holes}, the map $t\mapsto B_{\delta,t}$ is
norm-continuous for fixed $\delta$.
The uniform lower bound for $I-cB_{\delta,t}$ and the identity
\[
    \widetilde T(t)=\tau B_{\delta,t}(I-cB_{\delta,t})^{-1}
\]
therefore imply that $\widetilde T(t)$ is norm-continuous.
Applying Lemma~\ref{lem:stability}, there exists
$t_*\in\mathcal B$ such that
\[
    \dim\ker(\widetilde T(t_*)-\kappa^{-1} I)=m,\quad
\#\{j:\lambda_j(\widetilde T(t_*))>\kappa^{-1}\}=p.
\]
Thus $\widetilde T(t_*)$
has eigenvalue $\kappa^{-1}$ of multiplicity $\dim E$, with
exactly $n_-(\mathcal{M})$ eigenvalues above it.

For every Dirichlet eigenvalue $\lambda$ of $\Omega_{\delta,t_*}$,
the denominator $\tau^{-1}(\lambda-K)+\kappa$ is positive. Hence
\[
    \begin{aligned}
    \bigl[\tau^{-1}(\lambda-K)+\kappa\bigr]^{-1}=\kappa^{-1}
        &\iff \lambda=K,\\
    \bigl[\tau^{-1}(\lambda-K)+\kappa\bigr]^{-1}>\kappa^{-1}
        &\iff \lambda<K.
    \end{aligned}
\]
Since zero extension contributes only to the kernel, it follows that
\[
    \dim\ker(-\Delta_{\Omega_{\delta,t_*}}-K)=\dim E,\qquad
    \#\{j:\lambda_j(\Omega_{\delta,t_*})<K\}=n_-(\mathcal{M}).
\]The proof is complete.
\end{proof}

\subsection{Neumann and Robin boundary conditions}
\label{sec:robin-realization}

For a bounded smooth domain $\Omega\subset\R^2$ and
$\sigma\in\R$, let $H^\sigma_\Omega$ be the operator associated
with the quadratic form
\[
 u\longmapsto\int_\Omega|\nabla u|^2
   +\sigma\int_{\partial\Omega}u^2,
 \qquad u\in H^1(\Omega),\qquad \lambda_j^\sigma(\Omega)=\lambda_j(H^\sigma_\Omega).
\]
Thus $H^\sigma_\Omega$ has boundary condition
$\partial_\nu u+\sigma u=0$, and $\sigma=0$ gives the
Neumann Laplacian. The following approach is inspired by Colin de Verdi\`ere \cite{CdV}. For the convenience of the reader, we present the details.

\begin{proposition}\label{net:robin-realization}
For every fixed $\sigma\in\R$, if $n_-(\mathcal M)\ge1$,
there exists a bounded connected domain $\Omega\subset\R^2$
with $C^\infty$ boundary such that $H^\sigma_\Omega$
has an eigenvalue of multiplicity $\dim E$, with exactly
$n_-(\mathcal M)$ eigenvalues below it, counted with multiplicity.
\end{proposition}

\begin{proof}
\textbf{A weighted graph Laplacian.}
Let $G=(\mathcal V_G,\mathcal E)$ and $\mathcal M$
be as fixed at the beginning of this section,
with $N=|\mathcal V_G|$ and $n_-(\mathcal M)\ge1$. By the Perron--Frobenius theorem
\cite[p.~534]{HornJohnson}, its smallest eigenvalue is simple and admits an eigenvector with strictly positive entries. Since $n_-(\mathcal{M})\geq 1$, then its smallest eigenvalue is negative, we denote $a=-\lambda_{\min}(\mathcal{M})$. Write
\[
 \mathcal{M}h=-a h,\qquad a>0,\qquad h_i>0,
 \qquad \Lambda_h=\operatorname{diag}(h_1,\ldots,h_N).
\]
Retain the family $\mathcal{M}(t)=\mathcal{M}+L(t)$, here $L$ is the right inverse of $R_E$ as defined over \eqref{net:family}. After shrinking
$\mathcal B$, define
\begin{equation}\label{eq:neumann-vertex-masses}
 v_i(t)=-\frac{h_i(\mathcal{M}(t)h)_i}{a}>0,\qquad
 \Lambda_v(t)=\operatorname{diag}(v_1(t),\ldots,v_N(t)),
\end{equation}
and \begin{align*}
    c^{\mathcal{N}}_{ij}(t)=-\mathcal{M}(t)_{ij}h_i h_j>0
\end{align*}
for each edge $\{i,j\}\in\mathcal E$, and $v_i(0)=h_i^2$.
Let $L_{c^{\mathcal{N}}(t)}$ be the graph Laplacian defined by
\[
 x^TL_{c^{\mathcal{N}}(t)}x
 =\sum_{\{i,j\}\in\mathcal E}
       c^{\mathcal{N}}_{ij}(t)(x_i-x_j)^2.
\]
from the definition of \eqref{eq:neumann-vertex-masses}, we can get
\[
 L_{c^{\mathcal{N}}(t)}=\Lambda_h\mathcal{M}(t)\Lambda_h+a\Lambda_v(t).
\]
Consequently, the symmetric matrix
\[
 \mathcal M^{\mathcal{N}}(t)=\Lambda_v(t)^{-1/2}L_{c^{\mathcal{N}}(t)}
                         \Lambda_v(t)^{-1/2}
\]
is nonnegative and satisfies
\begin{equation}\label{eq:neumann-congruence}
 \mathcal M^{\mathcal{N}}(t)-aI
 =\Lambda_v(t)^{-1/2}\Lambda_h\mathcal{M}(t)\Lambda_h
                         \Lambda_v(t)^{-1/2},
 \qquad \mathcal M^{\mathcal{N}}(0)=\mathcal{M}+aI.
\end{equation}
Since $\mathcal M^{\mathcal{N}}(0)=\mathcal M+aI$ and
$E=\ker\mathcal M$, $a$ is the eigenvalue  of $\mathcal M^{\mathcal{N}}(0)$
with eigenspace $E$, and exactly $n_-(\mathcal M)$ eigenvalues lie
below it.
When \eqref{eq:neumann-congruence} is differentiated at zero
and restricted to $E\times E$, the derivatives of the two
outer factors give zero because $\mathcal{M}|_E=0$. Thus
\[
 D\mathcal M^{\mathcal{N}}(0)[s]|_{E\times E}
 =L(s)|_{E\times E}=s,
 \qquad s\in\Sym(E).
\]
Set
$
 T^{\mathcal{N}}(t)=(\mathcal M^{\mathcal{N}}(t)+I)^{-1}.
$
Then $T^{\mathcal{N}}(0)$ has eigenvalue $(a+1)^{-1}$ with eigenspace $E$,
exactly $n_-(\mathcal{M})$ eigenvalues above it, and
\[
 DT^{\mathcal{N}}(0)[s]|_{E\times E}=-(a+1)^{-2}s.
\]
Thus $T^{\mathcal{N}}(t)$ satisfies the transversality assumption of Lemma~\ref{lem:stability}.

\smallskip
\textbf{The Neumann approximation.}
Set $m=\dim E$, $p=n_-(\mathcal M)$ and $\mu=(a+1)^{-1}$.
In the planar construction of \cite[\S6(d), pp.~610--611]{CdV},
let $R_{i,\varepsilon}$ be the room at vertex $i$ and $T_{e,\varepsilon}$ the tube along
$e=\{i,j\}$. Denote their interfaces by $\Gamma_{ei,\varepsilon}$ and
$\Gamma_{ej,\varepsilon}$ as in Figure~\ref{fig:sec72-neumann}.

% Preamble: \usepackage{tikz,graphicx,capt-of}
% Insert between paragraphs, including between paragraphs of a proof.
\par\smallskip
\noindent\begin{minipage}{\linewidth}
 \centering
 \begin{minipage}[c]{.52\linewidth}\centering
 \resizebox{\linewidth}{!}{%
  \begin{tikzpicture}[font=\normalsize,>=stealth,line cap=round,line join=round,
      wall/.style={draw=black,line width=.85pt}]
    % Only exposed arcs of the rooms are boundary arcs. Tube mouths stay open.
    \def\NRfigureRoom#1#2#3#4{%
      \coordinate (#1) at #2;
      \begin{scope}[shift={(#1)}]
        \fill[black!7] (0,0) circle[radius=#3];
        \foreach \lo/\hi in {#4}{
          \draw[wall] (\lo:#3) arc[start angle=\lo,end angle=\hi,radius=#3];}
      \end{scope}}
    \NRfigureRoom{I}{(0,0)}{.64}{29/106,164/196,254/331}
    \NRfigureRoom{J}{(2.07,0)}{.61}{29/151,209/331}
    \NRfigureRoom{A}{(135:1.60)}{.40}{209/286,344/511}
    \NRfigureRoom{B}{(225:1.60)}{.42}{74/376}
    \NRfigureRoom{C}{(3.62,0)}{.40}{29/151,209/331}

    % Start room, direction, centre distance, and the two room radii.
    \foreach \site/\az/\dist/\ra/\rb in
      {I/0/2.07/.64/.61,I/135/1.60/.64/.40,
       I/225/1.60/.64/.42,J/0/1.55/.61/.40}{
      \begin{scope}[shift={(\site)},rotate=\az]
        \pgfmathsetmacro{\xl}{\ra*cos(29)}
        \pgfmathsetmacro{\yl}{\ra*sin(29)}
        \pgfmathsetmacro{\xr}{\dist-\rb*cos(29)}
        \pgfmathsetmacro{\yr}{\rb*sin(29)}
        \pgfmathsetmacro{\xm}{(\xl+\xr)/2}
        \def\topwall{(\xl,\yl)
          .. controls +(.07,-.126) and +(-.15,0) .. (\xm,.067)
          .. controls +(.15,0) and +(-.07,-.126) .. (\xr,\yr)}
        \def\bottomwall{(\xr,-\yr)
          .. controls +(-.07,.126) and +(.15,0) .. (\xm,-.067)
          .. controls +(-.15,0) and +(.07,.126) .. (\xl,-\yl)}
        \fill[black!7] \topwall -- \bottomwall -- cycle;
        \draw[wall] \topwall;
        \draw[wall] \bottomwall;
      \end{scope}}

    % A short open tube followed by dots indicates the omitted network.
    \foreach \site/\az in {A/180,C/0}{
      \begin{scope}[shift={(\site)},rotate=\az]
        \pgfmathsetmacro{\xl}{.40*cos(29)}
        \pgfmathsetmacro{\yl}{.40*sin(29)}
        \def\topwall{(\xl,\yl) .. controls +(.07,-.126)
                      and +(-.14,0) .. (.71,.055)}
        \def\bottomwall{(.71,-.055) .. controls +(-.14,0)
                         and +(.07,.126) .. (\xl,-\yl)}
        \fill[black!7] \topwall -- \bottomwall -- cycle;
        \draw[wall] \topwall;
        \draw[wall] \bottomwall;
        \foreach \xx in {.81,.94,1.07}
          \fill (\xx,0) circle[radius=.022];
      \end{scope}}

    \node at (I) {$R_{i,\varepsilon}$};
    \node at (J) {$R_{j,\varepsilon}$};
    \node (edge) at (1.06,-.71) {$T_{e,\varepsilon}$};
    \draw[thin,->] (edge.north) -- (1.06,-.08);
    \node (value) at (.88,1.28) {$y_i=x_i/\sqrt{v_i(t)}$};
    \draw[thin,->] (value.south west) -- (.16,.36);
    \node at (.90,-1.78) {$D_{\varepsilon,t}$};
  \end{tikzpicture}}
 \end{minipage}\hfill
 \begin{minipage}[c]{.44\linewidth}\centering
 \resizebox{\linewidth}{!}{%
  \begin{tikzpicture}[font=\small,>=stealth,line join=round,line cap=round]
    \def\NRx{1.10}
    \pgfmathsetmacro{\NRy}{.38*sqrt(\NRx*\NRx+.04)}
    % The lateral boundary arcs are graphs of +/- .38*sqrt(x^2+.04).
    \path[fill=black!7,draw=black,line width=.7pt]
      plot[domain=-\NRx:\NRx,samples=41] (\x,{.38*sqrt(\x*\x+.04)})
      .. controls (1.38,.53) and (1.48,.82) .. (1.91,.82)
      .. controls (2.92,.82) and (2.92,-.82) .. (1.91,-.82)
      .. controls (1.48,-.82) and (1.38,-.53) .. (\NRx,-\NRy)
      -- plot[domain=\NRx:-\NRx,samples=41] (\x,{-.38*sqrt(\x*\x+.04)})
      .. controls (-1.38,-.53) and (-1.48,-.82) .. (-1.91,-.82)
      .. controls (-2.92,-.82) and (-2.92,.82) .. (-1.91,.82)
      .. controls (-1.48,.82) and (-1.38,.53) .. cycle;
    \foreach \side in {-1,1}{
      \draw[black!65,densely dashed]
        ({\side*\NRx},\NRy) .. controls ({\side*1.20},.16)
          and ({\side*1.20},-.16) .. ({\side*\NRx},-\NRy);}
    \node at (-1.92,0) {$R_{i,\varepsilon}$};
    \node at (1.92,0) {$R_{j,\varepsilon}$};
    \node (gi) at (-1.11,-1.17) {$\Gamma_{ei,\varepsilon}$};
    \node (gj) at ( 1.11,-1.17) {$\Gamma_{ej,\varepsilon}$};
    \draw[thin,->] (gi.north) -- (-1.17,-.24);
    \draw[thin,->] (gj.north) -- ( 1.17,-.24);
    \node (edgeRegion) at (0,-.78) {$T_{e,\varepsilon}$};
    \draw[thin,->] (edgeRegion.north) -- (0,-.085);
    \draw[->,black] (0,.076) -- (0,.76) node[right] {$\nu_{\varepsilon}$};
    \node at (0,1.22) {$\partial_{\nu_{\varepsilon}}u=0$};
    \node at (0,-1.78) {$e=\{i,j\}$};
  \end{tikzpicture}}
 \end{minipage}
 \par\captionof{figure}{}
 \label{fig:sec72-neumann}
\end{minipage}\par\smallskip
For each $(\varepsilon,t)$, let $D_{\varepsilon,t}$ be the planar
domain obtained by joining the rooms $R_{i,\varepsilon}$ to the tubes $T_{e,\varepsilon}$
along the interfaces $\Gamma_{ei,\varepsilon}$. Thus
\[
  D_{\varepsilon,t}
  =\operatorname{int}_{\mathbb R^2}\!
    \left(
      \bigcup_{i=1}^{N}\overline{R_{i,\varepsilon}}
      \,\cup\!
      \bigcup_{e\in\mathcal E}\overline{T_{e,\varepsilon}}
    \right).
\]
 We give a precise parametrization of the tube $T_{e,\varepsilon}$ using a smooth variant of the conformal rectangle in \cite[Proposition~6.6, pp.~607--608]{CdV}. 
Fix $0<\theta<\pi/2$ and choose $d>0$ sufficiently small.
In local Euclidean coordinates $Z=X+\mathrm iY$, set
\[
    \ell_{\varepsilon}=\frac{\theta}{\varepsilon c^{\mathcal{N}}_{ij}(t)},
    \qquad
    \rho_{\varepsilon}=\frac{d}{\sinh\ell_{\varepsilon}},
    \qquad
    T_{e,\varepsilon}=F_{e,\varepsilon}(\mathcal R_{\ell_{\varepsilon},\theta}),
\]
where
\[
    \mathcal R_{\ell_{\varepsilon},\theta}
    =\{s+\mathrm i\phi:|s|<\ell_{\varepsilon},\ |\phi|<\theta\},
    \qquad
    F_{e,\varepsilon}(w)=\rho_{\varepsilon}\sinh w.
\]
Since $\operatorname{Re}F_{e,\varepsilon}'(s+\mathrm i\phi)
=\rho_{\varepsilon}\cosh s\cos\phi>0$ on $\overline{\mathcal R_{\ell_{\varepsilon},\theta}}$,
integration of $F_{e,\varepsilon}'$ along line segments shows that $F_{e,\varepsilon}$ is
injective on the closed rectangle and conformal in its interior.
In these coordinates,
\[
    X=\rho_{\varepsilon}\sinh s\cos\phi,
    \qquad
    Y=\rho_{\varepsilon}\cosh s\sin\phi.
\]

Thus the sides $s=-\ell_{\varepsilon}$ and $s=\ell_{\varepsilon}$ give the interfaces with
$R_{i,\varepsilon}$ and $R_{j,\varepsilon}$, respectively, parametrized by
\[
    \phi\longmapsto
    \bigl(\pm d\cos\phi,\sqrt{d^2+\rho_{\varepsilon}^2}\sin\phi\bigr),
    \qquad |\phi|\le\theta.
\]
The sides $\phi=\pm\theta$ give the two boundary arcs of
$D_{\varepsilon,t}$ along $T_{e,\varepsilon}$,
\[
    Y=\pm\tan\theta\sqrt{X^2+\rho_{\varepsilon}^2\cos^2\theta},
    \qquad |X|\le d\cos\theta.
\]
At the endpoints of the lateral arcs, $|X|=d\cos\theta>0$.
In fixed neighborhoods of these points, the arcs extend smoothly
across the interfaces and are joined to the room boundaries by a
smooth cutoff; see Figure \ref{fig:sec72-neumann} and Figure~\ref{fig:tube-parametrization}.
\begin{figure}[H]
\centering
\resizebox{\linewidth}{!}{%
\begin{tikzpicture}[
    x=1cm,y=1cm,font=\small,>=stealth,
    wall/.style={draw=blue!55!black,line width=1.10pt},
    interface/.style={draw=orange!78!black,line width=1.05pt,
                      dash pattern=on 4pt off 2pt},
    section/.style={draw=violet!80!black,line width=1pt,
                    dash pattern=on 3pt off 1.3pt on .7pt off 1.3pt},
    centerline/.style={draw=black!55,thin,densely dotted},
    gridline/.style={draw=blue!18!gray,opacity=.43,thin},
    extension/.style={draw=black!40,thin},
    heading/.style={font=\normalsize},
    note/.style={font=\footnotesize,inner sep=2pt},
    symbol/.style={font=\small,inner sep=2pt}
]
  % Values for the illustration. PGF's sin and cos take degrees.
  \def\ellval{2.05}
  \def\thetadeg{38}
  \def\dval{2.75}
  \pgfmathsetmacro{\rhoval}{\dval/sinh(\ellval)}
  \pgfmathsetmacro{\ymid}{\rhoval*sin(\thetadeg)}
  \pgfmathsetmacro{\rectx}{1.18*\ellval}
  \pgfmathsetmacro{\recty}{1.18*rad(\thetadeg)}
  \tikzset{plotcurve/.style={samples=55,smooth,variable=\u}}
  % Plane coordinates of F(s+i phi); the second argument is in degrees.
  \tikzset{declare function={
    FX(\s,\a)=\rhoval*sinh(\s)*cos(\a);
    FY(\s,\a)=\rhoval*cosh(\s)*sin(\a);
  }}
  \node[heading] at (0,2.47) {Parameter rectangle $\mathcal R_{\ell_{\varepsilon},\theta}$};
  \node[heading] at (8.70,2.47) {Tube $T_{e,\varepsilon}$ in the $(X,Y)$-plane};
  % Left panel: s increases to the right and phi increases upwards.
  \fill[blue!4] (-\rectx,-\recty) rectangle (\rectx,\recty);
  \foreach \k in {-.55,.55}{
      \draw[gridline] ({\k*\rectx},-\recty)--({\k*\rectx},\recty);
      \draw[gridline] (-\rectx,{\k*\recty})--(\rectx,{\k*\recty});}
  \draw[centerline] (-\rectx,0)--(\rectx,0);
  \draw[section] (0,-\recty)--(0,\recty);
  \foreach \k in {-1,1}{
      \draw[wall] (-\rectx,{\k*\recty})--(\rectx,{\k*\recty});
      \draw[interface] ({\k*\rectx},-\recty)--({\k*\rectx},\recty);}
  \node[symbol,above] at (0,\recty) {$\phi=\theta$};
  \node[symbol,below] at (0,-\recty) {$\phi=-\theta$};
  \node[symbol,left] at (-\rectx,0) {$s=-\ell_{\varepsilon}$};
  \node[symbol,right] at (\rectx,0) {$s=\ell_{\varepsilon}$};
  \node[note,fill=blue!4,text=violet!80!black] at (0,.32) {$s=0$};
  \node[note] at (0,-1.67) {$w=s+\mathrm i\phi$};
  \draw[->,semithick] (3.53,0)--(5.02,0)
      node[midway,above=4pt,symbol] {$F_{e,\varepsilon}$};
  \node[note] at (4.27,-.42) {$w\mapsto\rho_{\varepsilon}\sinh w$};
  \begin{scope}[shift={(8.70,0)}]
    % The sides phi=+-theta are hyperbolic arcs; the ends s=+-ell are elliptic arcs.
    \path[fill=blue!4] ({FX(-\ellval,-\thetadeg)},{FY(-\ellval,-\thetadeg)})
      -- plot[plotcurve,domain=-\ellval:\ellval] ({FX(\u,-\thetadeg)},{FY(\u,-\thetadeg)})
      -- plot[plotcurve,domain=-\thetadeg:\thetadeg] ({FX(\ellval,\u)},{FY(\ellval,\u)})
      -- plot[plotcurve,domain=\ellval:-\ellval] ({FX(\u,\thetadeg)},{FY(\u,\thetadeg)})
      -- plot[plotcurve,domain=\thetadeg:-\thetadeg] ({FX(-\ellval,\u)},{FY(-\ellval,\u)})
      -- cycle;
    \foreach \k in {-.55,.55}{
      \draw[gridline] plot[plotcurve,domain=-\thetadeg:\thetadeg]
                          ({FX(\k*\ellval,\u)},{FY(\k*\ellval,\u)});
      \draw[gridline] plot[plotcurve,domain=-\ellval:\ellval]
                          ({FX(\u,\k*\thetadeg)},{FY(\u,\k*\thetadeg)});}
    \draw[centerline] (-\dval,0)--(\dval,0);
    \foreach \k in {-1,1}{
      \draw[wall] plot[plotcurve,domain=-\ellval:\ellval]
                  ({FX(\u,\k*\thetadeg)},{FY(\u,\k*\thetadeg)});
      \draw[interface] plot[plotcurve,domain=-\thetadeg:\thetadeg]
                  ({FX(\k*\ellval,\u)},{FY(\k*\ellval,\u)});}
    \draw[section] (0,-\ymid)--(0,\ymid);
    \foreach \k in {-1,1}{
       \fill[violet!80!black] (0,{\k*\ymid}) circle[radius=1.15pt];
       \fill[black!65] ({\k*\dval},0) circle[radius=1pt];}
    \node[symbol] at (-3.23,0) {$R_{i,\varepsilon}$};
    \node[symbol] at ( 3.23,0) {$R_{j,\varepsilon}$};
    \node[note,text=orange!78!black] at (-3.00,1.13) {$\Gamma_{ei,\varepsilon}$};
    \node[note,text=orange!78!black] at ( 3.00,1.13) {$\Gamma_{ej,\varepsilon}$};
    \node[note] (widthlabel) at (0,1.26) {width $2\rho_{\varepsilon}\sin\theta$};
    \draw[->,thin,black!60] (widthlabel.south)--(0,{\ymid+.045});
    \node[note,below right] at (-\dval,-.03) {$(-d,0)$};
    \node[note,below left] at (\dval,-.03) {$(d,0)$};
    \node[note] at (0,-1.25) {$X=0$};
    \draw[->,thin,black!60] (0,-1.10)--(0,{-\ymid-.045});
    \draw[extension] (-\dval,-1.90)--(-\dval,-2.14);
    \draw[extension] (\dval,-1.90)--(\dval,-2.14);
    \draw[<->,thin,black!70] (-\dval,-2.03)--(\dval,-2.03)
       node[midway,fill=white,note] {centerline length $2d$};
  \end{scope}
  % Line styles distinguish all three kinds of curves in grayscale as well.
  \foreach \x/\sty/\lbl in {
       -1.60/wall/lateral boundary,
        3.14/interface/room interface,
        7.65/section/central section}{
    \draw[\sty] (\x,-2.83)--++(.70,0);
    \node[note,anchor=west] at ({\x+.81},-2.83) {\lbl};}
\end{tikzpicture}%
}
\caption{}
\label{fig:tube-parametrization}
\end{figure}

To compute the contribution of $T_{e,\varepsilon}$ to the quadratic form,
prescribe constant traces $y_i$ and $y_j$ on its two interfaces.
On $\mathcal R_{\ell_{\varepsilon},\theta}$ the Dirichlet integral with these
traces is minimized by
\[
    \chi_{\varepsilon}(s,\phi)=y_i+\frac{s+\ell_{\varepsilon}}{2\ell_{\varepsilon}}(y_j-y_i),
\]
since $\chi_{\varepsilon}$ is harmonic and $\partial_\phi\chi_{\varepsilon}=0$ on $\phi=\pm\theta$.
By conformal invariance of the planar Dirichlet integral,
$u_{\varepsilon}=\chi_{\varepsilon}\circ F_{e,\varepsilon}^{-1}$ is the corresponding minimizer on $T_{e,\varepsilon}$, and
\[
    \int_{T_{e,\varepsilon}}|\nabla u_{\varepsilon}|^2\,\mathrm dX\,\mathrm dY
    =\int_{-\theta}^{\theta}\int_{-\ell_{\varepsilon}}^{\ell_{\varepsilon}}
       \left|\frac{y_j-y_i}{2\ell_{\varepsilon}}\right|^2
       \,\mathrm ds\,\mathrm d\phi
    =\frac{\theta}{\ell_{\varepsilon}}(y_i-y_j)^2
    =\varepsilon c^{\mathcal{N}}_{ij}(t)(y_i-y_j)^2.
\]

Cut each $T_{e,\varepsilon}$ at $s=0$ and denote its half adjacent to $R_{i,\varepsilon}$ by
$T_{ei,\varepsilon}$. Let $S_{i,\varepsilon}$ be $R_{i,\varepsilon}$ with the halves $T_{ei,\varepsilon}$ attached along
$\Gamma_{ei,\varepsilon}$. Since
$|F_{e,\varepsilon}'(s+\mathrm i\phi)|^2=\rho_{\varepsilon}^2(\sinh^2s+\cos^2\phi)$,
\begin{align}\label{Sec 7.2 equ-a}
 |T_{ei,\varepsilon}|=\int_0^{\ell_{\varepsilon}}\int_{-\theta}^{\theta}
       |F_{e,\varepsilon}'(s+\mathrm i\phi)|^2\,\mathrm d\phi\,\mathrm ds
 =\theta d^2\coth\ell_{\varepsilon}+\tfrac12\rho_{\varepsilon}^2\ell_{\varepsilon}\sin(2\theta)
 \longrightarrow\theta d^2 .
\end{align}
Let $\deg_G(i)$ denote the number of edges incident to the vertex $i$. Scale a smooth planar embedding of the finite graph so that the
vertices have disjoint neighborhoods large enough to hold chambers
of area greater than $v_i(0)$. Take one small neck chart per edge
and join its two ends to the respective vertex chambers by disjoint
thin strips along the embedded half-edges. By first decreasing $d$
and then the strip widths, the limiting collars and strips at each
vertex can have total area less than
$v_i(0)-\deg_G(i)\theta d^2$. Smooth the joins away from the neck
charts. The chamber at each vertex can now be enlarged continuously,
keeping its attachment fixed, until the limiting room has area
$v_i(0)-\deg_G(i)\theta d^2$. These adjustments are independent.
For positive $\varepsilon$ retain these rooms away from the interface
collars and use the smoothly converging collars prescribed above.
Thus
\begin{align*}
    \lim_{\varepsilon\to 0}|R_{i,\varepsilon}(0)|=v_i(0)-\deg_G(i)\theta d^2,
\end{align*}
where $v_i(t)$ is given by \eqref{eq:neumann-vertex-masses}. Therefore, we can get $|S_{i,\varepsilon}(0)|\to v_i(0)$ as $\varepsilon \to0$. To get the uniform convergence, write $R_{i,\varepsilon}(0)$ locally as
$\{\zeta<f_i(\xi),\ \xi\in I_i\}$ in Euclidean coordinates
away from the interfaces. On disjoint such charts replace $f_i$ by
\[
 f_{i,t}(\xi)=f_i(\xi)+\bigl(v_i(t)-v_i(0)\bigr)\psi_i(\xi),
 \qquad \psi_i\in C_c^\infty(I_i),\qquad
 \int_{I_i}\psi_i(\xi)\,\mathrm d\xi=1.
\]
In the chart $U_i$, define $R_{i,\varepsilon}(t)$ by $\{\zeta<f_{i,t}(\xi),\ \xi\in I_i\}$. Then
\[
 |R_{i,\varepsilon}(t)|-|R_{i,\varepsilon}(0)|=v_i(t)-v_i(0)+o(1).
\]
Outside $U_i$, keep the room fixed except in the interface collars.
On a sufficiently small ball $\mathcal B=\overline B_{r_1}(0)$,
$\rho_{\varepsilon}\to0$ uniformly in $t$, so the collars converge smoothly
to their $t$-independent limits. Since
$|S_{i,\varepsilon}(t)|=|R_{i,\varepsilon}(t)|+\sum_{e\ni i}|T_{ei,\varepsilon}(t)|$
and the convergence in \eqref{Sec 7.2 equ-a} is uniform in $t$, we obtain
\[
 \sup_{t\in\mathcal B}\bigl|\,|S_{i,\varepsilon}(t)|-v_i(t)\bigr|\longrightarrow0
 \qquad\varepsilon \to 0.
\]

% Compile with pdfLaTeX. The tikzpicture can also be pasted into a paper
% whose preamble contains \usepackage{tikz}.

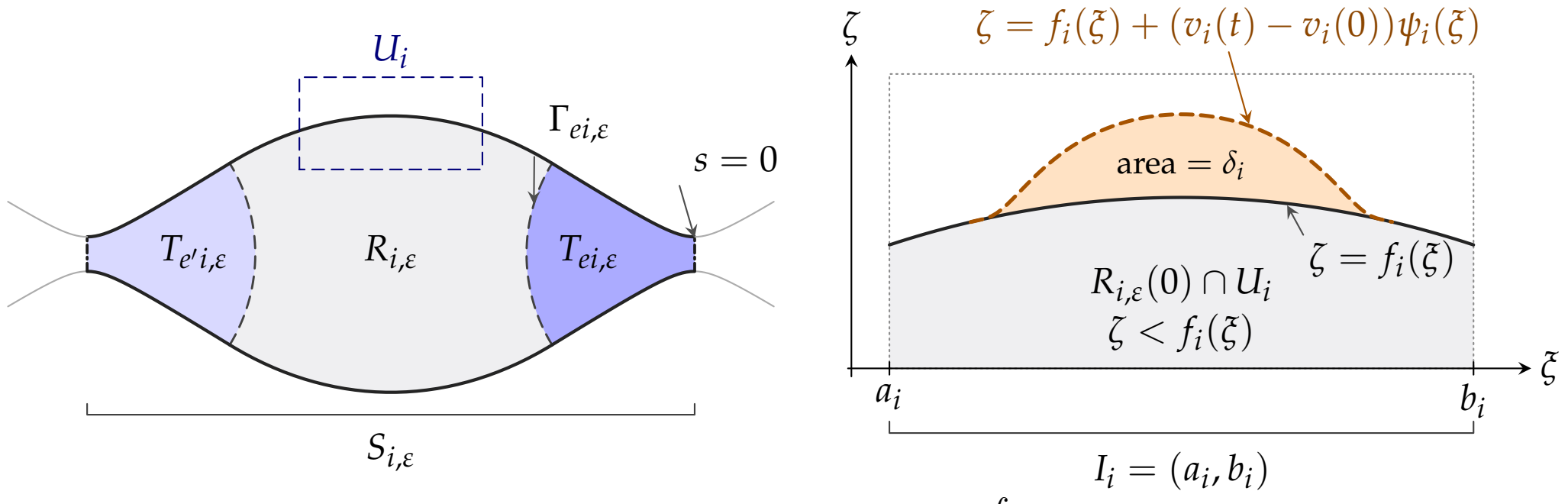
\begin{tikzpicture}[
  x=1cm,y=1cm,font=\small,line cap=round,line join=round,>=stealth,
  wall/.style={draw=black!86,line width=.8pt},
  interface/.style={draw=black!75,densely dashed,line width=.54pt},
  cut/.style={draw=black!90,dash pattern=on 2pt off 1pt on .5pt off 1pt,
              line width=.7pt},
  leader/.style={draw=black!70,thin,->},
  guide/.style={draw=black!45,thin,densely dotted},
  chart/.style={draw=blue!48!black,thin,dash pattern=on 3pt off 1.8pt},
  moved/.style={draw=orange!65!black,line width=.95pt,densely dashed},
  declare function={
    halfwidth(\x)=tan(32)*sqrt(\x*\x+(.30*cos(32))^2);
  }]
  \def\cutdist{2.80}
  \def\tubed{1.55}
  \pgfmathsetmacro{\joinu}{\tubed*cos(32)}
  \pgfmathsetmacro{\joinx}{\cutdist-\joinu}
  \pgfmathsetmacro{\joiny}{halfwidth(\joinu)}
  \pgfmathsetmacro{\ellipsey}{sqrt(\tubed^2+.30^2)}
  \pgfmathsetmacro{\slope}{tan(32)*\joinu/sqrt(\joinu^2+(.30*cos(32))^2)}
  \pgfmathsetmacro{\controlx}{\joinx-.95}
  \pgfmathsetmacro{\controly}{\joiny+.95*\slope}
  % The midpoint of the upper cubic wall is the center of the selected arc.
  % p in [-1,1] is its centered Bezier parameter. Select |p| < 0.56.
  \def\part{.56}
  \pgfmathdeclarefunction{arcX}{1}{\pgfmathparse{
    (3*(\joinx+\controlx)*#1+(\joinx-3*\controlx)*(#1)^3)/4}}
  \pgfmathdeclarefunction{arcY}{1}{\pgfmathparse{
    \joiny+.75*(\controly-\joiny)*(1-(#1)^2)}}

  % LEFT: original cell and the local chart away from the tube interfaces.
  \begin{scope}[scale=.94]
    \def\roomtop{(-\joinx,\joiny) .. controls
      (-\controlx,\controly) and (\controlx,\controly) .. (\joinx,\joiny)}
    \def\roombottom{(\joinx,-\joiny) .. controls
      (\controlx,-\controly) and (-\controlx,-\controly) .. (-\joinx,-\joiny)}
    \path[fill=blue!7!gray!12] \roomtop -- \roombottom -- cycle;
    \draw[wall] \roomtop;
    \draw[wall] \roombottom;
    \foreach \pos/\dir/\shade in
      {-\cutdist/1/blue!15!white,\cutdist/-1/blue!33!white}{
      \begin{scope}[shift={(\pos,0)},xscale=\dir]
        \path[fill=\shade]
          plot[domain=0:\joinu,samples=35] (\x,{halfwidth(\x)})
          -- plot[domain=32:-32,samples=35]
               ({\tubed*cos(\x)},{\ellipsey*sin(\x)})
          -- plot[domain=\joinu:0,samples=35] (\x,{-halfwidth(\x)}) -- cycle;
        \draw[interface] plot[domain=-32:32,samples=35]
          ({\tubed*cos(\x)},{\ellipsey*sin(\x)});
        \foreach \sgn in {-1,1}{
          \draw[wall] plot[domain=0:\joinu,samples=35]
            (\x,{\sgn*halfwidth(\x)});
          \draw[gray!65,thin] plot[domain=-.73:0,samples=24]
            (\x,{\sgn*halfwidth(\x)});
        }
        \draw[cut] (0,{-halfwidth(0)}) -- (0,{halfwidth(0)});
      \end{scope}
    }
    \node at (0,0) {$R_{i,\varepsilon}$};
    \node at (-1.82,0) {$T_{e'i,\varepsilon}$};
    \node at ( 1.82,0) {$T_{ei,\varepsilon}$};
    \node (gamma) at (1.74,1.20) {$\Gamma_{ei,\varepsilon}$};
    \draw[leader] (gamma.south west) --
      ({\cutdist-\tubed*cos(17)},{\ellipsey*sin(17)});
    \node (mid) at (3.18,.88) {$s=0$};
    \draw[leader] (mid.south west) -- (\cutdist,.15);
    \draw[black!75,thin] (-\cutdist,-1.38) -- ++(0,-.10)
                      -- (\cutdist,-1.48) -- ++(0,.10);
    \node[below=1pt] at (0,-1.48) {$S_{i,\varepsilon}$};
    \draw[chart] ({arcX(-\part)},.78) rectangle ({arcX(\part)},1.63);
    \node[above,text=blue!48!black,inner sep=2pt] at (0,1.63) {$U_i$};
  \end{scope}

  % RIGHT: enlarge exactly the bounded chart from the left-hand panel.
  % The local axes use a translated origin. Both coordinates are magnified
  % by the same factor. arcX is increasing on the displayed interval.
  \begin{scope}[shift={(3.96,0)}]
    \pgfmathdeclarefunction{cx}{1}{\pgfmathparse{2.89+3*arcX(\part*#1)}}
    \pgfmathdeclarefunction{cy}{1}{\pgfmathparse{3*(arcY(\part*#1)-1.11)}}
    % A compactly supported bump, shown with positive area increment.
    % The cutoff in the exponential only avoids numerical overflow in PGF.
    \pgfmathdeclarefunction{bump}{1}{\pgfmathparse{
      abs(#1)<.72 ? .72*exp(max(-9,1-1/max(.01,1-(#1/.72)^2))) : 0}}
    \pgfmathsetmacro{\xa}{cx(-1)}
    \pgfmathsetmacro{\xb}{cx(1)}
    \pgfmathsetmacro{\bottomchart}{3*(.78-1.11)}
    \pgfmathsetmacro{\topchart}{3*(1.63-1.11)}

    % The original room lies below the graph in U_i.
    \path[fill=blue!7!gray!12]
      plot[domain=-1:1,samples=65] ({cx(\x)},{cy(\x)})
      -- (\xb,\bottomchart) -- (\xa,\bottomchart) -- cycle;
    \draw[guide] (\xa,\bottomchart) rectangle (\xb,\topchart);
    % The added strip lies strictly inside the horizontal interval.
    \path[fill=orange!23]
      plot[domain=-.72:.72,samples=81] ({cx(\x)},{cy(\x)+bump(\x)})
      -- plot[domain=.72:-.72,samples=61] ({cx(\x)},{cy(\x)}) -- cycle;
    \draw[wall] plot[domain=-1:1,samples=65] ({cx(\x)},{cy(\x)});
    \draw[moved] plot[domain=-.72:.72,samples=81] ({cx(\x)},{cy(\x)+bump(\x)});

    \draw[->,thin] (0,\bottomchart) -- ({\xb+.50},\bottomchart)
      node[right,inner sep=2pt] {$\xi$};
    \draw[->,thin] (.03,{\bottomchart-.07}) -- (.03,{\topchart+.17})
      node[above,inner sep=2pt] {$\zeta$};
    \foreach \xx/\lab in {\xa/a_i,\xb/b_i}{
      \draw[thin] (\xx,{\bottomchart+.055}) -- (\xx,{\bottomchart-.055})
        node[below,inner sep=2pt] {$\lab$};
    }
    \draw[black!75,thin] (\xa,{\bottomchart-.46}) -- ++(0,-.10)
                 -- (\xb,{\bottomchart-.56}) -- ++(0,.10);
    \node[below=1pt] at (2.89,{\bottomchart-.56}) {$I_i=(a_i,b_i)$};

    \node[inner sep=1pt] at (2.89,-.27) {$R_{i,\varepsilon}(0)\cap U_i$};
    \node[inner sep=1pt] at (2.89,-.68) {$\zeta<f_i(\xi)$};
    \node[font=\footnotesize,inner sep=1pt] at (2.89,.77)
      {$\mathrm{area}=\delta_i$};
    \node[inner sep=1.2pt,anchor=east] (orig) at ({\xb-.1},-.05)
      {$\zeta=f_i(\xi)$};
    \draw[leader] (orig.north west) -- ({cx(.36)},{cy(.36)-.015});
    \node[inner sep=1.3pt,anchor=south,text=orange!55!black] (new) at (3.30,1.75)
      {$\zeta=f_i(\xi)+(v_i(t)-v_i(0))\psi_i(\xi)$};
    \draw[leader,draw=orange!65!black] (new.south) --
      ({cx(.23)},{cy(.23)+bump(.23)+.015});
  \end{scope}
  \node[font=\footnotesize] at (3.34,-2.45)
    {$\delta_i=v_i(t)-v_i(0)>0,\qquad
      \psi_i\in C_c^\infty(I_i),\qquad
      \displaystyle\int_{I_i}\psi_i(\xi)\,\mathrm{d}\xi=1.$};
\end{tikzpicture}

Write $D_{\varepsilon}=D_{\varepsilon,t}$ and
$q_{D_{\varepsilon}}^0(u,g)=\int_{D_{\varepsilon}}\nabla u\cdot\nabla g\,\mathrm dz$.
For $x\in\R^N$, set $y=\Lambda_v(t)^{-1/2}x$. Define the linear map
$J^{\mathcal{N}}_{\varepsilon}=J^{\mathcal{N}}_{\varepsilon,t}:\R^N\to H^1(D_{\varepsilon})$ by
\begin{align}\label{Sec 7 f-12}
 (J^{\mathcal{N}}_{\varepsilon}x)(z)=
 \begin{cases}
   y_i,
     & z\in R_{i,\varepsilon},\\[4pt]
   \displaystyle
   \frac{\ell_{\varepsilon}-s}{2\ell_{\varepsilon}}\,y_i
   +\frac{\ell_{\varepsilon}+s}{2\ell_{\varepsilon}}\,y_j,
     & z=F_{e,\varepsilon}(s+\mathrm i\phi)\in T_{e,\varepsilon},
       \quad e=\{i,j\}.
 \end{cases}
\end{align}
At $s=-\ell_{\varepsilon}$ and $s=\ell_{\varepsilon}$, the tube
formula gives $y_i$ and $y_j$, $J^{\mathcal{N}}_{\varepsilon}x $ is continuous on $D_{\epsilon}$ . By conformal invariance of the Dirichlet energy,
\[
 \int_{T_{e,\varepsilon}} |\nabla(J^{\mathcal{N}}_{\varepsilon}x)|^2\,\mathrm dz
 =\int_{-\ell_{\varepsilon}}^{\ell_{\varepsilon}}
   \int_{-\theta}^{\theta}
   \left(\frac{y_j-y_i}{2\ell_{\varepsilon}}\right)^2
   \,\mathrm d\phi\,\mathrm ds
 =\frac{\theta}{\ell_{\varepsilon}}(y_i-y_j)^2
 =\varepsilon c^{\mathcal{N}}_{ij}(t)(y_i-y_j)^2,
\]
where $\ell_{\varepsilon}=\theta/(\varepsilon c^{\mathcal{N}}_{ij}(t))$.
Using $y=\Lambda_v(t)^{-1/2}x$ and
the definition of $\mathcal M^{\mathcal{N}}(t)$, we obtain
\begin{align}\label{Sec 7 Lem 3 f-3}
 q_{D_{\varepsilon}}^0(J^{\mathcal{N}}_{\varepsilon}x,J^{\mathcal{N}}_{\varepsilon}x)
 =\varepsilon\sum_{\{i,j\}\in\mathcal E}
   c^{\mathcal{N}}_{ij}(t)(y_i-y_j)^2
 =\varepsilon x^T\mathcal M^{\mathcal{N}}(t)x .
\end{align}
where
\begin{align*}
    q_D^{0}(u,v)
=
\int_D \nabla u\cdot\nabla v.
\qquad u,v\in H^1(D).
\end{align*}
\textbf{$L^2$ norm of $J^{\mathcal{N}}_{\varepsilon}x$}. Fix an edge
$e=\{i,j\}$ and consider the half-tube $T_{ei,\varepsilon}$ adjacent
to $R_{i,\varepsilon}$. By the interpolation formula,
\[
 |J^{\mathcal{N}}_{\varepsilon}x-y_i|
 =\frac{\ell_{\varepsilon}-|s|}{2\ell_{\varepsilon}}\,|y_j-y_i|
 \le C|x|\frac{\ell_{\varepsilon}-|s|}{\ell_{\varepsilon}},
\]
since $y_i=x_i/\sqrt{v_i(t)}$ and the $v_i(t)$ are uniformly bounded
away from zero on $\mathcal B$. Moreover, using
$\rho_{\varepsilon}=d/\sinh\ell_{\varepsilon}$, we obtain
\[
 |F'_{e,\varepsilon}(s+\mathrm i\phi)|=\rho_{\varepsilon}\sqrt{\sinh^2s+\cos^2\phi}
 \le C e^{-(\ell_{\varepsilon}-|s|)}.
\]
Here
$dXdY=|F'_{e,\varepsilon}(s+\mathrm i\phi)|^2\,\mathrm ds\,\mathrm d\phi$.
Since the bounds above depend only on $|s|$, integration over either
half-tube gives
\[
 \begin{aligned}
 \|J^{\mathcal{N}}_{\varepsilon}x-y_i\|_{L^2(T_{ei,\varepsilon})}^2
 &\le \frac{C|x|^2}{\ell_{\varepsilon}^2}
   \int_0^{\ell_{\varepsilon}}
   (\ell_{\varepsilon}-s)^2 e^{-2(\ell_{\varepsilon}-s)}
   \,\mathrm ds \\
 &=\frac{C|x|^2}{\ell_{\varepsilon}^2}
   \int_0^{\ell_{\varepsilon}}s^2e^{-2s}\,\mathrm ds
 \le C\varepsilon^2|x|^2 .
 \end{aligned}
\]
Here  $\ell_{\varepsilon}^{-1}
=\varepsilon c^{\mathcal{N}}_{ij}(t)/\theta=O(\varepsilon)$ uniformly
for $t\in\mathcal B$. Since $J^{\mathcal{N}}_{\varepsilon}x=y_i$ on $R_{i,\varepsilon}$, summing over
the half-tubes attached to this room yields
\[
 \|J^{\mathcal{N}}_{\varepsilon}x-y_i\|_{L^2(S_{i,\varepsilon})}
 \le C\varepsilon |x|.
\]
Also, $|J^{\mathcal{N}}_{\varepsilon}x|\le C|x|$ by the interpolation formula,
and $|S_{i,\varepsilon}|$ is uniformly bounded. Thus, by the
Cauchy--Schwarz inequality,
\[
 \begin{aligned}
\left|\int_{S_{i,\varepsilon}}|J^{\mathcal{N}}_{\varepsilon}x|^2\,\mathrm dz
       -|S_{i,\varepsilon}|y_i^2\right|
 \le \|J^{\mathcal{N}}_{\varepsilon}x-y_i\|_{L^2(S_{i,\varepsilon})}
       \|J^{\mathcal{N}}_{\varepsilon}x+y_i\|_{L^2(S_{i,\varepsilon})} 
 \le C\varepsilon |x|^2.
 \end{aligned}
\]
Summing over $i$ and using the uniform convergence
$|S_{i,\varepsilon}|\to v_i(t)$, we conclude that
 \begin{align}\label{Sec 7 Lem-3 f-1}
\|J^{\mathcal{N}}_{\varepsilon}x\|_{L^2(D_{\varepsilon})}^2
 &=\sum_{i=1}^N |S_{i,\varepsilon}|y_i^2
   +O(\varepsilon)|x|^2 \notag\\
 &=\sum_{i=1}^N v_i(t)y_i^2+o(1)|x|^2
 =|x|^2+o(1)|x|^2 ,
 \end{align}
uniformly for $t\in\mathcal B$.

\textbf{Uniform Poincar\'e inequality} % In the preamble: \usepackage{tikz,capt-of}
% Body fragment for the manuscript. The figure is drawn below in place.
Each $S_{i,\varepsilon}$ satisfies a uniform Neumann Poincar\'e
inequality. By reflection, it suffices to describe the argument on
a half-tube $T_{ei,\varepsilon}$ lying in $X>0$.
Fix $0<X_0<d\cos\theta$ and write
\[
 A_\varepsilon
 =\left\{(X,Y):0<X<X_0,\ |Y|<\omega_{\rho_\varepsilon}(X):=\tan\theta\sqrt{X^2+\rho_\varepsilon^2\cos^2\theta} \right\}.
\]
Since $\omega_{\rho_\varepsilon}(X)\simeq X+\rho_\varepsilon$ and
$|\omega_{\rho_\varepsilon}'(X)|\le\tan\theta$, the map
\begin{align}\label{Sec 7 Lem-3 equ-d-2}
 (X,Y)\longmapsto
 \left(X,\frac{X+\rho_\varepsilon}{\omega_{\rho_\varepsilon}(X)}Y\right) 
\end{align}
is uniformly bi-Lipschitz on $A_\varepsilon$. Indeed, its vertical
stretch factor is bounded above and below, and on $A_\varepsilon$,
\[
 \left|Y\frac{\mathrm d}{\mathrm dX}
          \frac{X+\rho_\varepsilon}{\omega_{\rho_\varepsilon}(X)}\right|
 \le 1+\frac{X+\rho_\varepsilon}{\omega_{\rho_\varepsilon}(X)}
             |\omega_{\rho_\varepsilon}'(X)|\le C.
\]
Its image is the convex trapezoid
$K_{\varepsilon}:=\{(X,Y):0<X<X_0,\ |Y|<X+\rho_\varepsilon\}$, whose diameter is
uniformly bounded, see Figure~\ref{fig:positive-half-tube}(b). The Payne--Weinberger inequality on this trapezoid $K_{\varepsilon}$ (see \cite{PW1960}), for any $f\in H^1(K_{\varepsilon})$, we have
\begin{align*}
    \int_{K_{\varepsilon}}|f-\bar{f}_{K_{\varepsilon}}|^2\leq C\int_{K_{\varepsilon}}|\nabla f|^2,\qquad\mathrm{where}\qquad \bar f_A=|A|^{-1}\int_A f
\end{align*}
Set
$
g(X,Y)=f\left(X,\frac{(X+\rho_\varepsilon)Y}{\omega_{\rho_\varepsilon}(X)}\right).
$
Then the Jacobian of the map \eqref{Sec 7 Lem-3 equ-d-2} is uniformly bounded from above and below, namely $C_1\leq \frac{X+\rho_\varepsilon}{\omega_{\rho_\varepsilon}(X)}\leq C_2$, and hence we obtain
\begin{align}\label{Sec 7 equ-f-2}
    \int_{A_\varepsilon}|g-\bar g_{A_\varepsilon}|^2=&\inf_{c\in \R}\int_{A_\varepsilon}|g-c|^2\leq \int_{A_\varepsilon}|g-\bar{f}_{K_{\varepsilon}}|^2\nonumber\\
    =&\int_{K_\varepsilon}\frac{\omega_{\rho_\varepsilon}(X)}{X+\rho_\varepsilon} |f-\bar{f}_{K_{\varepsilon}}|^2\leq C_1^{-1} \int_{K_{\varepsilon}}|f-\bar{f}_{K_{\varepsilon}}|^2\nonumber\\
    \leq& C\int_{K_{\varepsilon}}|\nabla f|^2\leq  C\int_{A_{\varepsilon}}|\nabla g|^2
\end{align}
therefore gives a uniform Neumann Poincar\'e inequality on $A_\varepsilon$.

Let $B_\varepsilon$ consist of $R_{i,\varepsilon}$ together with the
portion of each incident half-tube beyond $X=X_0/2$, orienting $X$
towards the room on each half-tube,
$B_\varepsilon$ is uniformly bi-Lipschitz equivalent to a fixed
connected Lipschitz domain. It thus satisfies a uniform Neumann
Poincar\'e inequality:
\begin{align}\label{Sec 7 Lem-3 equ-s}
     \int_{B_\varepsilon}|g-\bar g_{B_\varepsilon}|^2\leq C\int_{B_{\varepsilon}}|\nabla g|^2.
\end{align}
 On the half-tube under consideration, the overlap
$O_\varepsilon=A_\varepsilon\cap B_\varepsilon$ has area
\[
 |O_\varepsilon|
 =2\int_{X_0/2}^{X_0}\omega_{\rho_\varepsilon}(X)\,\mathrm dX
 \ge cX_0^2.
\]
For $g\in H^1(S_{i,\varepsilon})$, write
$\bar g_A=|A|^{-1}\int_A g$ for its integral average on $A$.
% In the preamble: \usepackage{tikz,capt-of}
% Paste this block between paragraphs in the manuscript.
% In the preamble: \usepackage{tikz,capt-of}
% Paste this block between paragraphs in the manuscript.
\par\smallskip
\noindent\begin{minipage}{\linewidth}
\centering
\definecolor{halfTubeInk}{HTML}{223342}
\definecolor{halfTubeSoft}{HTML}{536877}
\definecolor{halfTubeBlue}{HTML}{28648A}
\definecolor{halfTubePblue}{HTML}{DDEEF7}
\definecolor{halfTubeGreen}{HTML}{236E66}
\definecolor{halfTubePgreen}{HTML}{D8EEE8}
\definecolor{halfTubeOrange}{HTML}{955116}
\definecolor{halfTubePorange}{HTML}{F8D8A6}
\resizebox{.75\linewidth}{!}{%
\begin{tikzpicture}[font=\normalsize,text=halfTubeInk,
  edge/.style={draw=halfTubeInk,line width=.9pt,line join=round},
  cut/.style={draw=halfTubeSoft,densely dashed,line width=.7pt},
  guide/.style={draw=halfTubeBlue,densely dashed,line width=.7pt},
  span/.style={|-|,line width=.9pt},
  title/.style={anchor=west,font=\bfseries\normalsize,inner sep=0pt},
  subtitle/.style={font=\small,inner sep=1.5pt},
  declare function={g(\x)=0.55*sqrt(\x*\x+0.28*0.28/(1+0.55*0.55));}]
% Example parameters used solely to draw the diagram:
% X_0=2.4, rho=0.28, tan(theta)=0.55.
\def\xend{3.45}
\pgfmathsetmacro{\dend}{\xend*sqrt(1+0.55*0.55)}
\pgfmathsetmacro{\ey}{g(\xend)}
\pgfmathsetmacro{\ang}{atan(0.55)}
\def\halfTubeSection{plot[domain=0:2.4,samples=60,variable=\z](\z,{g(\z)})
  -- plot[domain=2.4:0,samples=60,variable=\z](\z,{-g(\z)}) -- cycle}
\node[title] at (0,10.2) {(a) The half-tube near $X=0$};
\node[title] at (9.0,10.2) {(b) The comparison trapezoid};
\begin{scope}[shift={(4.906,7.58)},x=-1.69cm,y=.79cm]
  \fill[halfTubePblue] \halfTubeSection;
  \draw[edge] \halfTubeSection;
  \draw[guide] (1.4,{-g(1.4)}) -- (1.4,{g(1.4)});
  \fill[halfTubeBlue] (1.4,{g(1.4)}) circle[radius=1.7pt];
  \node[above=4pt,text=halfTubeBlue] at (1.4,{g(1.4)}) {$Y=\omega_{\rho_\varepsilon}(X)$};
  \node at (1.98,0) {$A_\varepsilon$};
  \node[below=5pt] at (0,{-g(0)}) {$0$};
  \node[below=5pt] at (2.4,{-g(2.4)}) {$X_0$};
\end{scope}
\draw[->,>=stealth,line width=1pt,color=halfTubeBlue] (5.65,7.58) -- (8.85,7.58)
 node[midway,below=4pt,subtitle] {Fixed $X$};
\begin{scope}[shift={(13.696,7.58)},x=-1.69cm,y=.79cm]
  \path[edge,fill=halfTubePblue] (0,-.28) -- (2.4,-2.68) -- (2.4,2.68) -- (0,.28) -- cycle;
  \draw[guide] (1.4,-1.68) -- (1.4,1.68);
  \fill[halfTubeBlue] (1.4,1.68) circle[radius=1.7pt];
  \node[above right=3pt,text=halfTubeBlue] at (1.4,1.68) {$Y=X+\rho_\varepsilon$};
  \node[below=5pt] at (0,-.28) {$0$};
  \node[left=3pt] at (2.4,-2.68) {$X_0$};
\end{scope}
\draw[halfTubeSoft!35,thin] (0,5.01) -- (14.4,5.01);
\node[title] at (0,4.56) {(c) The regions $A_\varepsilon$ and $B_\varepsilon$};
\begin{scope}[shift={(13.7116,2.30)},x=-1.84cm,y=.72cm]
  % The room is schematic. Only one incident half-tube is shown.
  \begin{scope}
    \clip (0,0) -- plot[domain=0:\xend,samples=80,variable=\z](\z,{g(\z)})
      .. controls (4.02,2.21) and (4.35,2.53) .. (5.45,2.53)
      .. controls (6.47,2.53) and (6.99,1.70) .. (6.99,0)
      .. controls (6.99,-1.70) and (6.47,-2.53) .. (5.45,-2.53)
      .. controls (4.35,-2.53) and (4.02,-2.21) .. (\xend,{-\ey})
      -- plot[domain=\xend:0,samples=80,variable=\z](\z,{-g(\z)}) -- cycle;
    \fill[halfTubePgreen] (0,-2.8) rectangle (7.05,2.8);
    \fill[halfTubePblue] (0,-2.8) rectangle (1.2,2.8);
    \fill[halfTubePorange] (1.2,-2.8) rectangle (2.4,2.8);
  \end{scope}
  \draw[edge] plot[domain=0:\xend,samples=80,variable=\z](\z,{g(\z)})
      .. controls (4.02,2.21) and (4.35,2.53) .. (5.45,2.53)
      .. controls (6.47,2.53) and (6.99,1.70) .. (6.99,0)
      .. controls (6.99,-1.70) and (6.47,-2.53) .. (5.45,-2.53)
      .. controls (4.35,-2.53) and (4.02,-2.21) .. (\xend,{-\ey})
      -- plot[domain=\xend:0,samples=80,variable=\z](\z,{-g(\z)}) -- cycle;
  \foreach \x in {1.2,2.4}{\draw[cut] (\x,{-g(\x)}) -- (\x,{g(\x)});}
  \node[text=halfTubeOrange] at (1.8,0) {$O_\varepsilon$};
  \node[below=5pt] at (0,{-g(0)}) {$0$};
  \node[below=5pt] at (1.2,{-g(1.2)}) {$X_0/2$};
  \node[below=5pt] at (2.4,{-g(2.4)}) {$X_0$};
  \draw[span,halfTubeBlue] (0,-3.11) -- (2.4,-3.11)
    node[midway,below=3pt,text=halfTubeBlue] {$A_\varepsilon$};
  \draw[span,halfTubeGreen] (1.2,-4.03) -- (6.99,-4.03)
    node[midway,below=3pt,text=halfTubeGreen] {$B_\varepsilon$};
\end{scope}
\end{tikzpicture}%
}\par
\captionof{figure}{}
 \label{fig:positive-half-tube}
\end{minipage}
\par\smallskip

\par\smallskip
\noindent Comparing the two means on the overlap in Figure~\ref{fig:positive-half-tube}(c), using \eqref{Sec 7 equ-f-2} and \eqref{Sec 7 Lem-3 equ-s}, we obtain
\[
 |O_\varepsilon|\,|\bar g_{A_\varepsilon}-\bar g_{B_\varepsilon}|^2
 \le 2\int_{O_\varepsilon}|g-\bar g_{A_\varepsilon}|^2
     +2\int_{O_\varepsilon}|g-\bar g_{B_\varepsilon}|^2
 \le C\int_{A_\varepsilon\cup B_\varepsilon}|\nabla g|^2.
\]
Comparing the mean on each incident half-tube portion $A_\varepsilon$
with the same mean $\bar g_{B_\varepsilon}$ and summing gives
\begin{align}\label{Sec 7.2 equ-c}
 \int_{S_{i,\varepsilon}}|g-\bar g_{S_{i,\varepsilon}}|^2
 \le \int_{S_{i,\varepsilon}}|g-\bar g_{B_\varepsilon}|^2
 \le C\int_{S_{i,\varepsilon}}|\nabla g|^2,
\end{align}
uniformly for $t\in\mathcal B$ and all sufficiently small $\varepsilon$.

Since the sets $S_{i,\varepsilon}$ partition $D_\varepsilon$ up to a
set of measure zero, summing their Poincar\'e inequalities gives
\begin{align}\label{Sec 7 Lem-3 f-2}
    \sum_{i=1}^N
    \|g-\bar g_i\|_{L^2(S_{i,\varepsilon})}^2
    \le C q_{D_\varepsilon}^0(g,g),
    \qquad g\in H^1(D_\varepsilon),\quad \bar g_i:=\bar g_{S_{i,\varepsilon}},
\end{align}
with $C$ independent of $\varepsilon$ and $t$.
Let $V\subset H^1(D_\varepsilon)$ be any subspace of dimension $N+1$.
The linear map
\[
    V\longrightarrow\R^N,\qquad
    g\longmapsto(\bar g_1,\ldots,\bar g_N)
\]
has a nontrivial kernel $\{g\in V:\bar g_1=\cdots=\bar g_N=0\}\ni g_V\ne0$. It satisfies $\overline{(g_V)}_i=0 $ for $i=1,2,\cdots,N$, then 
$\|g_{V} \|_{L^2(D_\varepsilon)}^2=\sum_{i=1}^N
    \|g_{V}-\overline{(g_V)}_i\|_{L^2(S_{i,\varepsilon})}^2 $. Hence the min--max principle and \eqref{Sec 7 Lem-3 f-2} give
\[
    \lambda_{N+1}^0(D_\varepsilon)
    =\inf_{\dim V=N+1}\;
      \sup_{g\in V\setminus\{0\}}
      \frac{q_{D_\varepsilon}^0(g,g)}
           {\|g\|_{L^2(D_\varepsilon)}^2} 
    \ge C^{-1}=:\gamma^{\mathcal{N}}.
\]

For the upper bound, the norm estimate (see \eqref{Sec 7 Lem-3 f-1}) for $J^{\mathcal{N}}_\varepsilon$ gives
$\|J^{\mathcal{N}}_\varepsilon x\|_{L^2(D_\varepsilon)}^2\ge \tfrac12|x|^2$
for all sufficiently small $\varepsilon$, uniformly in $t$.
Thus $J^{\mathcal{N}}_\varepsilon\R^N$ is an $N$-dimensional trial space.
Using \eqref{Sec 7 Lem 3 f-3} and the min--max
principle, we obtain
\[
    \lambda_N^0(D_\varepsilon)
    \le \sup_{x\in\R^N\setminus\{0\}}
       \frac{\varepsilon x^T\mathcal M^{\mathcal{N}}(t)x}
            {\|J^{\mathcal{N}}_\varepsilon x\|_{L^2(D_\varepsilon)}^2}
    \le 2\varepsilon\sup_{t\in\mathcal B}\|\mathcal M^{\mathcal{N}}(t)\|.
\]
Consequently, the first $N$ Neumann eigenvalues converge to zero
uniformly in $t$, while $\lambda_{N+1}^0(D_\varepsilon)\ge C^{-1}$,
as in \cite[\S6(d), p.~611]{CdV}.

On each tube $T_{e,\varepsilon}$, by \eqref{Sec 7 f-12}, the pullback of $J^{\mathcal{N}}_\varepsilon x$
by $F_{e,\varepsilon}$ is linear in $s$, with derivative
$(y_j-y_i)/(2\ell_\varepsilon)$.
By the conformal invariance, $J^{\mathcal{N}}_\varepsilon x$ is harmonic on $T_{e,\varepsilon}$
and has zero normal derivative on the sides $\phi=\pm\theta$.
At the two interfaces $s=\pm\ell_\varepsilon$, we have
\[
    |F_{e,\varepsilon}'(\pm\ell_\varepsilon+\mathrm i\phi)|^2
    =d^2+\rho_\varepsilon^2\cos^2\phi\ge d^2.
\]
Thus, with $\nu_\varepsilon$ denoting the outward unit normal to the tube,
\[
    |\partial_{\nu_\varepsilon}J^{\mathcal{N}}_\varepsilon x|
    =\frac{|y_j-y_i|}
           {2\ell_\varepsilon
            |F_{e,\varepsilon}'(\pm\ell_\varepsilon+\mathrm i\phi)|}
    \le \frac{|y_j-y_i|}{2d\ell_\varepsilon}
    \le C\varepsilon |x|,
\]
where the last inequality follows from
$\ell_\varepsilon^{-1}=\varepsilon c^{\mathcal{N}}_{ij}(t)/\theta$
and $|y_j-y_i|\le C|x|$.

Since $J^{\mathcal{N}}_\varepsilon x$ is constant on each room, integration by parts
on the tubes gives, for every $g\in H^1(D_\varepsilon)$,
\[
    q_{D_\varepsilon}^0(J^{\mathcal{N}}_\varepsilon x,g)
    =\sum_{e=\{i,j\}\in\mathcal E}
      \int_{\Gamma_{ei,\varepsilon}\cup\Gamma_{ej,\varepsilon}}
          g\,\partial_{\nu_\varepsilon}J^{\mathcal{N}}_\varepsilon x\,\mathrm dS.
\]
Each interface has uniformly bounded length, and the trace inequality
on its room-side collar gives
\[
    \int_{\Gamma_{ei,\varepsilon}} |g|\,\mathrm dS
    \le |\Gamma_{ei,\varepsilon}|^{1/2}
        \|g\|_{L^2(\Gamma_{ei,\varepsilon})}
    \le C\|g\|_{H^1(D_\varepsilon)}.
\]
Taking summation for all edges, we obtain
\[
    |q_{D_\varepsilon}^0(J^{\mathcal{N}}_\varepsilon x,g)|
    \le C\varepsilon |x|\,\|g\|_{H^1(D_\varepsilon)},
\]
uniformly in $t\in\mathcal B$ for all sufficiently small $\varepsilon$, 
 see also ~\cite[(6.5), p.~607]{CdV}.

Let $P^{\mathcal{N}}_{\varepsilon}=P^{\mathcal{N}}_{\varepsilon,t}$ be the spectral projection onto
the first $N$ Neumann eigenfunction space in $D_{\varepsilon}=D_{\varepsilon,t}$.
For $x\in\R^N$, write
\[
    J^{\mathcal{N}}_{\varepsilon}x=P^{\mathcal{N}}_{\varepsilon}J^{\mathcal{N}}_{\varepsilon}x+z_{\varepsilon},
    \qquad
    z_{\varepsilon}=(I-P^{\mathcal{N}}_{\varepsilon})J^{\mathcal{N}}_{\varepsilon}x.
\]
Since $z_{\varepsilon}\perp\operatorname{ran}P^{\mathcal{N}}_{\varepsilon}$ in $L^2$
and $\lambda_{N+1}^0(D_{\varepsilon})\ge \gamma^{\mathcal{N}}$, we have
\[
    \|z_{\varepsilon}\|_{L^2(D_{\varepsilon})}^2
    \le (\gamma^{\mathcal{N}})^{-1}q_{D_{\varepsilon}}^0(z_{\varepsilon},z_{\varepsilon}).
\]
The Neumann eigenfunction equation also gives
$q_{D_{\varepsilon}}^0(P^{\mathcal{N}}_{\varepsilon}J^{\mathcal{N}}_{\varepsilon}x,z_{\varepsilon})=0$.
Hence, by the mixed estimate with $g=z_{\varepsilon}$,
\[
    \begin{aligned}
    \|z_{\varepsilon}\|_{H^1(D_{\varepsilon})}^2
    &\le (1+(\gamma^{\mathcal{N}})^{-1})
        q_{D_{\varepsilon}}^0(z_{\varepsilon},z_{\varepsilon})\\
    &= (1+(\gamma^{\mathcal{N}})^{-1})
        q_{D_{\varepsilon}}^0(J^{\mathcal{N}}_{\varepsilon}x,z_{\varepsilon})
     \le C\varepsilon |x|\,\|z_{\varepsilon}\|_{H^1(D_{\varepsilon})}.
    \end{aligned}
\]
Thus $\|z_{\varepsilon}\|_{H^1(D_{\varepsilon})}\le C\varepsilon|x|$.
By orthogonality in $L^2$ and with respect to $q_{D_{\varepsilon}}^0$,
\begin{align*}  \|P^{\mathcal{N}}_{\varepsilon}J^{\mathcal{N}}_{\varepsilon}x\|_{L^2(D_{\varepsilon})}^2
    =\|J^{\mathcal{N}}_{\varepsilon}x\|_{L^2(D_{\varepsilon})}^2
      -\|z_{\varepsilon}\|_{L^2(D_{\varepsilon})}^2
      =|x|^2+o(1)|x|^2,
\end{align*}
and
    \begin{align}\label{Sec 7 Lem 3 f-2}
q_{D_{\varepsilon}}^0(P^{\mathcal{N}}_{\varepsilon}J^{\mathcal{N}}_{\varepsilon}x,P^{\mathcal{N}}_{\varepsilon}J^{\mathcal{N}}_{\varepsilon}x)
    =&\varepsilon x^T\mathcal M^{\mathcal{N}}(t)x
      -q_{D_{\varepsilon}}^0(z_{\varepsilon},z_{\varepsilon})\nonumber\\
      =&\varepsilon x^T\mathcal M^{\mathcal{N}}(t)x+O(\varepsilon^2)|x|^2,
    \end{align}
uniformly for $t\in\mathcal B$.

As in Lemma~\ref{lem:rooms}
(\eqref{Sec 7 Lem equ-a}, \eqref{Sec 7 Lem equ-b},
\eqref{Sec 7 Lem equ-d}), put
\[
    G^{\mathcal{N}}_{\varepsilon}=(P^{\mathcal{N}}_{\varepsilon}J^{\mathcal{N}}_{\varepsilon})^*P^{\mathcal{N}}_{\varepsilon}J^{\mathcal{N}}_{\varepsilon},
    \qquad
    U^{\mathcal{N}}_{\varepsilon,t}=P^{\mathcal{N}}_{\varepsilon}J^{\mathcal{N}}_{\varepsilon}(G^{\mathcal{N}}_{\varepsilon})^{-1/2}.
\]
The first estimate gives $G^{\mathcal{N}}_{\varepsilon}=I+o(1)$ in operator norm.
Thus, for sufficiently small $\varepsilon$, $G^{\mathcal{N}}_{\varepsilon}$ is
invertible and $(U^{\mathcal{N}}_{\varepsilon,t})^*U^{\mathcal{N}}_{\varepsilon,t}=I$.
Since $\dim\operatorname{ran}P^{\mathcal{N}}_{\varepsilon}=N$,
$U^{\mathcal{N}}_{\varepsilon,t}$ is an isometry from $\R^N$ onto
$\operatorname{ran}P^{\mathcal{N}}_{\varepsilon}$. By the estimate \eqref{Sec 7 Lem 3 f-2}, let $H^0_{D_{\varepsilon}}=-\Delta^{\mathcal{N}}_{D_{\varepsilon}}$
\[
    (U^{\mathcal{N}}_{\varepsilon,t})^*\varepsilon^{-1}H_{D_{\varepsilon}}^0
    U^{\mathcal{N}}_{\varepsilon,t}
    =(G^{\mathcal{N}}_{\varepsilon})^{-1/2}
      \bigl(\mathcal M^{\mathcal{N}}(t)+O(\varepsilon)\bigr)(G^{\mathcal{N}}_{\varepsilon})^{-1/2}.
\]
Since $(G^{\mathcal{N}}_{\varepsilon})^{-1/2}=I+o(1)$ uniformly on $\mathcal B$,
we obtain
\begin{equation}\label{eq:neumann-low-matrix}
    B^{\mathcal{N}}_{\varepsilon,t}
    :=(U^{\mathcal{N}}_{\varepsilon,t})^*H_{D_{\varepsilon}}^0
       U^{\mathcal{N}}_{\varepsilon,t},
    \qquad
    \sup_{t\in\mathcal B}\|\varepsilon^{-1}B^{\mathcal{N}}_{\varepsilon,t}-\mathcal M^{\mathcal{N}}(t)\|
    \longrightarrow0
    \qquad\varepsilon \to 0.
\end{equation}
Fix $\varepsilon>0$ and a parameter $t_0$. On each tube use the
identification
\[
 F_{e,\varepsilon,t_0}(s+i\phi)\ \longmapsto\
 F_{e,\varepsilon,t}\bigl((\ell_{\varepsilon,t}/
                 \ell_{\varepsilon,t_0})s+i\phi\bigr),
\]
where the subscript $t$ makes the parameter dependence explicit.
Continue the same formula a short distance into each room collar.
For $t$ near $t_0$ these maps are smoothly close to the identity
there. Choose the collar joins by cutting off these displacements
inside the collars, and use the prescribed graph deformations on the
disjoint room boundary charts.
This gives diffeomorphisms from $D_{\varepsilon,t_0}$ to
$D_{\varepsilon,t}$ that depend smoothly on $t$.
The pulled-back mass and energy forms are $C^1$ in form norm on a
fixed $H^1$ space and their sum is locally uniformly coercive.
Inverting the shifted form and using the Riesz formula at the
spectral gap gives a $C^1$ low spectral projection.
The pulled-back trial functions are $C^1$ as well, since the
normalized tube coordinate $s/\ell_{\varepsilon,t}$ is fixed.
The Gram normalization therefore makes $B^{\mathcal N}_{\varepsilon,t}$
and its rescaling $C^1$ in $t$.

Set $\widetilde T^{\mathcal{N}}_\varepsilon(t)=(\varepsilon^{-1}B^{\mathcal{N}}_{\varepsilon,t}+I)^{-1}$.
Since $\varepsilon^{-1}B^{\mathcal{N}}_{\varepsilon,t}$ and $\mathcal M^{\mathcal{N}}(t)$ are nonnegative, the resolvent
identity gives
\[
    \|\widetilde T^{\mathcal{N}}_\varepsilon(t)-T^{\mathcal{N}}(t)\|
    \le \|\varepsilon^{-1}B^{\mathcal{N}}_{\varepsilon,t}-\mathcal M^{\mathcal{N}}(t)\|.
\]
Choose $r_1<r_0$, where $r_0$ is given by  Lemma~\ref{lem:stability}.
Put $M_0=\sup_{t\in\mathcal B}\|\mathcal M^{\mathcal{N}}(t)\|$. The largest eigenvalue of
$\varepsilon^{-1}B^{\mathcal{N}}_{\varepsilon,t}$ is
$\varepsilon^{-1}\lambda_N^0(D_{\varepsilon,t})$, whereas
$\varepsilon^{-1}\lambda_{N+1}^0(D_{\varepsilon,t})\ge \gamma^{\mathcal{N}}/\varepsilon$.
Thus, by \eqref{eq:neumann-low-matrix}, we may fix $\varepsilon>0$
so small that $ \sup_{t\in\mathcal B}\|\widetilde T^{\mathcal{N}}_\varepsilon(t)-T^{\mathcal{N}}(t)\|<\frac{\delta_1(r_1)}{2}$, and 
\[
    \sup_{t\in\mathcal B}
        \varepsilon^{-1}\lambda_N^0(D_{\varepsilon,t})<M_0+1,
    \qquad
    \inf_{t\in\mathcal B}
        \varepsilon^{-1}\lambda_{N+1}^0(D_{\varepsilon,t})>M_0+2.
\]
\smallskip
For this fixed $\varepsilon$, the family
$\widetilde T^{\mathcal{N}}_\varepsilon(t)$ is in fact $C^1$ in $t$.
Apply Lemma~\ref{lem:stability} on $\mathcal H=\R^N$ with
\[
    T=T^{\mathcal{N}}(0),\qquad T(t)=T^{\mathcal{N}}(t),\qquad
    \widetilde T(t)=\widetilde T^{\mathcal{N}}_\varepsilon(t),
    \qquad \mu=(a+1)^{-1}.
\]
The compressed differential at $0$ is $s\mapsto-\mu^2s$ on
$\Sym(E)$, and the uniform approximation error is less than
$\delta_1(r_1)$. Thus there exists $t_*\in\mathcal B$ such that
\[
    \dim\ker\bigl(\widetilde T^{\mathcal{N}}_\varepsilon(t_*)-\mu I\bigr)=m,
    \qquad
    \#\{j:\lambda_j(\widetilde T^{\mathcal{N}}_\varepsilon(t_*))>\mu\}=p .
\]
The eigenvalues of $\widetilde T^{\mathcal{N}}_\varepsilon(t_*)$ are
$(1+\varepsilon^{-1}\lambda_j^0(D_{\varepsilon,t_*}))^{-1}$,
$1\le j\le N$. Since $\lambda\mapsto(1+\varepsilon^{-1}\lambda)^{-1}$
is strictly decreasing on $[0,\infty)$ and takes the value $\mu$ at
$\lambda=\varepsilon a$, exactly $m$ of these first $N$ Neumann
eigenvalues equal $\varepsilon a$ and exactly $p$ lie below it.
Moreover,
\[
    \varepsilon^{-1}\lambda_{N+1}^0(D_{\varepsilon,t_*})
    >M_0+2>a,
\]
since $a\le\|\mathcal M^{\mathcal{N}}(0)\|\le M_0$. Hence
\[
    \dim\ker\bigl(H^0_{D_{\varepsilon,t_*}}-\varepsilon a I\bigr)=\dim E,
    \qquad
    \#\{j:\lambda_j^0(D_{\varepsilon,t_*})<\varepsilon a\}
    =n_-(\mathcal M).
\]
This proves the proposition for $\sigma=0$ with
$\Omega=D_{\varepsilon,t_*}$.

\smallskip
\noindent\textbf{A fixed Robin parameter.}
Keep $\varepsilon>0$ as chosen in the Neumann case. For
$D_{\varepsilon}=D_{\varepsilon,t}$, the Robin and Neumann forms
have the common domain $H^1(D_{\varepsilon})$ and satisfy
\[
    q_{D_{\varepsilon}}^b(u,u)
    =q_{D_{\varepsilon}}^0(u,u)
     +b\int_{\partial D_{\varepsilon}}|u|^2\,\mathrm ds.
\]
Since the domains form a compact smooth family for fixed $\varepsilon$,
the trace inequality gives
\[
    \frac{|b|}{\varepsilon}
    \int_{\partial D_{\varepsilon}}|u|^2\,\mathrm ds
    \le C_{\varepsilon}|b|
       \left(\varepsilon^{-1}q_{D_{\varepsilon}}^0(u,u)
             +\|u\|_{L^2(D_{\varepsilon})}^2\right)
\]
uniformly for $t\in\mathcal B$.  Let \[
H_{D_\varepsilon}^{b}u=-\Delta u,\qquad
\operatorname{Dom}(H_{D_\varepsilon}^{b})
=\bigl\{u\in H^2(D_\varepsilon):
\partial_\nu u+bu=0\ \text{on }\partial D_\varepsilon\bigr\}.
\]
Thus, in the sense of quadratic forms,
\[
    (1-C_{\varepsilon}|b|)
       (I+\varepsilon^{-1}H^0_{D_{\varepsilon}})
    \le I+\varepsilon^{-1}H^b_{D_{\varepsilon}}
    \le (1+C_{\varepsilon}|b|)
       (I+\varepsilon^{-1}H^0_{D_{\varepsilon}}).
\]
For $C_{\varepsilon}|b|<1$, the middle operator is strictly positive.
Put
\[
    S_{\varepsilon,t}^b
    =(I+\varepsilon^{-1}H^b_{D_{\varepsilon,t}})^{-1}.
\]
Taking inverses in the preceding inequalities and using
$\|S_{\varepsilon,t}^0\|\le1$, we obtain
\[
    \sup_{t\in\mathcal B}
    \|S_{\varepsilon,t}^b-S_{\varepsilon,t}^0\|
    \le \frac{C_{\varepsilon}|b|}{1-C_{\varepsilon}|b|}
    \longrightarrow0
    \qquad(b\to0).
\]

By the min--max principle and the gap obtained in the Neumann case,
for all sufficiently small $|b|$ we have
\begin{align}\label{Sec 7 f-13}
    \sup_{t\in\mathcal B}
        \varepsilon^{-1}\lambda_N^b(D_{\varepsilon,t})
    <M_0+\tfrac32
    <\inf_{t\in\mathcal B}
        \varepsilon^{-1}\lambda_{N+1}^b(D_{\varepsilon,t}).
\end{align}
Let $P_{\varepsilon,t}^b$ be the spectral projection onto the first
$N$ Robin modes, so that $P_{\varepsilon,t}^0=P_{\varepsilon,t}^{\mathcal{N}}$.
The uniform gap and the Riesz projection formula give
\[
    \sup_{t\in\mathcal B}
    \|P_{\varepsilon,t}^b-P_{\varepsilon,t}^{\mathcal{N}}\|
    \longrightarrow0
    \qquad(b\to0).
\]
With $U=U_{\varepsilon,t}^{\mathcal{N}}$, define
\[
    U_{\varepsilon,t}^b
    =P_{\varepsilon,t}^bU
      \bigl(U^*P_{\varepsilon,t}^bU\bigr)^{-1/2}.
\]
Since $U^*P_{\varepsilon,t}^bU\to I$ uniformly in $t$, this is an
isometry from $\R^N$ onto $\operatorname{ran}P_{\varepsilon,t}^b$,
and $U_{\varepsilon,t}^b\to U_{\varepsilon,t}^{\mathcal{N}}$ uniformly
as $b\to0$. As in \eqref{eq:neumann-low-matrix}, set
\[
    B_{\varepsilon,t}^b
    =(U_{\varepsilon,t}^b)^*H^b_{D_{\varepsilon,t}}U_{\varepsilon,t}^b,
    \qquad
    \widetilde T_\varepsilon^b(t)
    =(\varepsilon^{-1}B_{\varepsilon,t}^b+I)^{-1}
    =(U_{\varepsilon,t}^b)^*S_{\varepsilon,t}^bU_{\varepsilon,t}^b.
\]
It follows that
\begin{equation}\label{eq:robin-matrix-continuity}
    \sup_{t\in\mathcal B}
    \|\widetilde T_\varepsilon^b(t)
      -\widetilde T_\varepsilon^{\mathcal{N}}(t)\|
    \longrightarrow0
    \qquad(b\to0).
\end{equation}
Under the same domain identifications as in the Neumann case, the
boundary measure pulls back to a smooth density on the fixed boundary.
Hence, for fixed $\varepsilon$ and $b$, the pulled-back Robin forms
depend $C^1$ on $t$ in form norm. The Riesz projection formula and
the Gram matrix normalization above then give $C^1$ dependence of
$\widetilde T_\varepsilon^b(t)$ on $t$.

Now fix the prescribed $\sigma\in\R$. For $v(x)=u(x/r)$ on $rD$,
\[
    (\partial_\nu v+\sigma v)(rx)
    =r^{-1}(\partial_\nu u+r\sigma u)(x),
    \qquad x\in\partial D.
\]
Thus the Robin parameters on $D$ and $rD$ are $r\sigma$ and $\sigma$,
respectively, and
\begin{equation}\label{eq:robin-scaling}
    \lambda_j^\sigma(rD)=r^{-2}\lambda_j^{r\sigma}(D).
\end{equation}

Choose $r>0$ so that $b=r\sigma$ satisfies the positivity and gap
estimates above and makes the error in
\eqref{eq:robin-matrix-continuity} less than $\delta_1(r_1)/2$.
For $\sigma=0$, take $r=1$. By the choice of $\varepsilon$ in the
Neumann case,
\[
    \begin{aligned}
    \sup_{t\in\mathcal B}
    \|\widetilde T_\varepsilon^b(t)-T^{\mathcal{N}}(t)\|
    &\le
    \sup_{t\in\mathcal B}
    \|\widetilde T_\varepsilon^b(t)
      -\widetilde T_\varepsilon^{\mathcal{N}}(t)\|
    +\sup_{t\in\mathcal B}
    \|\widetilde T_\varepsilon^{\mathcal{N}}(t)-T^{\mathcal{N}}(t)\|\\
    &<\delta_1(r_1).
    \end{aligned}
\]
Apply Lemma~\ref{lem:stability} on $\R^N$ to the operators
$T^{\mathcal{N}}(t)$ and  approximating family
$\widetilde T_\varepsilon^b(t)$. There exists $t_b\in\mathcal B$
such that
\[
    \dim\ker\bigl(\widetilde T_\varepsilon^b(t_b)-\mu I\bigr)=m,
    \qquad
    \#\{j:\lambda_j(\widetilde T_\varepsilon^b(t_b))>\mu\}=p.
\]
The eigenvalues of $\widetilde T_\varepsilon^b(t_b)$ are
$(1+\varepsilon^{-1}\lambda_j^b(D_{\varepsilon,t_b}))^{-1}$,
$1\le j\le N$. Their denominators are positive. By
$\mu=(1+a)^{-1}$, $a=-\lambda_{\min}(\mathcal{M})$ and \eqref{Sec 7 f-13}, then
\[
    \varepsilon^{-1}\lambda_{N+1}^b(D_{\varepsilon,t_b})
    >M_0+\tfrac32>a.
\]
Finally, set $\Omega=rD_{\varepsilon,t_b}$. By
\eqref{eq:robin-scaling}, $H_\Omega^\sigma$ has the eigenvalue
$\varepsilon a/r^2$ of multiplicity $m=\dim E$, with exactly
$p=n_-(\mathcal M)$ eigenvalues below it, counted with multiplicity.
\end{proof}
Finally, we discuss the case $k=2$. For the Dirichlet and Neumann cases, the sharpness of $m_2(\Omega)=3$ is given by \eqref{eq:maink=20}. The possibility of multiplicity three for the second Robin eigenvalue was discussed in \cite[\S6.1.5]{BH}, but a precise construction remains open. In the following, we give an affirmative answer.

\begin{cor}\label{cor:second-robin-sharp}
For every fixed $(\alpha,\beta)\in\mathbb R^2\setminus\{(0,0)\}$,
\begin{align}\label{Sec 7 Cor f-1}
    \max_{\Omega}
    \dim\ker\bigl(H^{(\alpha,\beta) }_\Omega-\lambda_2^{(\alpha,\beta) }(\Omega)I\bigr)=3,
\end{align}
where $\Omega$ ranges over all bounded connected planar domains with Lipschitz boundary.
\end{cor}

\begin{proof}
The upper bound has been proved at the end of Section \ref{sec:upper}. Therefore, we only need to construct a domain $\Omega$ satisfying $\dim\ker\bigl(H^{(\alpha,\beta) }_\Omega-\lambda_2^{(\alpha,\beta) }(\Omega)I\bigr)=3 $.

Take $G$ to be the complete graph $K_4$
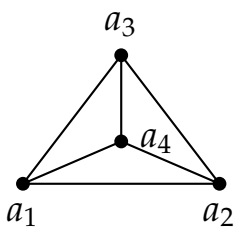
\begin{tikzpicture}[baseline={(0,.25)}, font=\scriptsize]
  \coordinate (a1) at (-.65,0);
  \coordinate (a2) at ( .65,0);
  \coordinate (a3) at (0,.85);
  \coordinate (a4) at (0,.28);

  \draw (a1) -- (a2) -- (a3) -- cycle;
  \draw (a1) -- (a4) -- (a2)  (a4) -- (a3);

  \foreach \v/\pos in {1/below,2/below,3/above,4/right}
    \fill (a\v) circle (1.3pt) node[\pos] {$a_{\v}$};
\end{tikzpicture},
and set $\mathcal M=-\mathbf 1\mathbf 1^T$, where
$\mathbf 1=(1,1,1,1)^T$. The graph $K_4$ is planar, and
$\mathcal M$ is well signed. Since $\mathcal M\mathbf 1=-4\mathbf 1$
and $\mathcal M$ vanishes on $\mathbf 1^\perp$, we have
\[
    n_-(\mathcal M)=1,\qquad
    E=\ker\mathcal M=\mathbf 1^\perp,\qquad \dim E=3.
\]
Moreover, $\mathcal L(K_4)=\Sym(\R^4)$, so every symmetric bilinear
form on $E$ extends to an element of $\mathcal L(K_4)$ by taking
its value to be zero whenever one argument lies in $E^\perp$.
Thus the restriction map $R_E$ is surjective.

Proposition~\ref{net:robin-realization} and Proposition \ref{net:realization} now give a bounded connected
smooth planar domain $\Omega$ for which
\[
    \lambda_1^{(\alpha,\beta) }(\Omega)
    <\lambda_2^{(\alpha,\beta) }(\Omega)
    =\lambda_3^{(\alpha,\beta) }(\Omega)
    =\lambda_4^{(\alpha,\beta) }(\Omega)
    <\lambda_5^{(\alpha,\beta) }(\Omega).
\]
This proves the equality.
\end{proof}

\begin{proof}[\bf Proof of the sharpness in Theorem~\ref{thm:main}]
Fix any $(\alpha,\beta)\in\R^2\setminus\{(0,0)\}$ and $k\ge3$.
Lemma~\ref{net:graph} gives a connected planar graph and a
well-signed matrix with $R_E(\mathcal L(G))=\Sym(E)$ and
\[
 n_-(\mathcal{M})=k-1,\qquad \dim\ker \mathcal{M}=k,
\]
If $\alpha=0$, Proposition~\ref{net:realization} gives a disk with
finitely many circular holes and Dirichlet boundary conditions.
If $\alpha\ne0$, put $\sigma=\beta/\alpha$. Then
$H_\Omega^{(\alpha,\beta)}=H_\Omega^\sigma$, and
Proposition~\ref{net:robin-realization} can construct a planar  domain consisting of rooms and tubes, see \ref{fig:sec72-neumann}.
In either case the eigenvalue has multiplicity $k$ and exactly $k-1$
eigenvalues below it. Hence, for domain $\Omega_k$ with $k\geq 3$,
\[
 \lambda_{k-1}(\Omega_k)<\lambda_k(\Omega_k)
 =\cdots=\lambda_{2k-1}(\Omega_k)<\lambda_{2k}(\Omega_k).
\]
Thus $m_k(\Omega_k)=k$. For $k=2$, see  Corollary \ref{cor:second-robin-sharp}.
\end{proof}
% Add the following entry to the bibliography if {CdV} is not already present.
% \bibitem{CdV}
% Y. Colin de Verdi\`ere,
% \emph{Construction de laplaciens dont une partie finie du spectre
% est donn\'ee},
% Ann. Sci. \'Ecole Norm. Sup. (4) \textbf{20} (1987), 599--615.

\subsection*{Acknowledgments}  Z. Chen is supported by National Key R\&D Program of China (No. 2023YFA1010002) and NSFC (No. 12222109). 
S. Zhang is supported by the Postdoctoral Fellowship Program and China Postdoctoral Science Foundation under Grant Numbers BX20250062 and 2026M793381.

\subsection*{AI use declaration}
Inspired by our recent work \cite{CLZ}, the authors conjectured the
optimal multiplicity bound and used ChatGPT to explore possible
proof strategies. The proof of sharpness builds on our previous
work \cite{CLZ} and the original approach of Colin de Verdi\`ere
\cite{CdV}, with suitable modifications. ChatGPT was also used
to improve the language and presentation of the manuscript.
The authors assume full responsibility for the mathematical
validity, accuracy, and integrity of the proofs presented in
the manuscript.

\end{document}